\documentclass[a4,10pt]{article}
\usepackage{amsmath,amscd,amssymb}
\usepackage[margin=.8in]{geometry}

\newcommand{\nbiga}{\mathcal{A}}

\newcommand{\nbigd}{\mathcal{D}}
\newcommand{\nbige}{\mathcal{E}}

\newcommand{\nbigh}{\mathcal{H}}

\newcommand{\nbigl}{\mathcal{L}}

\newcommand{\nbigo}{\mathcal{O}}
\newcommand{\nbigp}{\mathcal{P}}

\newcommand{\nbigr}{\mathcal{R}}
\newcommand{\nbigs}{\mathcal{S}}

\newcommand{\nbigu}{\mathcal{U}}
\newcommand{\nbigv}{\mathcal{V}}

\newcommand{\nbigx}{\mathcal{X}}
\newcommand{\nbigy}{\mathcal{Y}}
\newcommand{\nbigz}{\mathcal{Z}}

\newcommand{\proj}{\mathbb{P}}
\newcommand{\seisuu}{{\mathbb Z}}

\newcommand{\cnum}{{\mathbb C}}
\newcommand{\real}{{\mathbb R}}

\newcommand{\gbigr}{\mathfrak R}

\newcommand{\vecv}{{\boldsymbol v}}

\newcommand{\lrarr}{\longrightarrow}

\newcommand{\pf}{{\bf Proof}\hspace{.1in}}
\newcommand{\qed}{\mbox{\rule{1.2mm}{3mm}}}

\def\End{\mathop{\rm End}\nolimits}

\def\Re{\mathop{\rm Re}\nolimits}

\def\Gr{\mathop{\rm Gr}\nolimits}

\def\SL{\mathop{\rm SL}\nolimits}

\def\rank{\mathop{\rm rank}\nolimits}

\def\Gr{\mathop{\rm Gr}\nolimits}

\def\tr{\mathop{\rm tr}\nolimits}
\def\Tr{\mathop{\rm Tr}\nolimits}
\def\vol{\mathop{\rm dvol}\nolimits}

\def\id{\mathop{\rm id}\nolimits}

\newcommand{\del}{\partial}
\newcommand{\delbar}{\overline{\del}}

\newcommand{\nhom}{{\mathcal Hom}}

\newcommand{\barz}{\overline{z}}
\newcommand{\zbar}{\barz}

\newcommand{\lefttop}[1]{{}^{#1}\!}

\def\Harm{\mathop{\rm Harm}\nolimits}

\newcommand{\Etilde}{\widetilde{E}}

\newcommand{\thetatilde}{\widetilde{\theta}}

\newcommand{\nablatilde}{\widetilde{\nabla}}

\newcommand{\wbar}{\overline{w}}

\newcommand{\htilde}{\widetilde{h}}

\newcommand{\pitilde}{\widetilde{\pi}}

\newcommand{\stilde}{\widetilde{s}}

\newcommand{\Dtilde}{\widetilde{D}}
\newcommand{\Xtilde}{\widetilde{X}}

\def\Gal{\mathop{\rm Gal}\nolimits}

\newcommand{\gtilde}{\widetilde{g}}

\newcommand{\nbigxtilde}{\widetilde{\nbigx}}

\newcommand{\nbigdtilde}{\widetilde{\nbigd}}

\newcommand{\betabar}{\overline{\beta}}

\newcommand{\Kbar}{\overline{K}}

\newcommand{\veckappa}{{\boldsymbol \kappa}}

\newcommand{\Xbar}{\overline{X}}

\newcommand{\Sigmatilde}{\widetilde{\Sigma}}

\newtheorem{thm}{Theorem}[section]
\newtheorem{cor}[thm]{Corollary}
\newtheorem{rem}[thm]{Remark}
\newtheorem{lem}[thm]{Lemma}
\newtheorem{prop}[thm]{Proposition}
\newtheorem{df}[thm]{Definition}
\newtheorem{example}[thm]{Example}
\newtheorem{condition}[thm]{Condition}

\begin{document}

\title{Asymptotic behaviour of large-scale solutions
of Hitchin's equations in higher rank}

\author{Takuro Mochizuki\thanks{Research Institute for Mathematical Sciences, Kyoto University, Kyoto 606-8512, Japan, takuro@kurims.kyoto-u.ac.jp}
\and
Szil\'{a}rd Szab\'{o}\thanks{
Department of Geometry, Budapest University of Technology and 
Economics, M\H{u}egyetem rakpart 3, 1111, Budapest, Hungary,
szabosz@math.bme.hu}
\thanks{
Alfr\'ed R\'enyi Institute of Mathematics, Re\'altanoda utca 
13-15, 1053, Budapest, Hungary,
szabo.szilard@renyi.hu
}}
\date{}
\maketitle

\begin{abstract}

Let $X$ be a compact Riemann surface.
Let $(E,\theta)$ be a stable Higgs bundle
of degree $0$ on $X$.
Let $h_{\det(E)}$ denote a flat metric of
the determinant bundle $\det(E)$.
For any $t>0$, there exists a unique harmonic metric
$h_t$ of $(E,\theta)$
such that $\det(h_t)=h_{\det(E)}$.
We prove that
if the Higgs bundle is induced by
a line bundle on the normalization of the spectral curve,
then the sequence $h_t$ is convergent to the naturally defined
decoupled harmonic metric at the speed of the exponential order.
We also obtain a uniform convergence
for such a family of Higgs bundles.

\vspace{.1in}
\end{abstract}

\section{Introduction}

\subsection{Background}

Let $X$ be a Riemann surface.
Let $(E,\delbar_E,\theta)$ be a Higgs bundle of rank $r$ on $X$.
Let $h$ be a Hermitian metric of $E$.
We obtain the Chern connection $\nabla_h$ of $(E,\delbar_E,h)$
and the adjoint $\theta^{\dagger}_h$ of $\theta$
with respect to $h$.
Let $R(h)$ denote the curvature of $\nabla_h$.
The metric $h$ is called a harmonic metric of
$(E,\delbar_E,\theta)$
if
\[
 R(h)+[\theta,\theta^{\dagger}_h]=0.
\]
The metric $h$ is called a decoupled harmonic metric of
$(E,\delbar_E,\theta)$
if
\[
 R(h)=[\theta,\theta^{\dagger}_h]=0.
\]

Suppose that $X$ is compact,
and that $(E,\delbar_E,\theta)$ is stable of degree $0$.
Let $\Sigma_{E,\theta}$ denote the spectral curve of $(E,\theta)$.
We assume that $(E,\delbar_E,\theta)$ is generically regular semisimple,
i.e.,
$D(E,\theta)=\bigl\{
P\in X\,\big|\,|T^{\ast}_PX\cap\Sigma_{E,\theta}|<r
\bigr\}$
is a finite subset of $X$.

Let $h_{\det(E)}$ be a flat metric of $\det(E)$.
According to Hitchin \cite{Hitchin-self-duality}
and Simpson \cite{s1},
$(E,\delbar_E,\theta)$ has a unique harmonic metric $h$
such that $\det(h)=h_{\det(E)}$.
Because $(E,\delbar_E,t\theta)$ is stable of degree $0$
for any $t>0$,
there exists a unique harmonic metric $h_t$
of $(E,\delbar_E,t\theta)$ for any $t>0$
such that $\det(h_t)=h_{\det(E)}$.
We are interested in the behaviour of $h_t$ as $t\to\infty$.
See \cite{Gaiotto-Moore-Neitzke},
\cite{Katzarkov-Noll-Pandit-Simpson},
and \cite{Mazzeo-Swoboda-Weiss-Witt}
for the motivation of the study.
It is related with the geometric P=W conjecture
\cite{Szabo-2021, Szabo-2022}.
See also helpful survey papers \cite{Li-survey,Swoboda-survey}.

For any simply connected relatively compact open subset $K$
of $X\setminus D(E,\theta)$,
there exists a decomposition of the Higgs bundle
\begin{equation}
\label{eq;22.12.2.20}
 (E,\delbar_E,\theta)_{|K}
  =\bigoplus_{i=1}^{r}
  (E_{K,i},\delbar_{E_{K,i}},\theta_{K,i})
\end{equation}
such that $\rank E_{K,i}=1$.
According to \cite{Decouple},
there exist $C(K)>0,\epsilon(K)>0$
such that
\[
 \bigl|h_t(u,v)\bigr|
 \leq
 C(K)\exp(-\epsilon(K)t)|u|_{h_t}|v|_{h_t}
\]
for any local sections $u$ and $v$
of $E_{K,i}$ and $E_{K,j}$ $(i\neq j)$
in the decomposition (\ref{eq;22.12.2.20}).
It implies that there exist
$C'(K)>0$ and $\epsilon'(K)>0$ such that
\[
 \bigl|
 R(h_t)_{|K}
 \bigr|_{h_t}
 =\bigl|
 [\theta,\theta^{\dagger}_{h_t}]
 \bigr|_{h_t}
\leq C'(K)\exp(-\epsilon'(K)t).
\]
As a result, for any sequence $t(i)\to\infty$,
there exist a subsequence $t'(j)\to \infty$
and gauge transformations $g_{t'(j)}$
such that
the sequence $g_{t'(j)}^{\ast}h_{t'(j)}$ 
is convergent to a decoupled harmonic metric
of $(E,\delbar_E,\theta)_{|X\setminus D(E,\theta)}$
in the $C^{\infty}$-sense
locally on $X\setminus D(E,\theta)$.

We may ask the following questions
under appropriate assumptions.
\begin{description}
 \item[Q1]  Is there a sequence of gauge transformations
	    $g_t$
	    such that $g_t^{\ast}h_t$
	    convergent as $t\to\infty$ locally on $X\setminus D(E,\theta)$?
	    In other words, is the limit independent of the choice of
	    a subsequence?
 \item[Q2]  Let $K\subset X\setminus D(E,\theta)$
	   be any relatively compact open subset.
	   Then, is the order of the convergence on $K$
	    dominated by $e^{-\delta(K) t}$
	    for some $\delta(K)>0$?
\end{description}

In the rank two case,
under the assumption that $\Sigma_{E,\theta}$ is smooth,
Mazzeo-Swoboda-Weiss-Witt \cite{Mazzeo-Swoboda-Weiss-Witt}
solved the both questions completely.
In \cite{Decouple}, the question {\bf Q1} was solved
without assuming the smoothness of the spectral curve.
In the higher rank case,
Collier-Li \cite{Collier-Li}
solved the both questions for cyclic Higgs bundles.
Fredrickson \cite{Fredrickson-Generic-Ends} studied the both questions
when the spectral curve is smooth,
under a mild assumption on the ramification
of the spectral curve over $X$
(see Remark \ref{rem;22.12.27.2} and
\cite[Proposition 2.2, (2.9)]{Fredrickson-Generic-Ends}).

\begin{rem}
Chronologically,
the study {\rm\cite{Decouple}} was done 
inspired by the previous researches
{\rm\cite{Collier-Li}},
{\rm\cite{Katzarkov-Noll-Pandit-Simpson}}
and
{\rm\cite{Mazzeo-Swoboda-Weiss-Witt}}.
\hfill\qed
\end{rem}

\begin{rem}
\label{rem;22.12.27.2}
Let $Q\in\Sigma_{E,\theta}$ be a critical point of
$\pi:\Sigma_{E,\theta}\to X$.
Put $P=\pi(Q)$.
Let $(X_P,z)$ be a coordinate neighbourhood around $P$.
By using the holomorphic $1$-form $dz$,
we obtain the trivialization
$T^{\ast}X_P\simeq \cnum\times X_P$.
Let $\Sigma_{E,\theta,Q}$ denote the connected component of
$T^{\ast}X_P\cap \Sigma_{E,\theta}$ which contains $Q$.
We may assume that
$\Sigma_{E,\theta,Q}\cap T_P^{\ast}X_P=\{Q\}$
and that $\Sigma_{E,\theta,Q}$
is holomorphically isomorphic to a disc.
Let $r(Q)$ denote the degree of
$\Sigma_{E,\theta,Q}\to X_P$.
There exist holomorphic functions $a_j$
$(j=0,\ldots,r(Q)-1)$ on $X_P$
such that
\[
 \Sigma_{E,\theta,Q}
 =\left\{
 (y,z)\in \cnum\times X_P\,\left|\,
 y^{r(Q)}+\sum_{j=0}^{r(Q)-1}
 a_j(z)y^{j}=0
 \right.
 \right\}.
\]
Because $T^{\ast}X_P\cap \Sigma_{E,\theta,Q}=\{Q\}$,
there exists $\alpha\in\cnum$ such that 
\begin{equation}
\label{eq;22.12.18.1}
 y^{r(Q)}+\sum_{j=0}^{r(Q)-1}
 a_j(0)y^j 
=(y-\alpha)^{r(Q)}.
\end{equation}
The smoothness of $\Sigma_{E,\theta,Q}$ is equivalent to
the condition that
$a_0(z)-(-\alpha)^{r(Q)}$
has simple $0$ at $z=0$.
To study the local property of 
$\Sigma_{E,\theta,Q}$ around $Q$
and $\theta$ around $P$,
we may assume that $\alpha=0$
by considering 
$\theta_{|X_P}-\alpha dz\,\cdot \id_{E_{|X_P}}$.
Moreover,
we may assume that $a_{r(Q)-1}$ is constantly $0$
by considering
$\theta_{|X_P}-r(Q)^{-1}a_{r(Q)-1}\,dz\cdot \id_{E_{|X_P}}$.
By changing the coordinate $z$ to $w(z)$
satisfying $w(0)=0$
and $w(\del_zw)^{r(Q)}=-a_0(z)$,
we may assume that $a_0(z)=-z$.
In general,
$a_j$ $(1\leq j\leq r(Q)-2)$ are not constantly $0$.
\hfill\qed
\end{rem}

\subsection{Main results}

\subsubsection{The symmetric case}
\label{subsection;22.12.26.10}

As a first main result,
let us mention that
if $(E,\delbar_E,\theta)$
has a non-degenerate symmetric pairing $C$,
then both questions {\bf Q1} and {\bf Q2} are extremely easy.
As explained in \cite{Li-Mochizuki3},
there exists a unique decoupled harmonic metric $h^{C}$
of $(E,\theta)_{|X\setminus D(E,\theta)}$
which is compatible with $C$.
By using a variant of Simpson's main estimate
and an elementary linear algebraic argument
in \S\ref{subsection;22.12.25.30},
we can solve both questions {\bf Q1} and {\bf Q2},
and the limit is $h^{C}$ in this case.
The following theorem is a special case of
Corollary \ref{cor;22.12.20.31}.
\begin{thm}
\label{thm;22.12.27.1}
Let $K$ be any relatively compact open subset of $X\setminus D(E,\theta)$.
Let $s(h^C,h_t)$ denote the automorphism of $E_{|X\setminus D(E,\theta)}$
determined by $h_t=h^C\cdot s(h^C,h_t)$. 
For any $\ell\in\seisuu_{\geq 0}$,
there exist positive constants $C(\ell,K)$ and $\epsilon(\ell,K)$
such that
the $L_{\ell}^2$-norm of
$s(h^C,h_t)-\id$ on $K$
are dominated by 
$C(\ell,K)\exp(-\epsilon(\ell,K)t)$
as $t\to\infty$.
\end{thm}

For example,
we may apply this theorem to
a Higgs bundle contained in the Hitchin section
because it has a canonical non-degenerate symmetric pairing.

Indeed,
in Theorem \ref{thm;22.12.27.1},
we do not need to assume that $X$ is compact.
See Theorem \ref{thm;22.12.3.1}
and Corollary \ref{cor;22.12.20.31}
for the precise statements.
These results are also technically useful,
which will be applied to the third main result
(see \S\ref{subsection;22.12.26.1}--\ref{subsection;22.12.27.3}).

\subsubsection{The irreducible case}
\label{subsection;22.12.3.3}

The second main result in this paper is
an affirmative answer to {\bf Q1}
in the case that the spectral curve is
locally and globally irreducible.

\begin{thm}[Corollary
\ref{cor;22.12.24.30}]
Suppose that $\Sigma_{E,\theta}$ is locally irreducible
(see {\rm\cite[Page 8]{Grauert-Remmert}}).
Then, the sequence $h_t$ is convergent to
a decoupled harmonic metric $h_{\infty}$
in the $C^{\infty}$-sense
locally on $X\setminus D(E,\theta)$.
 \end{thm}
See Theorem \ref{thm;22.11.30.130} for the more general statement.
Note that $\Sigma_{E,\theta}$ is connected
because of the stability condition of $(E,\theta)$.

More precisely,
we canonically construct
a filtered bundle $\nbigp^{\star}_{\ast}(\nbigv)$
over $\nbigv=E(\ast D(E,\theta))$
in an algebraic way from $(E,\theta)$
such that
(i)  $(\nbigp^{\star}_{\ast}(\nbigv),\theta)$
is a decomposable filtered Higgs bundle
in the sense of Definition \ref{df;22.12.25.20},
(ii) $(\nbigp^{\star}_{\ast}(\nbigv),\theta)$
is stable of degree $0$,
(iii)
$\det(\nbigp^{\star}_{\ast}\nbigv)$
equals the filtered bundle naturally induced by $\det(E)$.
There exists a unique decoupled harmonic metric $h_{\infty}$
of $(E,\theta)_{|X\setminus D(E,\theta)}$
adapted to $\nbigp^{\star}_{\ast}(\nbigv)$
such that $\det(h_{\infty})=h_{\det(E)}$.
We shall prove that
the sequence $h_t$ is convergent to $h_{\infty}$
as $t\to\infty$ on $X\setminus D(E,\theta)$.

An outline of the proof is as follows.
Let $P\in D(E,\theta)$.
Let $X_P$ be a small neighbourhood of $P$ in $X$.
By a theorem of Donaldson \cite{Donaldson-boundary-value},
there exists a harmonic metric
$h_{P,t}$ of $(E,\delbar,t\theta)_{|X_P}$
such that
$h_{P,t|\del X_P}=h_{\infty|\del X_P}$.
According to Proposition \ref{prop;22.11.30.131},
the sequence $h_{P,t}$ is convergent to $h_{\infty|X_P\setminus\{P\}}$
in the $C^{\infty}$-sense
locally on $X_P\setminus\{P\}$ as $t\to\infty$.
As in \cite{Mazzeo-Swoboda-Weiss-Witt},
by patching $h_{P,t}$ and $h_{\infty}$,
we construct a family of Hermitian metrics $\htilde_t$ $(t>0)$
of $E$
such that
(i) $\det(\htilde_t)=h_{\det(E)}$,
(ii) $\lim_{t\to\infty}\htilde_t=h_{\infty}$ on $X\setminus D(E,\theta)$,
(iii)
$\int_X
\bigl|
R(\htilde_t)
+\bigl[
 t\theta,(t\theta)^{\dagger}_{\htilde_t}
 \bigr]\bigr|\to 0$.
Let $s(\htilde_t,h_t)$ denote the automorphism of
$E$ determined by $h_t=\htilde_t\cdot s(\htilde_t,h_t)$.
Then, we shall prove that
$\sup_X(s(\htilde_t,h_t)-\id_E)\to 0$
by the essentially same argument as that in \cite{Decouple}.

Because of the assumption of the local irreducibility
of $\Sigma_{E,\theta}$,
it is easy to find the candidate of
``the limiting configuration'' $h_{\infty}$.
In the rank two case,
the Higgs bundle $(E,\theta)_{|X_P}$ is easy to understand.
There is a homogeneous wild harmonic bundle
$(E'_P,\theta'_P,h'_P)$ on $(\proj^1,\infty)$
such that the restriction of $(E'_P,\theta'_P)$
to a neighbourhood of $0$ is isomorphic to $(E,\theta)_{|X_P}$,
where we consider an $S^1$-action on $\proj^1$
induced by $(a,z)\mapsto a^mz$ for some $m\in\seisuu_{>0}$.
(See \cite[\S8]{Mochizuki-KH-Higgs} for homogeneity of
harmonic bundles with respect to an $S^1$-action.)
The special case is a fiducial solution
in \cite{Mazzeo-Swoboda-Weiss-Witt}.
In \cite{Decouple},
the restriction of $h'_P$ was useful
in the construction of approximate solutions $\htilde_t$.
In the higher rank case,
the Higgs bundle $(E,\theta)_{|X_P}$ is more complicated
even under the assumption of the local irreducibility.
It does not seem that
the approximation by a homogeneous wild harmonic bundle
can work well.
Therefore, we develop a way to use the solutions of
the boundary-value problem in the construction of
approximate solutions.

\begin{rem}
Because we also study the question ${\bf Q1}$
for wild harmonic bundles
under a similar assumption on the spectral curve, 
we also study the Dirichlet problem  
for wild harmonic bundles
(Theorem {\rm\ref{thm;22.12.20.4}}).
\hfill\qed
\end{rem}

\subsubsection{The order of convergence in the smooth case}
\label{subsection;22.12.26.1}

We study question {\bf Q2}
under the following additional condition.
\begin{condition}
\label{condition;22.12.25.40}
Let $\rho:\Sigmatilde_{E,\theta}\to \Sigma_{E,\theta}$
be the normalization.
There exists a holomorphic line bundle $L$
with an isomorphism
$E\simeq (\pi\circ\rho)_{\ast}L$ 
such that
$\theta$ is induced by
the $\nbigo_{T^{\ast}X}$-action on $\rho_{\ast}L$.
\hfill\qed
\end{condition}

For example,
this condition is satisfied
if $\Sigma_{E,\theta}$ is smooth
according to \cite{Beauville-Narasimhan-Ramanan,Hitchin-self-duality}.
We shall prove the following theorem.

\begin{thm}[Theorem
\ref{thm;22.12.24.101}]
\label{thm;22.12.3.2}
Suppose that Condition {\rm\ref{condition;22.12.25.40}}
is satisfied.
Let $s(h_{\infty},h_{t})$  be the automorphism of
$(E,\theta)_{|X\setminus D(E,\theta)}$
determined by $h_t=h_{\infty}\cdot s(h_{\infty},h_t)$. 
Let $K\subset X\setminus D(E,\theta)$
be any relatively compact open subset.
For any $\ell\in\seisuu_{\geq 0}$,
there exist $C(\ell,K)>0$ and $\epsilon(\ell,K)>0$
such that the following holds as $t\to\infty$:
\[
 \bigl\|
(s(h_{\infty},h_t)-\id)_{|K}
 \bigr\|_{L_{\ell}^2}
\leq C(\ell,K)
\exp\bigl(-\epsilon(\ell,K)t\bigr).
\] 
 \end{thm}

To prove Theorem \ref{thm;22.12.3.2},
we refine the construction of $\htilde_t$
in \S\ref{subsection;22.12.3.3}.
For each $P\in D(E,\theta)$,
there exists a non-degenerate symmetric pairing $C_P$
of $(E,\delbar_E,\theta)_{|X_P}$
such that
$C_{P|X_P\setminus\{P\}}$ is compatible
with $h_{\infty|X_P\setminus\{P\}}$.
It is easy to see that
the harmonic metric $h_{P,t}$ of $(E,\theta)_{|X_P}$
satisfying $h_{P,t|\del X_P}=h_{\infty|\del X_P}$
is compatible with $C_{P}$.
Let $s(h_{\infty},h_{P,t})$ be the automorphism of
$E_{|X_P\setminus\{P\}}$
determined by
$h_{P,t}=h_{\infty|X_P\setminus\{P\}}\cdot s(h_{\infty},h_{P,t})$.
By the result in the symmetric case mentioned
in \S\ref{subsection;22.12.26.10},
on any relatively compact open subset $K$
of $X_P\setminus\{P\}$,
$s(h_{\infty},h_{P,t})-\id$
converges to $0$
at the speed of the order $e^{-\delta(K)t}$.
Then,
the following stronger condition is satisfied:
\[
 \int_X\bigl|
R(\htilde_t)
+\bigl[
 t\theta,(t\theta)^{\dagger}_{\htilde_t}
 \bigr]
 \bigr|_{\htilde_t}
 \leq Ce^{-\delta t}.
\]
Then, we can obtain the estimate
of $\sup|s(\htilde_t,h_t)-\id|$ on
any relatively compact open subset in
$X\setminus D(E,\theta)$.
By a general argument in \S\ref{subsection;22.12.3.10},
we can obtain the desired estimate
of the norms of
$s(\htilde_t,h_t)-\id$
and its higher derivatives on $X$
even around $D(E,\theta)$.

\subsubsection{A family case}
\label{subsection;22.12.27.3}

The result and the method in \S\ref{subsection;22.12.26.1}
can be generalized to the following family case.
Let $p_1:\nbigxtilde\to\nbigs$ be a smooth proper morphism
of complex manifolds such that each fiber is connected
and $1$-dimensional.
We also assume that $\nbigs$ is connected.
Let $\pi:\nbigs\times T^{\ast}X\to\nbigs\times X$
and $p_2:\nbigs\times X\to \nbigs$ denote the projections.
Let $\Phi_0:\nbigxtilde\to \nbigs\times T^{\ast}X$
be a morphism of complex manifolds
such that $p_2\circ\pi\circ\Phi_0=p_1$.
We set $\Phi_1:=\pi\circ\Phi_0:\nbigxtilde\to\nbigs\times X$.
We assume the following.
\begin{itemize}
 \item $\Phi_1$ is proper and finite.
 \item There exits a closed complex analytic hypersurface
       $\nbigd\subset \nbigs\times X$
       such that
       (i) $\nbigd$ is finite over $\nbigs$,
       (ii) the induced map
       $\nbigxtilde\setminus\Phi_1^{-1}(\nbigd)
       \lrarr (\nbigs\times X)\setminus\nbigd$
       is a covering map,
       (iii)
       $\Phi_0$ induces an injection
       $\nbigxtilde\setminus\Phi_1^{-1}(\nbigd)
       \lrarr
       \nbigs\times T^{\ast}X$.
\end{itemize}
We set $r:=|\Phi_1^{-1}(P)|$
for any $P\in(\nbigs\times X)\setminus\nbigd$.
Let $g(X)$ and $\gtilde$ denote the genus of
$X$ and $p_1^{-1}(x)$ $(x\in \nbigs)$, respectively.
We set $X_x=\{x\}\times X$
and $\nbigd_{x}=\nbigd\cap X_x$.
There exists a natural isomorphism $X_x\simeq X$.
We note that $\nbigd\to \nbigs$ is not assumed to be
a covering map,
and hence $|\nbigd_{x}|$ is not necessarily
constant on $\nbigs$.

Let $\nbigl$ be a holomorphic line bundle
on $\nbigxtilde$ such that
$\deg(\nbigl_{|p_1^{-1}(x)})=
\gtilde-rg(X)+r-1$.
We obtain a locally free $\nbigo_{\nbigs\times X}$-module
$\nbige=\Phi_{1\ast}(\nbigl)$.
It is equipped with the morphism
$\theta:
\nbige\to \nbige\otimes \Omega^1_{\nbigs\times X/\nbigs}$
induced by
the $\nbigo_{\nbigs\times T^{\ast}X}$-action
on $\Phi_{0\ast}\nbigl$.
For each $x\in\nbigs$,
we obtain the Higgs bundle
$(\nbige_x,\theta_x)=(\nbige,\theta)_{|X_x}$,
which is stable of degree $0$.

There exists a Hermitian metric $h_{\det\nbige}$ of $\det(\nbige)$
such that $h_{\det\nbige|X_x}$ are flat
for any $x\in\nbigs$.
There exist harmonic metrics
$h_{t,x}$ of $(\nbige_x,t\theta_x)$ $(x\in\nbigs)$
such that
$\det(h_{t,x})=h_{\det\nbige|X_x}$.
There also exist decoupled harmonic metrics
$h_{\infty,x}$ $(x\in\nbigs)$
of $(\nbige_x,\theta_x)_{|X_x\setminus \nbigd_{x}}$
such that
$\det(h_{\infty,x})=h_{\det(\nbige)|X_x\setminus\nbigd_x}$.

\begin{thm}[Theorem
\ref{thm;22.12.25.4}]
\label{thm;23.1.27.1}
Let $x_0\in\nbigs$.
Let $K$ be a relatively compact open subset
of $X_{x_0}\setminus\nbigd_{x_0}$.
Let $\nbigs_0$ be a neighbourhood of $x_0$
in $\nbigs$ such that
$\nbigs_0\times K$ is relatively compact in
$(\nbigs\times X)\setminus\nbigd$.
For any $\ell\in\seisuu_{\geq 0}$,
there exist $C(\ell),\epsilon(\ell)>0$
such that
the $L_{\ell}^2$-norm of
$s(h_{\infty,x},h_{t,x})-\id$ $(x\in \nbigs_0,t\geq 1)$
on $K$
are dominated by
$C(\ell)\exp(-\epsilon(\ell)t)$.
\end{thm}

\begin{rem}
Note that
for another Hermitian metric $h'_{\det\nbige}$ of $\det(\nbige)$
such that $h'_{\det\nbige|X_x}$ are flat
for any $x\in\nbigs$,
we obtain an $\real_{>0}$-valued
$C^{\infty}$-function $\beta$ on $\nbigs$
determined by
$h'_{\det(\nbige)}=\beta h_{\det(\nbige)}$,
and 
$\beta^{1/r}h_{t,x}$
(resp. $\beta^{1/r}h_{\infty,x}$)
are harmonic metrics (resp. decoupled harmonic metrics)
of $(\nbige_x,t\theta_x)$
(resp. $(\nbige_x,\theta_x)_{|X_x\setminus\nbigd_x}$)
such that 
$\det(\beta^{1/r}h_{t,x})=h'_{\det\nbige|X_x}$
(resp. 
$\det(\beta^{1/r}h_{\infty,x})=
h'_{\det(\nbige)|X_x\setminus\nbigd_x}$).
Hence, the claim of Theorem {\rm\ref{thm;23.1.27.1}} is
independent of the choice of $h_{\det\nbige}$.
\hfill\qed
\end{rem}

\begin{rem}
We may apply Theorem {\rm\ref{thm;23.1.27.1}}
to obtain a locally uniform estimate
for large scale solutions of the Hitchin equation
for a family of stable Higgs bundles of degree $0$
whose spectral curves are smooth. 
\hfill\qed
\end{rem}

\paragraph{Acknowledgements}

Both authors are grateful to Laura Fredrickson
for the answer to our question related to Remark \ref{rem;22.12.27.2}.

T.~M. thanks Qiongling Li for discussions and collaborations.
In particular, the ideas in \S\ref{section;23.1.28.1}
are variants of those in \cite{Li-Mochizuki2, Li-Mochizuki3}.
T.~M. is grateful to Rafe Mazzeo and Siqi He
for some discussion.
T.~M. is partially supported by
the Grant-in-Aid for Scientific Research (A) (No. 21H04429),
the Grant-in-Aid for Scientific Research (A) (No. 22H00094),
and the Grant-in-Aid for Scientific Research (C) (No. 20K03609),
Japan Society for the Promotion of Science.
This work was partially supported by the Research Institute for Mathematical
Sciences, an International Joint Usage/Research Center located in Kyoto
University

Sz.~Sz. benefited of support from the \emph{Lend\"ulet} Low Dimensional
Topology grant of the Hungarian Academy of Sciences
and of the grants K120697 and KKP126683 of NKFIH.
Sz.~Sz. partially enjoyed
the hospitality of the Max Planck Institute for Mathematics (Bonn).
Sz.~Sz. would like to thank Olivia Dumitrescu
and Richard Wentworth for inspiring discussions.

\section{Preliminaries}

\subsection{Some definitions}

\subsubsection{Decoupled harmonic bundles}

Let $Y$ be a Riemann surface.
Let $(V,\theta)$ be a Higgs bundle on $Y$.

\begin{df}
A Hermitian metric $h$ of $V$
is called a decoupled harmonic metric of $(V,\theta)$
if the following conditions are satisfied.
\begin{description}
 \item[(A1)] $h$ is a harmonic metric of
	     the Higgs bundle $(V,\delbar_V,\theta)$.
 \item[(A2)] $h$ is flat,
       i.e., the Chern connection
       $\nabla_h$ of $(V,\delbar_V,h)$ is flat.
\end{description}
Such $(V,\theta,h)$ is called a 
decoupled harmonic bundle.
\hfill\qed	    
\end{df}
Note that 
the conditions {\bf (A1)} and {\bf (A2)}
imply that 
$\theta$ and $\theta^{\dagger}_h$ are commuting.

\subsubsection{Symmetric Higgs bundles}

Let $C$ be a non-degenerate symmetric product of $V$.
It is called a non-degenerate symmetric product of
the Higgs bundle $(V,\theta)$
if $\theta$ is self-adjoint with respect to $C$.
Such a tuple $(V,\theta,C)$ is called a symmetric Higgs bundle.
Let $V^{\lor}$ denote the dual bundle of $V$.
Let $\Psi_C:V\to V^{\lor}$ be the isomorphism
induced by $C$.
Let $\theta^{\lor}$ be the induced Higgs field of $V^{\lor}$.
The condition is equivalent to that
$\Psi_C$ induces an isomorphism
of the Higgs bundles
$(V,\theta)\simeq (V^{\lor},\theta^{\lor})$.

A Hermitian metric $h$ of $V$ is called compatible with $C$
if $\Psi_C$ is isometric with respect to $h$
and its dual Hermitian metric $h^{\lor}$ of $V^{\lor}$.

\subsubsection{Generically regular semisimple Higgs bundles}

Let $\Sigma_{V,\theta}\subset T^{\ast}Y$
denote the spectral curve of $(V,\theta)$.
We say that $(V,\theta)$ is regular semisimple
if the projection $\Sigma_{V,\theta}\to Y$ is a covering map.
We say that $(V,\theta)$ is generically regular semisimple
if there exists a discrete subset $D\subset Y$
such that
$(V,\theta)_{|Y\setminus D}$ is regular semisimple.

Let $\pi:\Sigma_{V,\theta}\to Y$ denote the projection.
If $(V,\theta)$ is regular semisimple,
there exists a line bundle $L_V$ on $\Sigma_{V,\theta}$
with an isomorphism
$\pi_{\ast}L_V\simeq V$
such that $\theta$ is induced by
$\nbigo_{T^{\ast}Y}$-action on $L_V$.

\subsection{Regular semisimple case}

\subsubsection{Decoupled harmonic metrics}

Suppose that $(V,\theta)$ is regular semisimple.
We consider the following condition for a Hermitian metric $h$ of $V$.
\begin{description}
 \item[(A3)] For any $P\in Y$,
       the eigen decomposition of $\theta$ at $P$
       is orthogonal with respect to $h$.
\end{description} 
Note that {\bf (A3)} holds if and only if
$\theta$ and $\theta^{\dagger}_h$ are commuting.
The following lemma is easy to see.
\begin{lem}
If two of the conditions
{\rm \bf (A1), (A2), (A3)}
are satisfied for a Hermitian metric $h$ of $V$,
then $h$ is a decoupled harmonic metric of $(V,\theta)$.
\hfill\qed
\end{lem}

A flat metric $h_{L_V}$ of $L_V$
induces a Hermitian metric $\pi_{\ast}(h_{L_V})$ of $V$.
It is easy to check that
$\pi_{\ast}(h_{L_V})$ is a decoupled harmonic metric of
$(V,\theta)$.

\begin{prop} 
\label{prop;22.12.23.2}
This procedure induces an equivalence between
flat metrics of $L_V$ 
and a decoupled harmonic metrics of $(V,\theta)$.
\hfill\qed
\end{prop}

\begin{rem}
Let $(V,\theta,h)$ be a decoupled harmonic bundle.
Let $\Sigma_{V,\theta}=\coprod_{i\in\Lambda}\Sigma_{V,\theta,i}$
be the decomposition into connected components.
There exists the corresponding decomposition
of the Higgs bundle
$(V,\theta)=\bigoplus_{i\in\Lambda}(V_i,\theta_i)$
such that
$\Sigma_{V_i,\theta_i}=\Sigma_{V,\theta,i}$.
Because $h$ is a decoupled harmonic metric,
the decomposition is orthogonal with respect to $h$.
Hence, we obtain the decomposition of a decoupled harmonic bundle
$(V,\theta,h)=\bigoplus (V_i,\theta_i,h_i)$.
\hfill\qed
\end{rem}

\subsubsection{Symmetric products}
\label{subsection;22.12.23.100}

The multiplication of $\nbigo_{\Sigma_{V,\theta}}$
induces a multiplication
\[
\pi_{\ast}\nbigo_{\Sigma_{V,\theta}}
\otimes
\pi_{\ast}\nbigo_{\Sigma_{V,\theta}}
\lrarr
\pi_{\ast}\nbigo_{\Sigma_{V,\theta}}.
\]
Any local section $f$
of $\pi_{\ast}\nbigo_{\Sigma_{V,\theta}}$
induces an endomorphism $F_f$ of
the locally free $\nbigo_Y$-module
$\pi_{\ast}\nbigo_{\Sigma_{V,\theta}}$.
We obtain the local section $\tr(f):=\tr(F_f)$
of $\nbigo_Y$.

Let $C_{L_V}$ be a non-degenerate symmetric pairing of $L_V$.
We obtain the non-degenerate pairing $C$ of $V=\pi_{\ast}L_V$:
\begin{equation}
\label{eq;22.12.22.30}
\begin{CD}
 V \otimes_{\nbigo_Y} V
@>{\pi_{\ast}C_{L_V}}>>
 \pi_{\ast}\nbigo_{\Sigma_{V,\theta}}
@>{\tr}>>
 \nbigo_Y.
\end{CD}
\end{equation}

\begin{prop}
\label{prop;22.12.23.1}
This procedure induces an equivalence
between non-degenerate symmetric pairings of $L_V$ 
and non-degenerate symmetric pairings of $(V,\theta)$.
\hfill\qed
\end{prop}

We recall the following proposition.

\begin{prop}[\mbox{\cite[Proposition 2.30]{Li-Mochizuki3}}]
\label{prop;22.12.20.20}
For any non-degenerate symmetric pairing $C$ of $(V,\theta)$,
there exists a unique decoupled harmonic metric $h^C$
of $(V,\theta)$ which is compatible with $C$.
\hfill\qed
\end{prop}
Indeed, let $C_{L_V}$ be the non-degenerate symmetric pairing
of $L_V$ corresponding to $C$.
Let $h_{L_V}$ be the unique Hermitian metric of $L_V$
satisfying $h_{L_V}(s,s)=|C_{L_V}(s,s)|$.
We obtain the Hermitian metric $h^C$ corresponding to $h_{L_V}$.
Then, it is the decoupled harmonic metric compatible with $C$.

As for the converse,
the following holds.
\begin{lem}
\label{lem;22.12.23.20}
Let $h$ be a decoupled harmonic metric of $(V,\theta)$.
There exists a non-degenerate symmetric pairing of $(V,\theta)$
compatible with $h$
if and only if the following condition is satisfied.
\begin{itemize}
 \item Let $h_{L_V}$ be the corresponding Hermitian metric of
       $L_V$, whose Chern connection is flat.
       Let $\Sigma_{V,\theta,i}$ be any connected component of
       $\Sigma_{V,\theta}$.
       Let $\rho_i:\pi_1(\Sigma_{V,\theta,i})\to S^1$
       be the homomorphism
       obtained as the monodromy of
       $(L_V,h_{L_V})_{|\Sigma_{V,\theta,i}}$.       
       Then, the image of $\rho_i$ is contained in
       $\{\pm 1\}$.       
\end{itemize} 
\end{lem}
\pf
There exists a non-degenerate symmetric pairing of $(V,\theta)$
compatible with $h$
if and only if
there exists a non-degenerate symmetric pairing $C_{L_V}$ of $L_V$
compatible with $h_{L_V}$.
If such $C_{L_V}$ exists,
then each $\rho_i$ comes from an $\real$-representation.
(See \cite[\S2]{Li-Mochizuki3}.)
Hence, the image is contained in $\{\pm 1\}$.
Conversely, if the image of each $\rho_i$ is contained in $\{\pm 1\}$,
then it is easy to construct such a pairing $C_{L_V}$.
\hfill\qed

\subsection{Dirichlet problem for wild harmonic bundles on curves}
\label{subsection;22.12.23.50}

Let $Y$ be a Riemann surface
equipped with a K\"ahler metric $g_Y$.
Let $X\subset Y$ be a connected relatively compact connected open subset
whose boundary $\del X$ is smooth and non-empty.
Let $D\subset X$ be a finite subset.

Let $(\nbigp_{\ast}\nbigv,\theta)$
be a good filtered Higgs bundle on $(Y,D)$ of rank $r$.
We obtain
$(\det(\nbigp_{\ast}\nbigv),\tr(\theta))$.
We set
$(V,\theta)=(\nbigv,\theta)_{|Y\setminus D}$.
Let $h_{\del X}$ be a Hermitian metric of $V_{|\del X}$.

\begin{thm}
\label{thm;22.12.20.4}
There exists a unique harmonic metric $h$ of $(V,\delbar_V,\theta)_{|X}$
such that
(i) $h_{|\del X}=h_{\del X}$,
(ii) $\nbigp^h_{\ast}(V)=\nbigp_{\ast}\nbigv$.
\end{thm}
\pf
Let us study the case $r=1$.
There exists a Hermitian metric $h_0$ of $V$
such that
(i) $h_{0|\del X}=h_{\del X}$,
(ii) $h_0$ is flat around any point of $D$,
(iii) $\nbigp^{h_0}_{\ast}(V)=\nbigp_{\ast}V$.
There exists a $C^{\infty}$-function
$\alpha:X\to \real$
such that $\alpha_{|\del X}=0$
and that
$\delbar\del\alpha=R(h_0)$.
Then, $h=e^{-\alpha}h_0$ is a flat metric of $V$
satisfying the desired conditions.
Let $h'$ be another flat metric satisfying the same condition.
We obtain the $C^{\infty}$-function $s$ on $X$
determined by $h'=e^sh$.
Because $\Delta_{g_Y}s=0$ and $s_{|\del X}=0$,
we obtain that $s=0$ on $X$,
and hence $h'=h$.

Let us study the case $r\geq 2$.
At each point $P\in D$,
let $(X_P,z_P)$ be a holomorphic coordinate neighbourhood
around $P$ such that
(i) $X_P$ is relatively compact in $X\setminus (D\setminus\{P\})$,
(ii) $\Xbar_{P}\cap \Xbar_{P'}=\emptyset$ for any $P,P'\in D$,
(iii) the coordinate $z_P$ induces
$(X_P,P)\simeq (\{|z|<1\},0)$.
Let $h_{\det(V)}$ be a flat metric of $\det (V)$
adapted to $\det(\nbigp_{\ast}V)$
such that $h_{\det(V)|\del X}=\det(h_{\del X})$.
Let $h_0$ be a Hermitian metric of $V$
such that
(i) $h_{0|\del X}=h_{\del X}$,
(ii) $\det(h_0)=h_{\det(V)}$,
(iii) $\nbigp^{h_0}_{\ast}(V)=\nbigp_{\ast}\nbigv$,
(iv) around $P\in D$,
we have
$\bigl|
R(h_0)+[\theta,\theta^{\dagger}_{h_0}]
\bigr|_{h,g_Y}
=O(|z_P|^{-2+\epsilon})$
for some $\epsilon>0$.
(For example, see \cite{Mochizuki-KH-Higgs} for the construction.)
We set $F(h_0)=R(h_0)+[\theta,\theta^{\dagger}_{h_0}]$.
There exists $p>1$
such that $F(h_0)$ is $L^p$ on $X$.
There exists an $L_2^p$-function $\alpha$ on $X$
such that
(i) $\Delta_{g_Y}(\alpha)=|F(h_0)_{|X}|_{h_0,g_Y}$,
(ii) $\alpha_{|\del X}=0$.
There exists $C_0>0$
such that $|\alpha|<C_0$ on $X$.

For $0<\delta<1$,
we set $X_P(\delta)=\{|z_P|<\delta\}$
and $Z(\delta):=X\setminus\bigcup_{P\in D}X_P(\delta)$.
We have
$\del Z(\delta)=\del X\cup \bigcup_{P\in D}\del X_P(\delta)$.
By the Dirichlet problem for harmonic metrics
\cite{Donaldson-boundary-value,Li-Mochizuki2},
there exists a harmonic metric $h_{Z(\delta)}$
of $(V,\delbar_V,\theta)_{|Z(\delta)}$
such that
(i) $\det(h_{Z(\delta)})=h_{\det(V)|Z(\delta)}$,
(ii) $h_{Z(\delta)|\del Z(\delta)}=h_{0|\del Z(\delta)}$.
Let $s_{Z(\delta)}$ be the automorphism of $V_{|Z(\delta)}$
determined by
$h_{Z(\delta)}=h_{0|Z(\delta)}\cdot s_{Z(\delta)}$.
According to \cite[Lemma 3.1]{s1},
the following holds on $Z(\delta)$:
\[
 \Delta_{g_Y}\log\Tr(s_{Z(\delta)})
 \leq
 |F(h_0)_{|Z(\delta)}|_{h_0,g_Y}.
\]
Because
$\Delta_{g_Y}\bigl(
 \log \Tr(s_{Z(\delta)})-\alpha
 \bigr)
\leq 0$,
we obtain
\[
 \log\Tr(s_{Z(\delta)})
 \leq
 2C_0+
 \log r.
\]
Because $\det(s_{Z(\delta)})=1$,
there exists $C_1>0$,
which depends only on $C_0$ and $r$,
such that
\[
 |s_{Z(\delta)}|_{h_0}
 +|s_{Z(\delta)}^{-1}|_{h_0}
 \leq C_1.
\]
Then, there exists a sequence $\delta(i)\to 0$ $(i=1,2,\ldots,)$
such that the following holds
(see \cite[Proposition 2.6]{Li-Mochizuki2}):
\begin{itemize}
 \item The sequence $h_{Z(\delta(i))}$
       is convergent
       in the $C^{\infty}$-sense
       on any relatively compact open subset
       of $X\setminus D$.
       Let $h_{\infty}$ denote the limit,
       which is a harmonic metric.
 \item $h_{\infty}$ is mutually bounded with $h_0$.
       As a result,
       $\nbigp^{h_{\infty}}_{\ast}(V)=\nbigp_{\ast}\nbigv_{|X}$.
 \item $\det(h_{\infty})=h_{\det(V)}$.
\end{itemize}

Let $Z:=Z(1/2)$.
There exists a harmonic metric $h_{1,Z}$
of $(V,\delbar_V,\theta)_{|Z}$
such that
(i) $\det(h_{1,Z})=h_{\det(V)|Z}$,
(ii) $h_{1,Z|\del X_P(1/2)}=h_{\infty|\del X_P(1/2)}$
for any $P\in D$,
(iii) $h_{1,Z|\del X}=h_{0|\del X}$.
Let $i_0$ such that $\delta(i_0)<1/2$.
Let $s_{1,\delta(i)}$ be the automorphism of $V_{|Z}$
determined by
$h_{Z(\delta(i))|Z}=h_{1,Z}\cdot s_{1,\delta(i)}$.
We obtain
$\Delta_{g_Y}\log\Tr(s_{1,\delta(i)})\leq 0$
on $Z$.
Hence, we obtain
\[
 \log\bigl(
 \Tr(s_{1,\delta(i)})/r
 \bigr)
 \leq
 \max_{P\in D}
 \max_{Q\in \del X_P(1/2)}
 \bigl\{
 \log\bigl(
 \Tr(s_{1,\delta(i)|Q})/r
 \bigr)
 \bigr\}.
\]
Because $\log\Tr(s_{1,\delta(i)}/r)\to 0$
on $\bigcup_{P\in D}\del X_P(1/2)$,
we obtain that
$s_{1,\delta(i)}\to \id_{V}$ on $Z$.
Hence, we obtain
$h_{\infty|Z}=h_{1,Z}$,
which implies that
$h_{\infty}$ satisfies the condition
$h_{\infty|\del X}=h_{0|\del X}$.

Let $h'$ be another harmonic metric
satisfying the conditions (i) and (ii).
Note that $\det(h')=h_{\det(V)}$.
Let $s$ be the automorphism of $V$
determined by $h'=h\cdot s$.
By \cite[Lemma 3.1]{s1},
we have the following equality on $X\setminus D$:
\[
 \Delta_{g_Y}\Tr(s)
=-\bigl|
  \delbar_V(s)\cdot s^{-1/2}
  \bigr|_{h,g_Y}^2
  -\bigl|
  [\theta,s]s^{-1/2}
  \bigr|_{h,g_Y}^2.
\]
It implies that $\Tr(s)$ is subharmonic on $X\setminus D$.
Because $\Tr(s)$ is bounded,
we obtain that $\Tr(s)$ is a subharmonic function on $X$
(see \cite[Lemma 2.2]{s2}).
We obtain $\max_X\Tr(s)\leq \max_{\del X}\Tr(s)=r$.
Because $\det(s)=1$, we have $\Tr(s)\geq r$.
Hence, we obtain $\Tr(s)=r$ on $X$,
which implies $s=\id_V$.
\hfill\qed

\begin{cor}
\label{cor;22.12.20.40}
Suppose that $(\nbigp_{\ast}\nbigv,\theta)$
is equipped with a perfect symmetric pairing $C$.
If $h_{\del X}$ is compatible with $C_{|\del X}$,
then $h$ is also compatible with $C$.
\end{cor}
\pf
Let $h_{\del X}^{\lor}$ be the Hermitian metric of
$V^{\lor}_{|\del X}$ induced by $h_{\del X}$.
Let $h^{\lor}$ be the Hermitian metric of $V^{\lor}$
induced by $h$.
Then, $h^{\lor}$ is the unique harmonic metric
of $(V^{\lor},\theta^{\lor})$ satisfying
$h^{\lor}_{|\del X}=h_{\del X}^{\lor}$.

Let $\Psi_C:(V,\theta)\simeq (V^{\lor},\theta^{\lor})$
denote the isomorphism induced by $C$.
Because $h_{\del X}$ is compatible with $C$,
$h_{\del X}=\Psi_C^{\ast}h_{\del X}^{\lor}$ holds on $\del X$.
By the uniqueness, we obtain
$h=\Psi_C^{\ast}(h^{\lor})$,
i.e.,
$h$ is compatible with $C$.
\hfill\qed

\section{Large-scale solutions in the symmetric case}
\label{section;23.1.28.1}

\subsection{Preliminary from linear algebra}
\label{subsection;22.12.25.30}

\subsubsection{Hermitian metrics compatible with
a non-degenerate symmetric pairing}

Let $V$ be an $r$-dimensional $\cnum$-vector space.
The dual space is denoted by $V^{\lor}$.
An $\real$-structure of $V$ is
an $r$-dimensional $\real$-subspace $V_{\real}$
such that the natural morphism
$\cnum\otimes_{\real}V_{\real}\lrarr V$
is an isomorphism.
A positive definite symmetric bilinear form $C_{\real}$
of $V_{\real}$
induces a Hermitian metric $h$
and a non-degenerate symmetric bilinear form $C$ of $V$
by
$h(\alpha\otimes u,\beta\otimes v)
=\alpha\betabar h(u,v)$
and
$C(\alpha\otimes u,\beta\otimes v)
=\alpha\beta C(u,v)$
for any $\alpha,\beta\in\cnum$ and $u,v\in V_{\real}$.
An orthogonal decomposition
$V_{\real}=\bigoplus V_{\real,i}$ with respect to $C_{\real}$
induces
a decomposition $V=\bigoplus V_{\real,i}\otimes\cnum$
which is clearly orthogonal with respect to both $h$ and $C$.

Let $C$ be a non-degenerate symmetric pairing of $V$.
It induces a $\cnum$-linear morphism
$\Psi_C:V\to V^{\lor}$.
A Hermitian metric $h$ of $V$
is called compatible with $C$
if $\Psi_C$ is isometry between
$(V,h)$ and $(V^{\lor},h^{\lor})$,
where $h^{\lor}$ denote the Hermitian metric of $V^{\lor}$
induced by $h$.
If $h$ is compatible with $C$,
there uniquely exists an $\real$-structure $V_{\real}$
of $V$ equipped with a positive definite symmetric bilinear form $C_{\real}$
such that (i) $V_{\real}\otimes\cnum=V$,
(ii) $h$ and $C$ are induced by $C_{\real}$.

\subsubsection{An estimate}

Let $C$ be a non-degenerate symmetric form of $V$.
Let $V=\bigoplus_{i=1}^r V_i$ be an orthogonal decomposition
with respect to $C$ such that $\dim V_i=1$.
The following lemma is obvious.
\begin{lem}
There exists a unique Hermitian metric $h_0$ of $V$
such that
(i) $h_0$ is compatible with $C$,
(ii) the decomposition $V=\bigoplus V_i$ is orthogonal
with respect to $h_0$.
\hfill\qed
\end{lem}

For any Hermitian metric $h$ of $V$
compatible with $C$,
let $s(h_0,h)$ be the automorphism of $V$
determined by the condition
$h(u,v)=h_0(s(h_0,h)u,v)$
for any $u,v\in V$.
Note that $\det(s(h_0,h))=1$.
Let $\nbigh(C;\epsilon)$
be the set of Hermitian metrics $h$
of $V$ compatible with $C$
such that the following holds 
for any $u\in V_i$, $v\in V_j$ $(i\neq j)$:
\begin{equation}
\label{eq;22.12.1.1}
      |h(u,v)|\leq \epsilon |u|_h\cdot |v|_h.
\end{equation}

\begin{lem}
\label{lem;22.11.30.210}
There exists $C>0$,
depending only on $r$,
such that the following holds 
for any $0\leq \epsilon\leq (2r)^{-1}$
and any $h\in\nbigh(C;\epsilon)$:
\[
       \bigl|s(h_0,h)-\id_V\bigr|_{h_0}
       +\bigl|s(h_0,h)^{-1}-\id_V\bigr|_{h_0}
       \leq C\epsilon.
\]
\end{lem}
\pf
Let $e_i$ be a base of $V_i$
such that $C(e_i,e_i)=1$.
Note that the tuple $(e_1,\ldots,e_r)$
is an orthonormal base
with respect to $h_0$.
Let $H$ be the matrix determined by
$H_{i,j}=h(e_i,e_j)$.
Then, the linear map $s(h_0,h)$ is represented
by the matrix $\lefttop{t}H$ with respect to
the base $(e_1,\ldots,e_r)$.
Because $h$ is compatible with $C$,
$\lefttop{t}H\cdot H$ is the identity matrix.
We obtain
\begin{equation}
\label{eq;23.1.27.3}
 H_{i,i}^2-1
=\sum_{\substack{1\leq j\leq r \\ j\neq i}}H_{i,j}H_{j,i}.
\end{equation}
By the condition (\ref{eq;22.12.1.1}),
we have 
$|H_{i,j}|
 \leq
 \epsilon H_{i,i}^{1/2}H_{j,j}^{1/2}$
for $i\neq j$.
We obtain
\[
 H_{i,i}^2-1
 \leq
\epsilon \sum_{\substack{1\leq j\leq r \\ j\neq i}}H_{i,i}H_{j,j}.
\]
We set $A=\sum_{j=1}^r H_{j,j}$.
We obtain
\begin{equation}
\label{eq;23.1.27.2}
 H_{i,i}^2-1-\epsilon H_{i,i}A\leq 0.
\end{equation}
\begin{lem}
We obtain
$H_{i,i}\leq \epsilon A+1$.
\end{lem}
\pf
Let $a>0$.
Let us consider
the $\real$-valued function 
$f(s)=s^2-as-1$ $(s\in\real)$.
We set
$s_{\pm}=2^{-1}(a\pm\sqrt{a^2+4})$,
and then we have $f(s_{\pm})=0$ and $s_-<s_+$.
We obtain $f(s)>0$ for any $s>s_+$.
Hence, if $f(s)\leq 0$,
we obtain
\begin{equation}
\label{eq;23.1.27.4}
s\leq s_+\leq 2^{-1}(a+a+2)=a+1.
\end{equation}
By setting $a=\epsilon A$,
we obtain the claim of the lemma from (\ref{eq;23.1.27.2})
and (\ref{eq;23.1.27.4}).
\hfill\qed

\vspace{.1in}

We obtain
$A\leq \epsilon r A+r$,
and hence
$A\leq (1-\epsilon r)^{-1}r\leq 2r$.
By (\ref{eq;22.12.1.1}) and (\ref{eq;23.1.27.3}),
we obtain
\[
 |H_{i,i}^2-1|
 \leq
 \sum_{\substack{1\leq j\leq r\\ j\neq i}}
 |H_{i,j}|\cdot |H_{j,i}|
 \leq
 \epsilon
 \sum_{1\leq j\leq r}
 H_{i,i}\cdot H_{j,j}
 \leq
 \epsilon A^2\leq 4r^2\epsilon.
\]
Because $H_{i,i}$ are positive numbers,
we obtain 
$\bigl|H_{i,i}-1\bigr|\leq 4r^2\epsilon$.
We also obtain
$|H_{i,j}|\leq \epsilon(1+4r^2\epsilon)$.
\hfill\qed

\subsection{Harmonic metrics compatible with a non-degenerate symmetric pairing}

Let $Y$ be any Riemann surface.
Let $(V,\delbar_V,\theta)$ be
a Higgs bundle on $Y$ of rank $r$,
which is regular semisimple.
Let $C$ be a non-degenerate symmetric pairing of $(V,\theta)$.

For any $t>0$,
let $\Harm(V,\delbar_V,t\theta,C)$
denote the set of harmonic metrics of $(V,\delbar_V,t\theta)$
compatible with $C$.
Let $g_Y$ be a K\"ahler metric of $Y$.
For any non-negative integer $\ell$ and $p>1$,
and for any relatively compact open subset $K$ of $Y$,
we define
the $L_{\ell}^p$-norm
$\|f\|_{L_{\ell}^p,K}$
of a section $f$ of $\End(V)$ on $K$
by using
$g_Y$, $h^C$ and the Chern connection of $h^C$.

\begin{thm}
\label{thm;22.12.3.1}
Let $K$ be any relatively compact open subset of $Y$.
There exists $t(K)>0$ such that the following holds
\begin{itemize}
 \item 
For any $(\ell,p)\in\seisuu_{>0}\times\real_{>1}$,
there exist 
$A(\ell,p,K)>0$ and $\epsilon(\ell,p,K)>0$
such that the following holds
for any $h\in \Harm(V,\delbar_V,t\theta,C)$  $(t\geq t(K))$:
\begin{equation}
\label{eq;22.12.1.5}
 \bigl\|
s(h^C,h)-\id_E
 \bigr\|_{L_{\ell}^p,K}
 +
  \bigl\|
s(h^C,h)^{-1}-\id_E
 \bigr\|_{L_{\ell}^p,K}
 \leq A(\ell,p,K)\exp(-\epsilon(\ell,p,K)t).
\end{equation}
\end{itemize}
\end{thm}
\pf
To simplify the description,
we set $s(h):=s(h^C,h)$ in this proof.
By \cite[Corollary 2.6]{Decouple}
and Lemma \ref{lem;22.11.30.210},
there exist $A(K)>0$, $\epsilon(K)>0$ and $t(K)>0$
such that 
the following holds
for any $h\in\Harm(E,\delbar_E,t\theta,C)$ $(t\geq t(K))$:
\begin{equation}
\label{eq;23.1.27.6}
 \sup_K\bigl|
 s(h)-\id_E
 \bigr|_{h^C}
+\sup_K
  \bigl|
  s(h)^{-1}-\id_E
  \bigr|_{h^C}
\leq A(K)\exp(-\epsilon(K)t).
\end{equation}

Let $R(h)$ denote the curvature of
the Chern connection of $(V,\delbar_V,h)$.
By \cite[Theorem 2.9]{Decouple},
there exist $A^{(1)}(K)>0$ and $\epsilon^{(1)}(K)>0$
such that 
the following holds
for any $h\in\Harm(E,\delbar_E,t\theta,C)$ $(t\geq t(K))$:
\begin{equation}
\label{eq;22.12.1.3}
 \sup_K\bigl|
 R(h)
 \bigr|_{h^C,g_Y}
\leq A^{(1)}(K)\exp(-\epsilon^{(1)}(K)t).
\end{equation}
Note that
$R(h)=\delbar_V\bigl(s(h)^{-1}\del_{h^C}s(h)\bigr)$.

Because $s(h)$ is self-adjoint with respect to $h^C$
and satisfies $\det s(h)=1$,
we have $\Tr(s(h)-\id)\geq 0$,
and $\Tr(s(h)-\id_E)=0$ holds if and only if
$s(h)=\id_E$.
The following holds on $Y$ (see \cite[Lemma 3.1]{s1}): 
\[
 \Delta_{g_Y}
 \Tr\bigl(s(h)-\id_E\bigr)
=
 \Delta_{g_Y}\bigl(
 \Tr(s(h))-r
 \bigr)
=-\bigl|
 s(h)^{-1/2}\del_{h^C}s(h)
 \bigr|^2_{g_Y,h^C}
 -\bigl|
 [\theta,s(h)]
 s(h)^{-1/2}
 \bigr|^2_{g_Y,h^C}.
\]
Let $K_1$ be a relatively compact open neighbourhood of
$\Kbar$ in $Y$.
Let $\chi:Y\to \real_{\geq 0}$ be a function such that
$\chi=1$ on $\Kbar$
and $\chi=0$ on $Y\setminus K_1$.
We obtain
the following:
\[
\int_K
\bigl|
 s(h)^{-1/2}\del_{h^C}s(h)
 \bigr|^2_{g_Y,h^C}
 \leq
 \int_{Y}
 \Tr(s(h)-\id_E)
 \cdot
 \bigl|
 \Delta_{g_Y}\chi
 \bigr|.
\]
There exist $A^{(2)}(K)>0$ and $\epsilon^{(2)}(K)>0$
such that the following holds
for any $h\in\Harm(V,\delbar_V,t\theta,C)$
$(t\geq t(K))$:
\begin{equation}
\label{eq;22.12.1.4}
 \int_K
\bigl|
 s(h)^{-1}\del_{h^C}s(h)
 \bigr|^2_{g_Y,h^C}
 \leq
 A^{(2)}(K)\exp(-\epsilon^{(2)}(K)t).
\end{equation}
By (\ref{eq;22.12.1.3}) and (\ref{eq;22.12.1.4}),
there exist 
$A^{(3)}(p,K)>0$ and $\epsilon^{(3)}(p,K)>0$
such that the following holds
for any $h\in\Harm(V,\delbar_V,t\theta,C)$
$(t\geq t(K))$:
\begin{equation}
\label{eq;23.1.27.5}
\bigl\|
 s(h)^{-1}\del_{h^C}s(h)
 \bigr\|_{L_1^p,K}
 \leq
 A^{(3)}(K)\exp(-\epsilon^{(3)}(p,K)t).
\end{equation}
By (\ref{eq;23.1.27.6}) and (\ref{eq;23.1.27.5}),
there exist 
$A^{(4)}(p,K)>0$ and $\epsilon^{(4)}(p,K)>0$
such that the following holds
for any $h\in\Harm(V,\delbar_V,t\theta,C)$
$(t\geq t(K))$:
\begin{equation}
\label{eq;23.1.27.7}
\bigl\|
 s(h)-\id
 \bigr\|_{L_1^p,K}
 \leq
 A^{(4)}(K)\exp(-\epsilon^{(4)}(p,K)t).
\end{equation}
By (\ref{eq;23.1.27.5}) and (\ref{eq;23.1.27.7}),
there exist 
$A^{(5)}(p,K)>0$ and $\epsilon^{(5)}(p,K)>0$
such that the following holds
for any $h\in\Harm(V,\delbar_V,t\theta,C)$
$(t\geq t(K))$:
\begin{equation}
\bigl\|
 s(h)-\id
 \bigr\|_{L_2^p,K}
 \leq
 A^{(5)}(K)\exp(-\epsilon^{(5)}(p,K)t).
\end{equation}
Then, by using a standard bootstrapping argument,
we obtain the claim of the proposition.
\hfill\qed

\begin{cor}
\label{cor;22.12.20.31}
Let $t(i)>0$ be any sequence such that $\lim_{i\to\infty} t(i)=\infty$.
Let $h_{t(i)}\in\Harm(V,\delbar_V,t(i)\theta,C)$.
Then, $h_{t(i)}$ is convergent to $h^C$
in the $C^{\infty}$-sense
on any relatively compact open subsets of $Y$.
The order of the convergence is estimated as in
{\rm(\ref{eq;22.12.1.5})}.
\hfill\qed
\end{cor}

\section{Some estimates for harmonic bundles on a disc}

This section is preliminary for 
Theorem \ref{thm;22.11.28.20}.

\subsection{Universal boundedness of higher derivatives of Higgs fields}

For any $R>0$,
we set $B(R)=\bigl\{z\in\cnum\,\big|\,|z|<R\bigr\}$.
Let $R_0>0$.
Let $(E,\delbar_E,\theta)$ be a Higgs bundle on $B(R_0)$ of rank $r$.
Let $f$ be the endomorphism of $E$ determined by $\theta=f\,dz$.
Let $C_0$ be a constant such that
\[
 |\tr(f^j)|<C_0
 \quad
 (j=1,\ldots,r).
\]

Let $h$ be a harmonic metric of $(E,\delbar_E,\theta)$.
Let $\nabla_h$ denote the Chern connection of $h$.
Let $R(h)$ denote the curvature of $\nabla_h$.
We obtain the endomorphism $\gbigr$ determined by
$R(h)=\gbigr\,dz\,d\zbar$.
Let $f^{\dagger}_h$ denote the adjoint of $f$
with respect to $h$.
Because
$R(h)+[\theta,\theta^{\dagger}_h]=0$,
we have
$\gbigr+[f,f^{\dagger}_h]=0$.

Let $g_0=dz\,d\zbar$ denote the standard Euclidean metric.
We consider the $L_{\ell}^p$-norm of
sections of $\End(E)$ with respect to
$g_0$, $h$ and the derivatives with respect to $\nabla_h$.

\begin{prop}
\label{prop;22.12.2.10}
Let $0<R_1<R_0$.
For any $\ell\in\seisuu_{\geq 0}$ and $p\geq 1$,
there exist $C(\ell,p)$, depending only on
$r$, $R_0$, $R_1$ and $C_0$,
such that
\[
 \|f_{|B(R_1)}\|_{L_{\ell}^p}
+\|f^{\dagger}_{h|B(R_1)}\|_{L_{\ell}^p}
+\|\gbigr_{|B(R_1)}\|_{L_{\ell}^p}
 \leq C(\ell,p).
\]
\end{prop}
\pf
Let $R_2=(R_0+R_1)/2$.
By Simpson's main estimate \cite{s2,s5},
there exists $C_1$, depending only on 
$r$, $R_0$, $R_1$ and $C_0$,
such that
$|f|_{h}=|f^{\dagger}_h|_h\leq C_1$ on $B(R_2)$.
We also obtain
$|R(h)|_{h,g_0}=|\gbigr|_{h}\leq 2C_1^2$ on $B(R_2)$.

We recall a result due to Uhlenbeck.
\begin{thm}[\mbox{\cite[Theorem 1.3]{u1}}]
\label{thm;22.12.2.1}
Let $V$ be a vector bundle on $B(1)$
equipped with a Hermitian metric $h_V$
and a unitary connection $\nabla_V$.
Let $R(\nabla_V)$ denote the curvature of $\nabla_V$.
For $p\geq 1$,
let $\|R(\nabla_V)\|_{L^p,h_V}$ denote the $L^p$-norm
with respect to $g_0$ and $h_V$.
Then, there exist positive constants $c$  and $\kappa$
depending only on $r$ and $p$ such that
the following holds
\begin{itemize}
 \item If $\|R(\nabla_V)\|_{L^p,h_V}\leq \kappa$,
       then
       there exists an orthonormal frame
       $\vecv$ of $V$
       such that the connection form $A$
       of $\nabla_V$ with respect to $\vecv$
       satisfies
       (i) $d^{\ast}A=0$,
       (ii) $\|A\|_{L_1^p}\leq c\|R(\nabla_V)\|_{L^p}$.
       \hfill\qed
\end{itemize}
\end{thm}

We choose $T>0$
such that
$100T^{-1}C_1^2<\kappa$
and $T(R_0-R_2)>100$.
Let $\varphi_T:\cnum_w\to\cnum_z$
be defined by $\varphi_T(w)=T^{-1}z$.
We consider
$(\Etilde,\delbar_{\Etilde},\thetatilde,\htilde)
=\varphi_T^{\ast}(E,\delbar_E,\theta,h)$ on $B(TR_0)$.
Let $w_0\in B(TR_2-1)$.
Let $p>2$.
Let $\vecv^{(w_0)}$
be an orthonormal frame of
$\Etilde_{|D(w_0,1)}$
as in Theorem \ref{thm;22.12.2.1}
for the metric $\htilde$
and the connection $\nabla_{\htilde}$.
Let $\nbiga^{(w_0)}$ and $\nbigr^{(w_0)}$
denote the connection form
and the curvature form of $\nabla_{\htilde}$
with respect to $\vecv^{(w_0)}$.
We have
\begin{equation}
\label{eq;22.12.2.2}
 d^{\ast}\nbiga^{(w_0)}=0,
 \quad
 d\nbiga^{(w_0)}+\nbiga^{(w_0)}\wedge\nbiga^{(w_0)}
 =\nbigr^{(w_0)},
\end{equation}
\begin{equation}
\label{eq;22.12.2.3}
 \|\nbiga^{(w_0)}\|_{L_1^p(D(w_0,1))}
 \leq
 c\|\nbigr \|_{L^p(D(w_0,1))}. 
\end{equation}
Let $\Theta^{(w_0)}$
denote the matrix valued $(1,0)$-form determined by
$\thetatilde\vecv^{(w_0)}=\vecv^{(w_0)}\Theta^{(w_0)}$.
We have the decomposition
$\nbiga^{(w_0)}=\nbiga^{(w_0)}_w\,dw+
\nbiga^{(w_0)}_{\wbar}\,d\wbar$.
We have
$\nbiga^{(w_0)}_{w}=
-\lefttop{t}\overline{\nbiga^{(w_0)}_{\wbar}}$.
Because $\delbar\thetatilde=0$,
the following holds.
\begin{equation}
\label{eq;22.12.2.4}
\del_{\wbar}\Theta^{(w_0)}
+[\nbiga^{(w_0)}_{\wbar},\Theta^{(w_0)}]=0.
\end{equation}
We also have
\begin{equation}
 \nbigr^{(w_0)}+[\Theta^{(w_0)},\lefttop{t}\overline{\Theta^{(w_0)}}]=0.
\end{equation}
Then, by a standard bootstrapping argument,
we can prove that
for any $\ell$ there exists
$C_2(\ell)$, depending only on $\ell$ and $r$
such that
\[
\bigl\|
\Theta^{(w_0)}
\bigr\|_{L^p_{\ell}(D(w_0,1/2))}
+
\bigl\|
\nbiga^{(w_0)}
\bigr\|_{L^p_{\ell+1}(D(w_0,1/2))}
\leq C_2(\ell).
\]
Then, we obtain a desired estimate for
$\|f_{|B(R_1)}\|_{L_{\ell}^p}$,
which implies a desired estimate for
$\|f^{\dagger}_{h|B(R_1)}\|_{L_{\ell}^p}$.
Because
$\gbigr+[f,f^{\dagger}_h]=0$,
we also obtain a desired estimate for
$\|\gbigr_{|B(R_1)}\|_{L_{\ell}^p}$.
\hfill\qed

\subsection{Difference of two families of large-scale solutions on a disc}
\label{subsection;22.12.3.10}

Let $R_0>0$.
Let $(E,\delbar_E,\theta)$ be a Higgs bundle on $B(R_0)$ of rank $r$.
Let $f$ be the endomorphism of $E$ determined by $\theta=f\,dz$.
Let $C_0$ be a constant such that
\[
 |\tr(f^j)|<C_0
 \quad
 (j=1,\ldots,\rank(E)).
\]

Let $h_{\det(E)}$ be a flat metric of $\det(E)$.
Let $h_{0,t}$ $(t>0)$ be harmonic metrics of $(E,\delbar_E,t\theta)$
such that $\det(h_{0,t})=h_{\det(E)}$.
Let $\nabla^{0,t}$ denote the Chern connection of
$(E,\delbar_E,h_{0,t})$.
For any section $u$ of $\End(E)$
and for
any element $\veckappa=(\kappa_1,\kappa_2,\ldots,\kappa_{\ell})
\in\{z,\zbar\}^{\ell}$,
we set
\[
 \nabla^{0,t}_{\veckappa}u
 =\nabla^{0,t}_{\kappa_1}\circ
 \nabla^{0,t}_{\kappa_2}\circ
 \cdots\circ
 \nabla^{0,t}_{\kappa_{\ell}}(u).
\]

\begin{thm}
\label{thm;22.12.2.11}
Let $0<R_1<R_0$. Let $C_1,\epsilon_1>0$.
For any $\ell\in\seisuu_{\geq 0}$,
there exist positive constants $C(\ell),\epsilon(\ell)>0$,
depending only on $r$, $C_0,C_1,\epsilon_1$ and $\ell$
such that the following holds.
\begin{itemize}
 \item  Let $t(i)>0$ be an increasing sequence such that $t(i)\to\infty$ as
	$i\to\infty$.
	We also assume that $t(1)(R_0-R_1)>100$.
	Let $h_{t(i)}$ be harmonic metrics of
	$(E,\delbar_E,t(i)\theta)$
	such that
	$\det(h_{t(i)})=h_{\det(E)}$.
	Assume the following on
	$B(R_0)\setminus B(R_1)$:
\begin{equation}
\label{eq;22.12.1.10}
	\bigl|s(h_{0,t(i)},h_{t(i)})-\id
	\bigr|_{h_{0,t(i)}}
	\leq
	C_1\exp(-\epsilon_1t(i)).
\end{equation}
	Then, the following holds on $B(R_1)$
	for any $\veckappa\in\{z,\zbar\}^{\ell}$:
\[
	\Bigl|
	\nabla^{0,t(i)}_{\veckappa}\bigl(
	s(h_{0,t(i)},h_{t(i)})-\id
	\bigr)
	\Bigr|_{h_{0,t(i)}}
	\leq
	C(\ell)
	\exp\bigl(-\epsilon(\ell)t(i)\bigr).
\]
\end{itemize}
\end{thm}

\subsubsection{The case $\ell=0$}

To simplify the notation
we set $s_i=s(h_{0,t(i)},h_{t(i)})$.
By (\ref{eq;22.12.1.10}),
there exist $C'(0),\epsilon'(0)>0$,
depending only on
$r$, $C_1$, and $\epsilon_1$
such that the following holds
on $B(R_0)\setminus B(R_1)$:
\begin{equation}
\label{eq;22.12.1.11}
\Tr(s_i-\id_{E})
 \leq C'(0)\exp(-\epsilon'(0)t(i))
\end{equation}
By \cite[Lemma 3.1]{s1},
we have
\begin{equation}
\label{eq;22.12.1.31}
 -\del_z\del_{\zbar}\Tr(s_i-\id_E)
 =-\bigl|
 \delbar(s_i) s_i^{-1/2}
 \bigr|^2_{h_{0,t(i)}}
 -\bigl|
 [t\theta,s_i]s_i^{-1/2}
 \bigr|^2_{h_{0,t(i)}}.
\end{equation}
In particular,
$\Tr\bigl(s_i-\id_E\bigr)$
is a subharmonic function on $B(R_0)$.
By the maximum principle
of subharmonic functions,
(\ref{eq;22.12.1.11}) holds on $B(R_0)$.
Because $\det(s_i)=1$,
we obtain the claim in the case $\ell=0$.

\subsubsection{Estimates for $L^2$-norms}

We set $R_2=(R_0+R_1)/2$
and $R_3=(R_0+R_2)/2$.
Let $\chi:\cnum\to\real_{\geq 0}$
be a $C^{\infty}$-function such that
$\chi(z)=1$ $(|z|\leq R_2)$
and
$\chi(z)=0$ $(|z|\geq R_3)$.
Let $g_z=dz\,d\zbar$ be the standard Euclidean metric.
By using \cite[Lemma 3.1]{s1},
we obtain
\[
\int_{B(R_2)}
\Bigl(
 \bigl|
 \delbar(s_i) s_i^{-1/2}
 \bigr|^2_{h_{0,t(i)},g_z}
+\bigl|
 [t\theta,s_i]s_i^{-1/2}
 \bigr|^2_{h_{0,t(i)},g_z}
 \Bigr)\vol_{g_z}
 \leq
 \int_{B(R_3)\setminus B(R_2)}
 \bigl|
 \del_z\del_{\zbar}(\chi)
 \bigr|\cdot
 \bigl(
 \Tr(s_i-\id_E)
 \bigr)\vol_{g_z}.
\]
Hence, there exist $C_5>0,\epsilon_5>0$
such that
\begin{equation}
\label{eq;22.12.1.32}
\int_{B(R_2)}
\Bigl(
 \bigl|
 \delbar(s_i)
  s_i^{-1}
 \bigr|^2_{h_{0,t(i)},g_z}
+\bigl|
 s_i^{-1}[t\theta,s_i]
 \bigr|^2_{h_{0,t(i)},g_z}
 \Bigr)\vol_{g_z}
 \leq
C_5\exp(-\epsilon_5t(i)).
\end{equation}

\subsubsection{Rescaling}

To study the derivatives,
for any $t>t(1)$,
we define the map
$\rho_t:\cnum_w\to\cnum_z$
by $\rho_t(w)=t^{-1}w$.
We have $\rho_t^{-1}(B(R))=B(tR)$.
We use the standard Euclidean metric $g_w=dw\,d\wbar$
on $\cnum_w$.

We set
$\Etilde_{t}=\rho_t^{\ast}(E)$ on $B(tR_0)$.
It is equipped with the Higgs field
$\thetatilde_t=\rho_t^{\ast}(t\theta)$.
We have
$\thetatilde_t=\rho_t^{\ast}(f)\,dw$.
We have the harmonic metrics
$\htilde_{0,t}=\varphi_t^{\ast}(h_{0,t})$
of the Higgs bundles
$(\Etilde_t,\delbar_{\Etilde_t},\thetatilde_t)$.
Let $\nablatilde^{0,t}$ denote the Chern connection of
$(\Etilde_t,\delbar_{\Etilde_t},\htilde_{0,t})$.

By Simpson's main estimate,
there exists $C_{10}>0$,
depending only on $r$ and $C_0$
such that the following holds
on $B(tR_0-1)$:
\begin{equation}
\label{eq;22.12.1.30}
 \bigl|
 \thetatilde_{t}
 \bigr|_{\htilde_{0,t},g_w}
 \leq
 C_{10}.
 \end{equation}

Let $R(\htilde_{0,t})$
denote the curvature of the Chern connection of
$(\Etilde_{t},\delbar_{\Etilde_{t}},
\htilde_{0,t})$.
We have the following equality:
\begin{equation}
\label{eq;22.12.1.50}
 R(\htilde_{0,t})
 +\bigl[
 \thetatilde_{t},
 (\thetatilde_{t})^{\dagger}_{\htilde_{0,t}}
 \bigr]
 =0.
\end{equation}
By (\ref{eq;22.12.1.30}) and (\ref{eq;22.12.1.50}),
we have the following on $B(tR_0-1)$:
\begin{equation}
\label{eq;22.12.1.41}
 \bigl|
 R(\htilde_{0,t})
 \bigr|_{\htilde_{0,t},g_w}
 \leq
 2C_{10}^2.
\end{equation}
We also have the universal estimates
for higher derivatives
of $\thetatilde$ and $R(\htilde_{0,t})$
as in Proposition \ref{prop;22.12.2.10}.

\subsubsection{Estimates for higher derivatives}

We also have the harmonic metrics
$\htilde_{t(i)}$
of
$(\Etilde_{t(i)},\delbar_{\Etilde_{t(i)}},\thetatilde_{t(i)})$.
Let $\stilde_i=\varphi_{t(i)}^{\ast}(s_i)$.
We have
$\htilde_{t(i)}=\htilde_{0,t(i)}\stilde_i$.
By (\ref{eq;22.12.1.32}),
we have
\begin{equation}
\label{eq;22.12.1.33}
\int_{B(tR_2)}
\Bigl(
\bigl|
\delbar(\stilde_i)
\stilde_i^{-1}
 \bigr|^2_{\htilde_{0,t(i)},g_w}
+\bigl|
\stilde_i^{-1}
 [\thetatilde_t,\stilde_i]
 \bigr|^2_{\htilde_{0,t(i)},g_w}
 \Bigr)\vol_{g_w}
 \leq
C_5\exp(-\epsilon_5t(i)).
\end{equation}
It implies
\begin{equation}
 \label{eq;22.12.1.34}
\int_{B(tR_2)}
 \bigl|
 \stilde_i^{-1}\del_{\htilde_{0,t(i)}}(\stilde_i)
 \bigr|^2_{\htilde_{0,t(i)},g_w}\vol_{g_w}
 \leq
C_5\exp(-\epsilon_5t(i)).
\end{equation}

Let $R(\htilde_{t(i)})$
denote the curvature of the Chern connection of
$(\Etilde_{t(i)},\delbar_{\Etilde_{t(i)}},\htilde_{t(i)})$.
We have
\[
 R(\htilde_{t(i)})
 +\bigl[
 \thetatilde_{t(i)},
 (\thetatilde_{t(i)})^{\dagger}_{\htilde_{t(i)}}
 \bigr] 
=0.
\]
Note that
\[
 (\thetatilde_{t(i)})^{\dagger}_{\htilde_{t(i)}}
=\stilde_{t(i)}^{-1}
 (\thetatilde_{t(i)})^{\dagger}_{\htilde_{0,t(i)}}
 \stilde_{t(i)}.
\]
We obtain
\begin{multline}
 \delbar\bigl(
 \stilde_i^{-1}
 \del_{\htilde_{0,t(i)}}
 \stilde_i
 \bigr)
=R(\htilde_{t(i)})-
 R(\htilde_{0,t(i)})
=-\Bigl[
\thetatilde_{t(i)},
\stilde_i^{-1}
(\thetatilde_{t(i)})^{\dagger}_{\htilde_{0,t(i)}}
\stilde_i
-(\thetatilde_{t(i)})^{\dagger}_{\htilde_{0,t(i)}}
 \Bigr]
 \\
=
-\Bigl[
\thetatilde_{t(i)},
\stilde_i^{-1}
\bigl[
(\thetatilde_{t(i)})^{\dagger}_{\htilde_{0,t(i)}},
\stilde_i-\id
\bigr]
\Bigr].
\end{multline}
Hence, there exist $C_{11}>0$ and $\epsilon_{11}>0$
such that the following holds
on $B(tR_0-1)$:
\begin{equation}
\label{eq;22.12.1.40}
\bigl|
 \delbar\bigl(
 \stilde_i^{-1}\del_{\htilde_{0,t(i)}}\stilde_i
 \bigr)
\bigr|_{\htilde_{0,t(i)},g_w}
\leq
C_{11}\exp(-\epsilon_{11}t(i)).
\end{equation}
For any $w_0\in \cnum_w$,
we set
$D(w_0,T)=\{|w-w_0|<T\}$.
By (\ref{eq;22.12.1.41}),
(\ref{eq;22.12.1.34}),
and (\ref{eq;22.12.1.40}),
for any $p\geq 2$,
there exist $C_{12}(p)>0,\epsilon_{12}(p)>0$
such that the following holds
for any $w_0\in B(tR_2-1)$:
\begin{equation}
\label{eq;22.12.1.42}
 \bigl|
 \stilde_i^{-1}
 \del_{\htilde_{0,t(i)}}
 \stilde_i
 \bigr|_{L_1^p(D(w_0,2/3)),\htilde_{0,t(i)},g_w}
 \leq
 C_{12}(p)\exp\bigl(-\epsilon_{12}(p)t(i)\bigr).
\end{equation}
By (\ref{eq;22.12.1.42}) and the estimate in the case $\ell=0$,
for any $p>1$,
there exist $C_{13}(p)>0,\epsilon_{13}(p)>0$
such that the following holds
for any $w_0\in B(tR_2-1)$:
\begin{equation}
 \bigl|
 \stilde_i-\id
 \bigr|_{L_2^p(D(w_0,2/3)),\htilde_{0,t(i)},g_w}
 \leq
 C_{13}(p)\exp\bigl(-\epsilon_{13}(p)t(i)\bigr).
\end{equation}
By a standard bootstrapping argument,
for any
for any $p>1$ and $\ell\in\seisuu_{\geq 2}$,
there exist $C_{14}(\ell,p)>0,\epsilon_{14}(\ell,p)>0$
such that the following holds
for any $w_0\in B(tR_2-1)$:
\begin{equation}
 \bigl|
 \stilde_i-\id
 \bigr|_{L_{\ell}^p(D(w_0,1/2)),\htilde_{0,t(i)},g_w}
 \leq
 C_{14}(\ell,p)\exp\bigl(-\epsilon_{14}(\ell,p)t(i)\bigr).
\end{equation}
Then, we obtain the claim of Theorem \ref{thm;22.12.2.11}.
\hfill\qed

\section{Decomposable filtered extensions}

\subsection{Meromorphic extensions and filtered extensions}

\subsubsection{Vector bundles}
\label{subsection;22.12.22.10}

Let $U\subset\cnum$ be a simply connected open neighbourhood of $0$.
We set $U^{\ast}=U\setminus\{0\}$.
Let $\iota:U^{\ast}\to U$ denote the inclusion.
Let $V$ be a locally free $\nbigo_{U^{\ast}}$-module of rank $r$.
We obtain a locally free $\iota_{\ast}\nbigo_{U^{\ast}}$-module
$\iota_{\ast}(V)$.
A meromorphic (resp. smooth) extension of $V$ is defined to be
a locally free $\nbigo_{U}(\ast 0)$-submodule
(resp. $\nbigo_U$-submodule)
$\nbigv\subset \iota_{\ast}(V)$
such that
$\nbigv_{|U^{\ast}}=V$.
A filtered extension of $V$
is defined to be
a meromorphic extension $\nbigv$
equipped with a filtered bundle
$\nbigp_{\ast}(\nbigv)$ over $\nbigv$.

\begin{example}
The $\nbigo_U(\ast 0)$-submodule
$\nbigo_U(\ast 0)\exp(z^{-1})\subset \iota_{\ast}(\nbigo_{U^{\ast}})$
is a  meromorphic extension of $\nbigo_{U^{\ast}}$,
which is different from
$\nbigo_U(\ast 0)\subset \iota_{\ast}(\nbigo_{U^{\ast}})$.
\hfill\qed
\end{example}

For a positive integer $\ell$,
let $\varphi_{\ell}:\cnum\to\cnum$ be defined by
$\varphi_{\ell}(\zeta)=\zeta^{\ell}$.
We set $U^{(\ell)}=\varphi_{\ell}^{-1}(U)$
and $U^{(\ell)\ast}=U^{(\ell)}\setminus\{0\}$.
The induced morphisms
$U^{(\ell)}\to U$
and $U^{(\ell)\ast}\to U^{\ast}$
are also denoted by $\varphi_{\ell}$.
Let $\Gal(\ell)$ denote the Galois group
of the ramified covering $\varphi_{\ell}$.
Namely,
we put $\Gal(\ell)=\{a\in\cnum^{\ast}\,|\,a^{\ell}=1\}$,
and we consider the action of $\Gal(\ell)$
on $U^{(\ell)}$ by the multiplication on the coordinate $\zeta$.
Let $\iota^{(\ell)}:U^{(\ell)\ast}\to U^{(\ell)}$
denote the inclusion.
We set $V^{(\ell)}:=\varphi_{\ell}^{\ast}(V)$,
which is naturally $\Gal(\ell)$-equivariant.
The $(\iota^{(\ell)})_{\ast}\nbigo_{U^{(\ell)\ast}}$-module
$(\iota^{(\ell)})_{\ast}(V^{(\ell)})$
is also $\Gal(\ell)$-equivariant.
A $\Gal(\ell)$-equivariant meromorphic extension
of $V^{(\ell)}$ is defined to be
a locally free $\nbigo_{U^{(\ell)}}(\ast 0)$-submodule
$\nbigv^{(\ell)}\subset\iota^{(\ell)}_{\ast}(V^{(\ell)})$
which is preserved by the $\Gal(\ell)$-action.
A $\Gal(\ell)$-equivariant filtered extension
of $V^{(\ell)}$ is defined to be
a filtered bundle
$\nbigp_{\ast}(\nbigv^{(\ell)})$
over a $\Gal(\ell)$-equivariant meromorphic extension
$\nbigv^{(\ell)}$ of $V^{(\ell)}$
such that
each $\nbigp_{a}\nbigv^{(\ell)}$
is preserved by the $\Gal(\ell)$-action.

A meromorphic extension $\nbigv$ of $V$
induces a $\Gal(\ell)$-equivariant meromorphic extension
$\varphi_{\ell}^{\ast}(\nbigv)$ of $V^{(\ell)}$.
Conversely,
for any $\Gal(\ell)$-equivariant meromorphic extension
$\nbigv^{(\ell)}$ of $V^{(\ell)}$,
we obtain
the $\nbigo_{U}(\ast 0)$-module
$\varphi_{\ell\ast}(\nbigv^{(\ell)})$
equipped with the $\Gal(\ell)$-action.
The $\Gal(\ell)$-invariant part 
$\varphi_{\ell\ast}(\nbigv^{(\ell)})^{\Gal(\ell)}$
is called the descent of $\nbigv^{(\ell)}$
which is a meromorphic extension of $V$.
\begin{lem}
For a meromorphic extension $\nbigv$ of $V$,
the descent of $\varphi_{\ell}^{\ast}(\nbigv)$
equals $\nbigv$.
For a $\Gal(\ell)$-equivariant meromorphic extension
$\nbigv^{(\ell)}$ of $V^{(\ell)}$,
 $\varphi_{\ell}^{\ast}\bigl(
 \varphi_{\ell\ast}(\nbigv^{(\ell)})^{\Gal(\ell)}
 \bigr)$
equals $\nbigv^{(\ell)}$.
These procedures induce
an equivalence between
meromorphic extensions of $V$ and
$\Gal(\ell)$-equivariant meromorphic extensions of $V^{(\ell)}$.
\hfill\qed
\end{lem}

For a filtered extension $\nbigp_{\ast}\nbigv$ of $V$,
we obtain a $\Gal(\ell)$-equivariant filtered extension
$\nbigp_{\ast}(\varphi_{\ell}^{\ast}(\nbigv))$
over $\nbigv^{(\ell)}$
as follows:
\[
 \nbigp_a(\varphi_{\ell}^{\ast}\nbigv)
 =\sum_{\substack{b\in\real, k\in\seisuu\\
 \ell b+k\leq a}}
 \zeta^{-k}
 \varphi_{\ell}^{\ast}(\nbigp_{b}\nbigv)
 \subset
 \varphi_{\ell}^{\ast}(\nbigv).
\]
The filtered bundle
$\nbigp_{\ast}(\varphi_{\ell}^{\ast}(\nbigv))$
is denoted by
$\varphi_{\ell}^{\ast}(\nbigp_{\ast}\nbigv)$.

For a $\Gal(\ell)$-equivariant filtered extension
$\nbigp_{\ast}(\nbigv^{(\ell)})$ of $V^{(\ell)}$,
we obtain a filtered extension
$\nbigp_{\ast}(\varphi_{\ell\ast}(\nbigv^{(\ell)})^{\Gal(\ell)})$
as follows:
\[
\nbigp_{a}(\varphi_{\ell\ast}(\nbigv^{(\ell)})^{\Gal(\ell)})
=\varphi_{\ell\ast}(\nbigp_{\ell a}\nbigv^{(\ell)})^{\Gal(\ell)}.
\]
It is called the descent of $\nbigp_{\ast}(\nbigv^{(\ell)})$.

\begin{lem}
These procedures induce an equivalence
between
filtered extension of $V$ 
and $\Gal(\ell)$-equivariant filtered extension of $V^{(\ell)}$.
\hfill\qed 
\end{lem}

\subsubsection{Non-degenerate symmetric pairing}

For any $b\in\real$,
let $\nbigp^{(b)}_{\ast}(\nbigo_U(\ast 0))$
denote the filtered bundle over $\nbigo_U(\ast 0)$
defined by
\[
 \nbigp^{(b)}_a(\nbigo_U(\ast 0))
 =z^{-[a-b]}\nbigo_U.
\]

Let $C:V\otimes V\to\nbigo_{U^{\ast}}$
be a holomorphic non-degenerate symmetric pairing.
We say that
a meromorphic extension $\nbigv$ is compatible with $C$
if $C$ extends to a pairing
$\nbigv\otimes\nbigv\to\nbigo_{U}(\ast 0)$.
We say that
a filtered extension $\nbigp_{\ast}\nbigv$ is compatible with $C$
if $C$ induces 
$\nbigp_{\ast}\nbigv\otimes\nbigp_{\ast}\nbigv
\to\nbigp^{(0)}_{\ast}(\nbigo_U(\ast 0))$.
We say that $C$ is perfect with respect to $\nbigp_{\ast}\nbigv$
if $C$ induces an isomorphism
$\nbigp_{\ast}(\nbigv)\simeq\nbigp_{\ast}(\nbigv^{\lor})$.

We have the induced symmetric pairing $\det(C)$ of $\det(V)$.
If $\nbigv$ (resp. $\nbigp_{\ast}\nbigv$) is compatible with $C$,
then $\det(\nbigv)$ (resp. $\det(\nbigp_{\ast}\nbigv)$)
is compatible with $\det(C)$.
\begin{lem}[\cite{Li-Mochizuki3}]
Suppose that $\nbigp_{\ast}\nbigv$ is compatible with $C$.
Then, $C$ is perfect with respect to $\nbigp_{\ast}(\nbigv)$
if and only if $\det(C)$ is perfect with respect to
$\det(\nbigp_{\ast}\nbigv)$.
\hfill\qed
\end{lem}

\begin{lem}
\label{lem;22.12.22.40}
There exists a unique meromorphic extension $\nbigl$
of $\det(V)$
which is compatible with $\det(C)$.
There exists a unique filtered bundle
$\nbigp^{C}_{\ast}\nbigl$ over $\nbigl$
such that 
$\det(C)$ is perfect with respect to
$\nbigp^C_{\ast}\nbigl$.
\end{lem}
\pf
We may assume that $U$ is a disc.
Let $v_0$ be a frame of $\det(V)$ on $U^{\ast}$.
We obtain a holomorphic function
$(\det C)(v_0,v_0)$ on $U^{\ast}$.
There exist an integer $k$
and a holomorphic function $g_1$
such that
$(\det C)(v_0,v_0)=z^{-k}\exp(g_1)$.
We obtain a frame $v_1=\exp(-g_1/2)v_0$
of $\det(V)$ on $U^{\ast}$.
We set
$\nbigl=\nbigo_{U}(\ast 0)v_1
\subset \iota_{\ast}(V)$.
Then, $\nbigl$ is compatible with $\det(C)$.

We have 
$\det(C)(v_1,v_1)=z^{-k}$.
We define
\[
 \nbigp^C_{a}(\nbigl)
=z^{-[a-k/2]}\nbigo_U\cdot v_1.
\]
Then, $\nbigp^C_{\ast}\nbigl$
satisfies the desired condition.
The uniqueness is clear.
\hfill\qed

\vspace{.1in}

We set $C^{(\ell)}:=\varphi_{\ell}^{\ast}C$
which is a non-degenerate symmetric pairing of $V^{(\ell)}$.
\begin{lem}
$\nbigv$ (resp. $\nbigp_{\ast}\nbigv$)
is compatible with $C$
if and only if $\varphi_{\ell}^{\ast}(\nbigv)$
(resp. $\varphi_{\ell}^{\ast}(\nbigp_{\ast}\nbigv)$)
is compatible with $C^{(\ell)}$.
When $\nbigp_{\ast}\nbigv$ and $C$ are compatible,
$C$ is perfect with respect to $\nbigp_{\ast}\nbigv$
if and only if $C^{(\ell)}$ is perfect
with respect to $\varphi_{\ell}^{\ast}(\nbigp_{\ast}\nbigv)$.
\hfill\qed
\end{lem}

\subsubsection{Higgs bundles}

Let $\theta$ be a Higgs field of $V$,
i.e.,
$\theta:V\to V\otimes\Omega^1_{U^{\ast}}$.
We obtain
$\iota_{\ast}(\theta):\iota_{\ast}(V)\to
\iota_{\ast}(V)\otimes\Omega^1_{U}$.
A meromorphic (resp. smooth) extension of $(V,\theta)$
is defined to be a meromorphic (resp. smooth) extension $\nbigv$ of $V$
such that
$\iota_{\ast}(\theta)(\nbigv)\subset\nbigv\otimes\Omega^1_U$.
The induced Higgs field of $\nbigv$
is denoted by $\theta$.
A filtered extension of
$(V,\theta)$ is defined to be
a filtered extension $\nbigp_{\ast}(\nbigv)$
over a meromorphic extension $\nbigv$
of $(V,\theta)$.
A filtered extension $(\nbigp_{\ast}\nbigv,\theta)$
is called regular
(resp. good, unramifiedly good)
if $(\nbigp_{\ast}\nbigv,\theta)$
is a regular (resp. good, unramifiedly good)
filtered Higgs bundle.
(See \cite[\S2.4]{Mochizuki-KH-Higgs}
for the notion of
good filtered Higgs bundles and
unramifiedly good filtered Higgs bundles.)

\begin{lem}
Let $f$ be the endomorphism of $V$
defined by $\theta=f\,dz/z$.
Let $a_j(z)$ be the holomorphic functions on $U^{\ast}$
obtained as the coefficients of
the characteristic polynomial
$\det(t\id_V-f)=\sum_{j=0}^r a_j(z)t^j$.
\begin{itemize}
 \item
      A meromorphic extension of $(V,\theta)$ exists
      if and only if
      the Higgs bundle $(V,\theta)$ is wild,
      i.e.,
      $a_j(z)$ are meromorphic at $z=0$.
      In that case,
      there exists a good filtered extension.
 \item A regular filtered extension exists
       if and only if $(V,\theta)$ is tame,
       i.e.,
       $a_j(z)$ are holomorphic at $z=0$.
\hfill\qed
\end{itemize}
\end{lem}

We obtain the Higgs field $\theta^{(\ell)}$ of $V^{(\ell)}$.
The following lemmas are clear.
\begin{lem}
The pull back and the descent induce
an equivalence between
meromorphic extensions of $(V,\theta)$
and $\Gal(\ell)$-equivariant meromorphic extensions of
$(V^{(\ell)},\theta^{(\ell)})$.
\hfill\qed
 \end{lem}

\begin{lem}
The pull back and the descent induce
an equivalence between
regular (resp. good) filtered 
extensions of $(V,\theta)$
and $\Gal(\ell)$-equivariant regular (resp. good) filtered
meromorphic extensions of
$(V^{(\ell)},\theta^{(\ell)})$.
\hfill\qed
\end{lem}

\subsection{Decomposable filtered extensions
of regular semisimple Higgs bundles}

\label{subsection;22.12.22.20}

\subsubsection{Decomposable filtered extensions}

We continue to use the notation in
\S\ref{subsection;22.12.22.10}.
Let $(V,\theta)$ be a regular semisimple Higgs bundle
on $U^{\ast}$. Assume that $\theta$ is wild.
There exist
$\ell\in\seisuu_{>0}$
and the decomposition
\begin{equation}
 \label{eq;22.12.22.1}
 \varphi_{\ell}^{\ast}(V,\theta)
  =\bigoplus_{i=1}^r
  (V_i,\theta_i),
\end{equation}
where $\rank V_i=1$,
and $\theta_i-\theta_j$ $(i\neq j)$ are nowhere vanishing on
$U^{(\ell)\ast}$.
Let $\nbigv$ be a meromorphic extension of
$(V,\theta)$.
The decomposition (\ref{eq;22.12.22.1})
extends to 
\begin{equation}
\label{eq;22.12.21.10}
 \varphi_{\ell}^{\ast}(\nbigv,\theta)
  =\bigoplus_{i=1}^r
  (\nbigv_i,\theta_i),
\end{equation}
where each $\nbigv_i$ is a meromorphic extension of $V_i$.

\begin{df}
\label{df;22.12.25.20}
A filtered bundle $\nbigp_{\ast}\nbigv$ over $\nbigv$
is called a decomposable filtered extension of $(V,\theta)$
if the filtered bundle $\varphi_{\ell}^{\ast}(\nbigp_{\ast}\nbigv)$
is compatible with the decomposition
{\rm(\ref{eq;22.12.21.10})},
i.e., the following holds for any $a\in\real$:
\[
 \nbigp_{a}(\varphi_{\ell}^{\ast}\nbigv)
 =\bigoplus_{i=1}^r
 \nbigp_a(\varphi_{\ell}^{\ast}\nbigv)
\cap \nbigv_{i}.
\]
Such $(\nbigp_{\ast}\nbigv,\theta)$ is called
a decomposable filtered Higgs bundle.
\hfill\qed
\end{df}

The following lemma is obvious by definition.
\begin{lem}
Suppose that $(\nbigp_{\ast}\nbigv,\theta)$ is decomposable,
\begin{itemize}
 \item
 $(\nbigp_{\ast}\nbigv,\theta)$ is a good filtered Higgs bundle.
 \item
Any decomposition
$(\nbigv,\theta)_{|U^{\ast}}
=(V_1,\theta_1)\oplus(V_2,\theta_2)$
extends to a decomposition
$(\nbigp_{\ast}\nbigv,\theta)
=(\nbigp_{\ast}\nbigv_1,\theta_1)\oplus
(\nbigp_{\ast}\nbigv_2,\theta_2)$
such that
$\nbigv_{i|U^{\ast}}=V_i$.      
\hfill\qed      
\end{itemize}
\end{lem}

\subsubsection{Filtered line bundles and
Decomposable filtered Higgs bundles}
\label{subsection;22.12.22.31}

There exists the decomposition
\begin{equation}
\label{eq;22.12.22.2}
 (V,\theta)=\bigoplus_{k\in S}(V^{[k]},\theta^{[k]})
\end{equation}
such that $\Sigma_{V^{[k]},\theta^{[k]}}$ are connected.
We set $r_k=\rank V^{[k]}$.
For each $k$, there exists the decomposition of
the Higgs bundle
\begin{equation}
 \label{eq;22.12.22.3}
 \varphi_{r_k}^{\ast}(V^{[k]},\theta^{[k]})
 =\bigoplus_{i=1}^{r_k}
 (V^{[k]}_i,\theta^{[k]}_i),
\end{equation}
where $\rank V^{[k]}_i=1$,
and $\theta^{[k]}_i$ are $1$-forms such that
$\theta^{[k]}_i-\theta^{[k]}_j$ $(i\neq j)$
are nowhere vanishing on $U^{(r_k)\ast}$.
A decomposable filtered extension
$\nbigp_{\ast}\nbigv$ of $(V,\theta)$
induces a decomposable filtered extension
$\nbigp_{\ast}(\nbigv^{[k]}_i)$
of $(V^{[k]}_i,\theta^{[k]}_i)$.
Note that
$\nbigp_{\ast}(\nbigv^{[k]}_i)
=\sigma^{\ast}\nbigp_{\ast}(\nbigv^{[k]}_1)$
for $\sigma\in\Gal(r_k)$
such that $\sigma^{\ast}\theta^{[k]}_1=\theta^{[k]}_i$.
Conversely,
a filtered extension
$\nbigp_{\ast}\nbigv^{[k]}_1$ of $V^{[k]}_1$
induces a $\Gal(r_k)$-equivariant filtered extension
$\bigoplus_{\sigma\in\Gal(r_k)} \sigma^{\ast}\nbigp_{\ast}\nbigv^{[k]}_1$
of $\varphi_{r_k}^{\ast}(V^{[k]})=\bigoplus_{i=1}^{r_k} V^{[k]}_i$,
and hence a decomposable filtered extension
$\nbigp_{\ast}\nbigv^{[k]}$ of
$(V^{[k]},\theta^{[k]})$.
Thus, we obtain
a decomposable filtered extension
$\bigoplus_{k\in S}\nbigp_{\ast}\nbigv^{[k]}$
of $(V,\theta)$.
Note that 
$\nbigp_{\ast}\nbigv^{[k]}$ is also obtained as
$(\varphi_{r_k})_{\ast}(\nbigp_{\ast}\nbigv^{[k]}_1)$
by the natural identification
$(\varphi_{r_k})_{\ast}(V^{[k]}_1)=V^{[k]}$.
The following proposition is easy to see.
\begin{prop}
\label{prop;22.12.22.4}
This procedure induces an equivalence
between decomposable filtered extensions of $(V,\theta)$
 and a tuple of filtered extensions of $V^{[k]}_1$
 $(k\in S)$.
\hfill\qed
\end{prop}

\subsubsection{Decomposable filtered extension determined by
determinant bundles}

Let $\nbigv$ be a meromorphic extension of $(V,\theta)$.
The decomposition (\ref{eq;22.12.22.2})
extends to a decomposition
\begin{equation}
\label{eq;22.12.22.22}
 (\nbigv,\theta)
  =\bigoplus_{k\in S}
 (\nbigv^{[k]},\theta^{[k]}).
\end{equation}
The decomposition (\ref{eq;22.12.21.10})
extends to a decomposition
\begin{equation}
 \label{eq;22.12.23.11}
 \varphi_{r_k}^{\ast}(\nbigv^{[k]},\theta^{[k]})
=\bigoplus_{i=1}^{r_k}
 (\nbigv^{[k]}_i,\theta^{[k]}_i).
\end{equation}

\begin{prop}
\label{prop;22.12.20.1}
For a tuple of filtered bundles
$\nbigp_{\ast}\det(\nbigv^{[k]})$ 
over $\det(\nbigv^{[k]})$,
there uniquely exists
a decomposable filtered bundle
$\nbigp^{\star}_{\ast}(\nbigv)
=\bigoplus_{k\in S}\nbigp^{\star}_{\ast}(\nbigv^{[k]})$
over $\nbigv$
such that
$\det(\nbigp^{\star}_{\ast}\nbigv^{[k]})
=\nbigp_{\ast}\det(\nbigv^{[k]})$
for any $k\in S$.
Moreover, the following holds for any $k\in S$.
\begin{itemize}
 \item $\dim \Gr^{\nbigp^{\star}}_a(\nbigv^{[k]})\leq 1$
       for any $a\in\real$.
 \item Let $d_k$ be a real number
       such that
       $\Gr^{\nbigp}_{d_k}(\det(\nbigv^{[k]}))\neq 0$.
       Then,
       $\Gr^{\nbigp^{\star}}_{a}(\nbigv^{[k]})\neq 0$
       if and only if
       $r_ka-d_k\in\seisuu$ ($r_k$ is odd),
       or
       $r_ka-d_k\in\frac{1}{2}\seisuu\setminus \seisuu$ ($r_k$ is even).
 \item $\Gr^{\nbigp^{\star}}_a(\nbigv^{[k]}_i)\neq 0$
       if and only if
       $a-d_k\in \seisuu$ ($r_k$ is odd),
       or
       $a-d_k\in \frac{1}{2}\seisuu$ ($r_k$ is even).
\end{itemize}
\end{prop}
\pf
It is enough to consider the case where
$\Sigma_{V,\theta}$ is connected,
i.e. $|S|=1$.
We omit the superscript $[k]$ and the subscript $k$.
We set
$(V^{(r)},\theta^{(r)})=\varphi_r^{\ast}(V,\theta)$
and $\nbigv^{(r)}=\varphi_{r}^{\ast}(\nbigv)$.
There exist the following decomposition
of the Higgs bundle on $U^{(r)\ast}$:
\begin{equation}
\label{eq;22.11.7.2}
 (V^{(r)},\theta^{(r)})
=\bigoplus_{i=1}^r (V_{\beta(i)},\beta(i)\,d\zeta).
\end{equation}
Here, $\beta(i)$ are meromorphic functions
on $(U^{(r)},0)$
such that
$\beta(i)-\beta(j)$ $(i\neq j)$
are nowhere vanishing on $U^{(r)\ast}$.
It extends to a decomposition on $U^{(r)}$:
\begin{equation}
\label{eq;22.11.7.1}
 (\nbigv^{(r)},\theta^{(r)})
 =\bigoplus_{i=1}^{r}
 (\nbigv_{\beta(i)},\beta(i)\,d\zeta).
\end{equation}
We have
$\sigma^{\ast}\nbigv_{\beta(i)}
=\nbigv_{\sigma^{\ast}(\beta(i))}$
for any $\sigma\in\Gal(r)$.

Let $v_{\beta(1)}$ be a frame of $\nbigv_{\beta(1)}$.
We obtain frames $v_{\sigma^{\ast}(\beta(1))}=\sigma^{\ast}v_{\beta(1)}$
of $\nbigv_{\sigma^{\ast}\beta(1)}$,
and the tuple $v_{\beta(1)},\ldots,v_{\beta(r)}$
is a frame of $\nbigv^{(r)}$.
We set
\[
b:=\min\Bigl\{
c\in\real\,\big|\,
  v_{\beta(1)}\wedge\cdots
 \wedge
 v_{\beta(r)}
\in\nbigp_c(\varphi_r^{\ast}\det\nbigv)
\Bigr\}.
\]
We define the filtered bundles
$\nbigp^{\star}_{\ast}(\nbigv_{\beta(i)})$
as follows:
\[
 \nbigp^{\star}_a(\nbigv_{\beta(i)})
 =\zeta^{-[a-b/r]}
 \nbigo_{U^{(r)}}v_{\beta(i)}. 
\]
They are independent of the choice of $v_{\beta(1)}$.
We set
$\nbigp^{\star}_{\ast}(\nbigv^{(r)})
=\bigoplus\nbigp^{\star}_{\ast}(\nbigv_{\beta(i)})$,
which is $\Gal(r)$-equivariant.
As the descent,
we obtain a filtered bundle
$\nbigp^{\star}_{\ast}(\nbigv)$ over $\nbigv$,
which satisfies the desired condition.
The uniqueness is clear.
By the construction,
$(\nbigp^{\star}_{\ast}(\nbigv),\theta)$
is clearly a good filtered Higgs bundle.

Let $\tau$ be a frame of $\nbigp_d(\det\nbigv)$.
There exist an integer $m$
and a nowhere vanishing holomorphic function $g$ on $U^{(r)}$
such that
\[
v_{\beta(1)}\wedge\cdots\wedge v_{\beta(r)}
=\zeta^mg(\zeta)\varphi_r^{\ast}\tau.
\]
Because a generator $\sigma_0$ of $\Gal(r)$
acts on the set $\{\sigma(i)\}$ in a cyclic way,
we have
$\sigma_0^{\ast}(v_{\beta(1)}\wedge\cdots\wedge v_{\beta(r)})
=(-1)^{(r-1)}v_{\beta(1)}\wedge\cdots\wedge v_{\beta(r)}$.
Hence, we obtain that
$\sigma_0^{\ast}(\zeta^m)=(-1)^{r-1}\zeta^m$
and $\sigma_0^{\ast}g=g$.
It implies that
$m/r\in\seisuu$ if $r$ is odd,
or
that
$m/r\in \frac{1}{2}\seisuu\setminus \seisuu$
if $r$ is even.
By our choice of $b$,
we have $b=-m+rd$.
It is easy to see that
$\Gr^{\nbigp^{\star}}_c(\nbigv_{\beta(i)})\neq 0$
if and only if $c-b/r\in\seisuu$.
For each $p\in\seisuu$,
we have the $\Gal(r)$-invariant sections
$\sum_{\sigma\in\Gal(r)}
\sigma^{\ast}(\zeta^pv_{\beta(1)})$
of $\nbigv^{(r)}$
which induces a section
of $\nbigp^{\star}_{b/r^2-p/r}(\nbigv)$.
Moreover, it induces a frame of
$\Gr^{\nbigp^{\star}}_{b/r^2-p/r}(\nbigv)$.
Hence, it is easy to see that
$\Gr^{\nbigp^{\star}}_a(\nbigv)\neq 0$
if and only if
$ra-b/r\in \seisuu$,
and that $\dim\Gr^{\nbigp^{\star}}_a(\nbigv)\leq 1$.
Then, we obtain the last two claims.
\hfill\qed

\subsection{Non-degenerate pairings and decomposable filtered extensions}

\subsubsection{Non-degenerate symmetric pairings of regular semisimple
Higgs bundles}
 
We continue to use the notation in \S\ref{subsection;22.12.22.20}.
Let $C$ be a non-degenerate symmetric pairing of $(V,\theta)$.
For any $z_0\in U^{\ast}$,
the eigen decomposition of $\theta$ at $z_0$ 
is orthogonal with respect to $C$.
The decomposition {\rm(\ref{eq;22.12.22.1})}
is orthogonal with respect to $\varphi_{\ell}^{\ast}C$.

The decomposition (\ref{eq;22.12.22.2})
is orthogonal with respect to $C$.
Let $C^{[k]}$ denote the restriction of $C$
to $V^{[k]}$.
The decomposition (\ref{eq;22.12.22.3})
is orthogonal with respect to
$\varphi_{r_k}^{\ast}C^{[k]}$.
Let $C^{[k]}_i$ denote the induced symmetric pairing
of $V^{[k]}_i$.
We have $C^{[k]}_i=\sigma^{\ast}C^{[k]}_1$
for $\sigma\in\Gal(r_k)$ such that
$\sigma^{\ast}\theta^{[k]}_1=\theta^{[k]}_i$.
Conversely,
for any non-degenerate symmetric pairings
$C^{[k]}_1$ $(k\in S)$,
we obtain a $\Gal(r_k)$-equivariant
non-degenerate symmetric pairing
$\bigoplus_{\sigma\in\Gal(r_k)} \sigma^{\ast}C^{[k]}_1$
of $\varphi_{r_k}^{\ast}V^{[k]}$.
It induces a non-degenerate symmetric pairing $C^{[k]}$
of $(V^{[k]},\theta^{[k]})$,
and a non-degenerate pairing
$\bigoplus C^{[k]}$ of $(V,\theta)$.
The following lemma is a special case of
Proposition \ref{prop;22.12.23.1}.
\begin{lem}
\label{lem;22.12.22.30}
These procedures induce an equivalence between
non-degenerate symmetric pairings $C$ of
$(V,\theta)$
and  
a tuple $(C^{[k]}_1)_{k\in S}$ of non-degenerate symmetric pairings
of $V^{[k]}_1$.
\hfill\qed 
\end{lem}

\subsubsection{Canonical decomposable filtered extensions
in the symmetric case}

We recall the following \cite[\S4.1]{Li-Mochizuki3}.

\begin{prop}
\label{prop;22.12.23.30}
For a non-degenerate symmetric pairing $C$ of $(V,\theta)$,
there uniquely exists a meromorphic extension
$\nbigv^C$ of $(V,\theta)$
compatible with $C$. 
Moreover, there uniquely exists a filtered bundle
$\nbigp^C_{\ast}(\nbigv^C)$
over $\nbigv^C$
satisfying the following conditions. 
\begin{itemize}
 \item $C$ is perfect with respect to $\nbigp^C_{\ast}(\nbigv^C)$.
 \item $\nbigp^C_{\ast}(\nbigv^C)$
       is a decomposable filtered extension of
       $(V,\theta)$.
       \hfill\qed
\end{itemize}
\end{prop}

We have the non-degenerate symmetric pairing $C^{[k]}_1$ $(k\in S)$
of $V^{[k]}_1$ corresponding to $C$
as in Lemma {\rm\ref{lem;22.12.22.30}}.
There exist unique filtered extensions
$\nbigp^C_{\ast}\bigl((\nbigv^{[k]}_1)^C\bigr)$ 
of $V^{[k]}_1$ compatible with $C^{[k]}_1$
as in Lemma {\rm\ref{lem;22.12.22.40}}.
The decomposable filtered extension 
$\nbigp^C_{\ast}(\nbigv^C)$ of $(V,\theta)$
corresponds to the tuple 
$\nbigp^C_{\ast}\bigl((\nbigv^{[k]}_1)^C\bigr)$ $(k\in S)$
(Proposition {\rm\ref{prop;22.12.22.4}}).

\subsubsection{Comparison of two canonical extensions}

Let $C$ be a non-degenerate symmetric pairing of $(V,\theta)$.
We have the unique filtered extension $\nbigp^C_{\ast}\nbigv^C$
of $(V,\theta)$ compatible with $C$.
We have the decomposition
\[
 (\nbigv^C,\theta)
 =\bigoplus_{k\in S}
 ((\nbigv^C)^{[k]},\theta^{[k]}).
\]
Let $\det(C^{[k]})$ denote the induced symmetric pairings of
$(\det(V^{[k]}),\tr(\theta^{[k]}))$.
Note that
$\det((\nbigv^C)^{[k]})$ is a meromorphic extension of
$(\det(V^{[k]}),\tr(\theta^{[k]}))$
compatible with $\det(C^{[k]})$.
We have the unique filtered extension
$\nbigp^C_{\ast}\det((\nbigv^C)^{[k]})$
of $(\det(V^{[k]}),\tr(\theta^{[k]}))$ compatible with $\det(C^{[k]})$.
We obtain the decomposable filtered Higgs bundle
$(\nbigp^{\star}_{\ast}(\nbigv^C),\theta)$
determined by
the tuple $\nbigp^C_{\ast}\det((\nbigv^C)^{[k]})$
as in Proposition \ref{prop;22.12.20.1}.

\begin{prop}
$\nbigp^C_{\ast}(\nbigv^C)=\nbigp^{\star}_{\ast}(\nbigv^C)$.
\end{prop}
\pf
The filtered Higgs bundle
$(\nbigp^C_{\ast}(\nbigv^C),\theta)$ is decomposable.
We have
$\det(\nbigp^C_{\ast}(\nbigv^C)^{[k]})
=\nbigp^C_{\ast}\det((\nbigv^C)^{[k]})
=\det\nbigp^{\star}_{\ast}((\nbigv^C)^{[k]})$.
Hence, we obtain
$\nbigp^C_{\ast}(\nbigv^C)=\nbigp^{\star}_{\ast}(\nbigv^C)$
by the uniqueness.
\hfill\qed

\begin{cor}
\label{cor;22.12.23.10}
Let $\nbigp_{\ast}(\nbigv^C)$ be a filtered extension of $(V,\theta)$
satisfying the following conditions.
\begin{itemize}
 \item $C$ is perfect with respect to $\nbigp_{\ast}(\nbigv^C)$.
 \item $\nbigp_{\ast}\nbigv^C=\bigoplus_{k\in S}
        \nbigp_{\ast}((\nbigv^C)^{[k]})$.
\end{itemize}
Let $\nbigp^{\star}_{\ast}(\nbigv)$
be the decomposable filtered extension of $(V,\theta)$
determined by
the filtered bundles 
$\det\bigl(\nbigp_{\ast}((\nbigv^C)^{[k]})\bigr)$
$(k\in S)$.
Then, $\nbigp^C_{\ast}(\nbigv^C)=\nbigp^{\star}_{\ast}(\nbigv^C)$. 
\end{cor}
\pf
It follows from
$\det\bigl(\nbigp_{\ast}((\nbigv^C)^{[k]})\bigr)
=\nbigp^C_{\ast}\det\bigl((\nbigv^C)^{[k]}\bigr)$.
\hfill\qed

\vspace{.1in}

Let $C$ and $C'$ be non-degenerate symmetric pairings of $(V,\theta)$.
Let $C^{[k]}$ and $C^{\prime[k]}$ $(k\in S)$ be the induced 
non-degenerate symmetric pairings of $(V^{[k]},\theta^{[k]})$.
We have the corresponding symmetric pairings
$C^{[k]}_1$ and $C^{\prime[k]}_1$
of $V^{[k]}_1$.

\begin{cor}
Suppose that $\det(C^{[k]})=\det(C^{\prime[k]})$ for any $k\in S$.
Then, $\nbigv^{C}=\nbigv^{C'}$ holds
if and only if
$\nbigp^C_{\ast}\nbigv^C=\nbigp^{C'}_{\ast}\nbigv^{C'}$ holds.
It is equivalent to the condition that
there exist holomorphic functions
$\gamma_1^{[k]}$  $(k\in S)$ on $U^{(r_k)}$ satisfying
$C^{\prime[k]}_1=\exp(\gamma_1^{[k]})C^{[k]}_1$
and
$\sum_{\sigma\in\Gal(r_k)}\sigma^{\ast}\gamma^{[k]}_1=0$.
\end{cor}
\pf
The ``if'' part of the claim is clear.
The ``only if'' part of the claim follows from
Corollary \ref{cor;22.12.23.10}.
\hfill\qed

\subsection{Prolongation of decoupled harmonic bundles}

Let $(V,\theta)$ be a Higgs bundle on $U^{\ast}$,
which is regular semisimple and wild.
Let $h$ be a decoupled harmonic metric of $(V,\theta)$.
We obtain the good filtered Higgs bundle
$(\nbigp^h_{\ast}V,\theta)$ on $(U,0)$.

\begin{lem}
\label{lem;22.12.23.22}
$(\nbigp^h_{\ast}V,\theta)$ is decomposable.
\end{lem}
\pf
Because the decomposition (\ref{eq;22.12.22.1})
is orthogonal with respect to $\varphi_{\ell}^{-1}(h)$,
the claim is clear.
\hfill\qed

\begin{rem}
If $h$ is a decoupled harmonic metric of $(V,\theta)$,
then we obtain that $\nbigp^h_{\ast}V$ is a filtered bundle
without assuming $\theta$ is wild.
\hfill\qed 
\end{rem}

We have the decomposition
$\nbigp^h_{\ast}(V)=\bigoplus_{k\in S}\nbigp^h_{\ast}(V^{[k]})$.
We obtain the filtered extensions
$\det(\nbigp^h_{\ast}V^{[k]})
=\nbigp^{\det(h)}_{\ast}\det(V^{[k]})$
of $\det(V^{[k]})$.
We have the filtered bundle
$\nbigp^{\star}_{\ast}(\nbigv)$
over $\nbigv=\nbigp^hV$
determined by $\det(\nbigp^h_{\ast}V^{[k]})$
as in Proposition \ref{prop;22.12.20.1}.
\begin{lem}
\label{lem;22.11.10.11}
We have $\nbigp^h_{\ast}(V)
=\nbigp^{\star}_{\ast}(\nbigv)$.
\end{lem}
\pf
This follows from the uniqueness of
the decomposable filtered extension
$\nbigp^{\star}_{\ast}(\nbigv)$ of $(V,\theta)$
satisfying the condition in Proposition \ref{prop;22.12.20.1}.
\hfill\qed

\vspace{.1in}

The decomposition (\ref{eq;22.12.22.2})
is orthogonal with respect to $h$.
Let $h^{[k]}$ denote the induced decoupled harmonic metric
of $(V^{[k]},\theta^{[k]})$ $(k\in S)$.
The decomposition (\ref{eq;22.12.22.3})
is orthogonal with respect to
$\varphi_{r_k}^{\ast}(h^{[k]})$.
Let $h^{[k]}_1$ denote the induced flat metric of $V^{[k]}_1$.

Let $h'$ be another decoupled harmonic metric of $(V,\theta)$.
Similarly, we obtain
the induced decomposable harmonic metric
$h^{\prime[k]}$ of $(V^{[k]},\theta^{[k]})$
and the induced flat metric
$h^{\prime[k]}_1$ of $V^{[k]}_1$.

\begin{cor}
\label{cor;22.12.25.10}
Suppose that $\det(h^{[k]})=\det(h^{\prime[k]})$
for any $k\in S$.
Then, $\nbigp^hV=\nbigp^{h'}V$ hold if and only if
$\nbigp^{h}_{\ast}(V)=\nbigp^{h'}_{\ast}(V)$ holds.
It is equivalent to the condition that
there uniquely exist holomorphic functions
$\gamma_1^{[k]}$ $(k\in S)$ on $U^{(r_k)}$
such that
(i) $h^{\prime[k]}_{1}=\exp(2\Re(\gamma_1^{[k]}))h^{[k]}_1$,
(ii) $\sum_{\sigma\in\Gal(r_k)} \sigma^{\ast}\gamma_1^{[k]}=0$.
\end{cor}
\pf
The ``if'' part of the claim is clear.
The ``only if'' part of the claim follows from
Lemma \ref{lem;22.11.10.11}.
The second claim is clear.
\hfill\qed

\subsection{Decoupled harmonic metrics and symmetric products}

\subsubsection{Comparison of extensions}

Let $(V,\theta)$ be a Higgs bundle on $U^{\ast}$
which is regular semisimple and wild.
Let $C$ be a non-degenerate symmetric pairing of $(V,\theta)$.
There exists a unique decoupled harmonic metric $h^C$
of $(V,\theta)$ compatible with $C$.

\begin{lem}
\label{lem;22.12.20.30}
We have $\nbigp^{h^C}_{\ast}(V)=\nbigp^C_{\ast}(\nbigv^C)$.
\end{lem}
\pf
By the pull back via $\varphi_{\ell}$,
it is enough to consider the case $\rank V=1$,
which is easy to check.
\hfill\qed

\subsubsection{Symmetric products compatible with a decoupled harmonic metric}

The following lemma is a special case of Lemma \ref{lem;22.12.23.20}.
\begin{lem}
\label{lem;22.12.23.12}
Suppose $\rank V=1$.
Let $h$ be a flat metric of $V$.
There exists a holomorphic non-degenerate symmetric product $C$
of $V$ which is compatible with $h$
if and only if
the monodromy of the Chern connection of $h$ is $1$ or $-1$.
It is equivalent to the condition
\[
 \bigl\{d\in\real\,\big|\,
 \Gr^{\nbigp^h}_d(V)\neq 0
 \bigr\}
\subset\frac{1}{2}\seisuu.
\]
If $C'$ is another non-degenerate symmetric pairing of
$V$ which is compatible with $h$,
there exists a non-zero constant $\alpha$
such that (i) $C'=\alpha C$,
(ii) $|\alpha|=1$.
\hfill\qed
 \end{lem}

\begin{prop}
\label{prop;22.12.20.2}
Let $h$ be a decoupled harmonic metric of $(V,\theta)$.
Suppose that there exist
non-degenerate symmetric products $C_{\det(V^{[k]})}$ $(k\in S)$
of $\det(V^{[k]})$
which are compatible with $\det(h^{[k]})$.
\begin{itemize}
\item 
There exists a non-degenerate symmetric pairing
$C$ of $(V,\theta)$ such that
(i)  $C$ is compatible with $h$, (ii) $\det(C^{[k]})=C_{\det(V^{[k]})}$.
\item If $C'$ is another non-degenerate symmetric pairing
       of $(V,\theta)$ satisfying the above conditions (i) and (ii).
       Then, there exist $r_k$-roots $\mu_k$ of $1$
       such that $C^{\prime[k]}=\mu_k C^k$.
\end{itemize}
\end{prop}
\pf
Let $h^{[k]}_1$ $(k\in S)$ be the induced flat metrics of $V^{[k]}_1$.
By Lemma \ref{lem;22.12.23.20} and Lemma \ref{lem;22.12.23.12},
there exist non-degenerate symmetric products
$C^{[k]}_1$ of $V^{[k]}_1$
compatible with $h^{[k]}_1$ for any $k\in S$.
They induce non-degenerate symmetric products
$C^{[k]}$ of $(V^{[k]},\theta^{[k]})$.
Because $\det(C^{[k]})$ is compatible with $\det h^{[k]}$,
there exist constants $\alpha_k$
such that
$\det(C^{[k]})=\alpha_k\cdot C_{\det(V^{[k]})}$
and $|\alpha_k|=1$.
By replacing $C^{[k]}$
with $\alpha_k^{1/r_k}C^{[k]}$,
we obtain the first claim.
The second claim is also clear.
\hfill\qed

\subsubsection{Existence}

Let $\nbigv$ be a meromorphic extension of $(V,\theta)$.
\begin{lem}
Let $C_{\det(V^{[k]})}$ be non-degenerate symmetric pairings of
$\det(V^{[k]})$
such that
$\det(\nbigv^{[k]})$ is compatible with $C_{\det(V^{[k]})}$.
Then, there exists a non-degenerate symmetric pairing $C$
of $(V,\theta)$ such that
(i) $\det(C^{[k]})=C_{\det(\nbigv^{[k]})}$,
(ii) $\nbigv^C=\nbigv$.
\end{lem}
\pf
It is enough to consider the case $|S|=1$.
We omit the superscript $[k]$ and the subscript $k$.
We use the notation in the proof of Proposition \ref{prop;22.12.20.1}.
Let $C_{1,\beta(1)}'$ be a non-degenerate symmetric pairing of
$\nbigv_{\beta(1)}$.
We obtain
a $\Gal(r)$-invariant non-degenerate symmetric pairing
$\bigoplus_{\sigma\in\Gal(r)}\sigma^{\ast}C_{1,\beta(1)}'$
of $\nbigv^{(r)}$.
It induces a non-degenerate symmetric pairing $C'$
of $\nbigv$.
From $C_{1,\beta(1)}''=\zeta C_{1,\beta(1)}'$,
we obtain another non-degenerate symmetric pairing $C''$,
for which we have
$\det(C'')=z\det(C')$.

Let $\alpha$ be the holomorphic function on $U^{\ast}$
determined by
$\det(C')=\alpha\cdot C_{\det(V)}$.
By the above consideration,
we may assume that
$\alpha$ induces a nowhere vanishing holomorphic function
on $U$.
By choosing an $r$-th root $\alpha^{1/r}$ of $\alpha$,
and by setting $C=\alpha^{-1/r}C$,
we obtain a desired non-degenerate pairing $C$.
\hfill\qed

\vspace{.1in}
We can prove the following lemma similarly.
\begin{lem}
Let $h_{\det(V^{[k]})}$ be flat metrics of $\det(V^{[k]})$
such that $\det(\nbigv^{[k]})=\nbigp^{h_{\det(V^{[k]})}}(\det(V^{[k]}))$.
There exists a decoupled harmonic metric $h$ of $(V,\theta)$
such that
(i) $\det(h^{[k]})=h_{\det(V^{[k]})}$,
(ii) $\nbigp^h(V)=\nbigv$.
\hfill\qed
\end{lem}

\subsection{Global case}

\subsubsection{Meromorphic extensions and filtered extensions}

Let $Y$ be a Riemann surface with a discrete subset $D$.
Let $\iota_{Y\setminus D}\colon Y\setminus D\to Y$
denote the inclusion.
For a holomorphic vector bundle $V$ on $Y\setminus D$,
a meromorphic extension of $V$ to $(Y,D)$
is defined to be a locally free $\nbigo_Y(\ast D)$-submodule $\nbigv$
of $(\iota_{Y\setminus D})_{\ast}V$
such that $\nbigv_{|Y\setminus D}=V$.
A filtered extension of $V$ to $(Y,D)$
is a filtered bundle $\nbigp_{\ast}\nbigv$
over a meromorphic extension $\nbigv$ of $V$.
We use similar terminology for non-degenerate symmetric parings
and Higgs bundles in this situation.

\subsubsection{Decomposable filtered extensions}
\label{subsection;22.12.23.40}

Let $(V,\theta)$ be a regular semisimple Higgs bundle
on $Y\setminus D$
which is wild along $D$.
Let $\nbigp_{\ast}\nbigv$ be a filtered extension
of $(V,\theta)$ to $(Y,D)$.

\begin{df}
$\nbigp_{\ast}\nbigv$ 
is called a decomposable filtered extension of $(V,\theta)$
if the restriction
to a neighbourhood of any $P\in D$
is decomposable.
\hfill\qed 
\end{df}
The following lemma is clear.
\begin{lem}
\label{lem;22.12.23.23}
A decomposable filtered Higgs bundle $(\nbigp_{\ast}\nbigv,\theta)$
is a good filtered Higgs bundle.
Any decomposition
$(\nbigv,\theta)_{|Y\setminus D}=
(V_1,\theta_1)\oplus(V_2,\theta_2)$
extends to a decomposition
$(\nbigp\nbigv,\theta)=
 (\nbigp_{\ast}\nbigv_1,\theta_1)
 \oplus(\nbigp_{\ast}\nbigv_2,\theta_2)$.
\hfill\qed 
\end{lem}

We have the line bundle $L_V$ on $\Sigma_{V,\theta}$
corresponding to $(V,\theta)$.
Let $\proj(T^{\ast}Y)$
be the projective completion of $T^{\ast}Y$.
Let $Z$ be the closure of
$\Sigma_{V,\theta}\subset T^{\ast}(Y\setminus D)$
in $\proj(T^{\ast}Y)$.
Let $\Sigmatilde_{V,\theta}\to Z$
denote the normalization.
We may naturally regard
$\Sigmatilde_{V,\theta}$
as a partial compactification of
$\Sigma_{V,\theta}$.
We set $\Dtilde=\Sigmatilde_{V,\theta}\setminus\Sigma_{V,\theta}$.
The morphism $\pi:\Sigma_{V,\theta}\to Y\setminus D$
uniquely extends to
a morphism $\pitilde:(\Sigmatilde_{V,\theta},\Dtilde)\to (Y,D)$.
From a meromorphic extension $\nbigl_V$ of $L_V$
to $(\Sigmatilde_{V,\theta},\Dtilde)$,
we obtain a meromorphic extension
$\pitilde_{\ast}(\nbigl_V)$ of $(V,\theta)$
to $(Y,D)$.
From a filtered extension $\nbigp_{\ast}\nbigl_V$
of $L_V$ to $(\Sigmatilde_{V,\theta},\Dtilde)$,
we obtain a decomposable filtered extension
$\pitilde_{\ast}(\nbigp_{\ast}\nbigl_V)$ of $(V,\theta)$
to $(Y,D)$.
The following proposition is a reformulation of
Proposition \ref{prop;22.12.22.4}.

\begin{prop}
The above procedure induce an equivalence between 
filtered extensions (resp. meromorphic extensions) of $L_V$
to $(\Sigmatilde_{V,\theta},\Dtilde)$
and 
decomposable filtered extensions (resp. meromorphic extensions)
of $(V,\theta)$ to $(Y,D)$.
\hfill\qed
\end{prop}

\subsubsection{Symmetric products}

Let $C$ be a non-degenerate symmetric pairing of $(V,\theta)$.
We restate Proposition \ref{prop;22.12.23.30}
in the global setting.

\begin{prop}
For a non-degenerate symmetric pairing $C$ of $(V,\theta)$,
there uniquely exists a meromorphic extension
$\nbigv^C$ of $(V,\theta)$ to $(Y,D)$
compatible with $C$.
Moreover, there uniquely exists a filtered bundle
$\nbigp^C_{\ast}(\nbigv^C)$
over $\nbigv^C$
satisfying the following conditions. 
\begin{itemize}
 \item $C$ is perfect with respect to $\nbigp^C_{\ast}(\nbigv^C)$.
 \item $\nbigp^C_{\ast}(\nbigv^C)$
       is a decomposable filtered extension of
       $(V,\theta)$.
       \hfill\qed
\end{itemize}
\end{prop}

The decomposable filtered extension
$\nbigp^C_{\ast}(\nbigv^C)$
is described as follows.
Let $C_0$ be the non-degenerate symmetric pairing of
$L_V$ corresponding to $C$.
There exists the unique filtered extension
$\nbigp^{C_0}_{\ast}(\nbigl^{C_0}_V)$ of $L_V$
to $(\Sigmatilde_{V,\theta},\Dtilde)$.
Then,
$\nbigp^C_{\ast}(\nbigv^C)=
\pitilde_{\ast}(\nbigp^{C_0}_{\ast}(\nbigl^{C_0}_V))$.

\subsubsection{Decoupled harmonic bundles}

Let $h$ be a decoupled harmonic metric of $(V,\theta)$.
We obtain the good filtered Higgs bundle
$(\nbigp^h_{\ast}V,\theta)$ on $(Y,D)$.
We obtain the following lemma from Lemma \ref{lem;22.12.23.22}.
\begin{lem}
\label{lem;22.12.21.12}
$(\nbigp^h_{\ast}V,\theta)$ is decomposable.
\hfill\qed
\end{lem}

We obtain the following lemma from Lemma \ref{lem;22.12.20.30}.
\begin{lem}
For a non-degenerate symmetric pairing $C$ of $(V,\theta)$,
we have $\nbigp^{h^C}_{\ast}(V)=\nbigp^C_{\ast}(\nbigv^C)$.
\hfill\qed
\end{lem}

\subsection{Kobayashi-Hitchin correspondence for decoupled harmonic bundles}

Let $X$ be a compact Riemann surface.
Let $D\subset X$ be a finite subset.
Let $(V,\theta)$ be a regular semisimple Higgs bundle
on $X\setminus D$,
which is wild along $D$.
For any decoupled harmonic metric $h$ of $(V,\theta)$
we obtain a good filtered Higgs bundle
$(\nbigp^h_{\ast}\nbigv,\theta)$ on $(X,D)$
which is polystable of degree $0$.
According to Lemma \ref{lem;22.12.21.12},
it is decomposable.

Conversely,
let $(\nbigp_{\ast}\nbigv,\theta)$
be a polystable decomposable filtered Higgs bundle of degree $0$ on $(X,D)$
such that $(V,\theta)=(\nbigv,\theta)_{|X\setminus D}$
is regular semisimple.
There exists a harmonic metric $h$ of $(V,\theta)$
adapted to $\nbigp_{\ast}\nbigv$
by \cite{Biquard-Boalch,Mochizuki-KH-Higgs,s2}.

\begin{prop}
\label{prop;22.12.21.20}
$h$ is a decoupled harmonic metric.
\end{prop}
\pf
It is enough to consider the case where
$(\nbigp_{\ast}\nbigv,\theta)$ is stable.
By Lemma \ref{lem;22.12.23.23},
$\Sigma_{V,\theta}$ is connected.
Let $\proj(T^{\ast}X)$ denote the projective completion of
$T^{\ast}X$.
Let $Z$ denote the closure of
$\Sigma_{V,\theta}$ in $\proj(T^{\ast}X)$.
Let $\Sigmatilde_{V,\theta}\to Z$
denote the normalization.
Let $\rho:\Sigmatilde_{V,\theta}\to X$
denote the induced morphism.
We set $\Dtilde=\rho^{-1}(D)$.
Let $L_V$ be the line bundle on $\Sigma_{V,\theta}$
corresponding to $(V,\theta)$.
Because $\nbigp_{\ast}\nbigv$ is a decomposable
filtered extension of $(V,\theta)$,
there exists the corresponding filtered extension
$\nbigp_{\ast}\nbigl_V$ of $L_V$
on $(\Sigmatilde_{V,\theta},\Dtilde)$.
We have $\rho_{\ast}(\nbigp_{\ast}\nbigl)=\nbigp_{\ast}\nbigv$.
By Proposition \ref{prop;22.12.21.21} below,
we have
$\deg(\nbigp_{\ast}\nbigl_V)
=\deg(\nbigp_{\ast}\nbigv)=0$.
There exists a flat metric $h_{L_V}$
of $L_V$ adapted to $\nbigp_{\ast}\nbigl_V$.
We obtain a decoupled harmonic metric $h_1$
of $(V,\theta)$
corresponding to $h_{L_V}$,
which is adapted to $\nbigp_{\ast}\nbigv$.
By the stability,
there exists a positive constant $h=ah_1$,
and hence $h$ is also a decoupled harmonic metric.
\hfill\qed

\subsubsection{Degree}

Let $\rho:X_1\to X_2$ be a non-constant morphism of
compact Riemann surfaces.
Let $D_2\subset X_2$ be a finite subset.
We set $D_1=\rho^{-1}(D_2)$.
Let $\nbigp_{\ast}\nbigv$ be a filtered bundle
on $(X_1,D_1)$.
We obtain a filtered bundle
$\rho_{\ast}(\nbigp_{\ast}\nbigv)$ on $(X_2,D_2)$.
Let $m(P)$ denote the ramification index of $\rho$ at $P\in X_1$.

\begin{prop}
\label{prop;22.12.21.21}
The following holds.  
\[
 \deg(\rho_{\ast}(\nbigp_{\ast}\nbigv))
=\deg(\nbigp_{\ast}\nbigv)
 -\frac{\rank\nbigv}{2}
 \sum_{P\in X_1\setminus D_1}
 (m(P)-1).
\]
\end{prop}
\pf
We have
$\nbigp_0(\rho_{\ast}\nbigv)=\rho_{\ast}(\nbigp_0\nbigv)$.
By the Grothendieck-Riemann-Roch theorem
and the Riemann-Hurwitz formula,
we have
\[
 \deg(\rho_{\ast}\nbigp_0\nbigv)
 =\deg(\nbigp_0\nbigv)
-\frac{\rank\nbigv}{2}\sum_{P\in X_1}(m(P)-1).
\]
By the construction of $\rho_{\ast}(\nbigp_{\ast}\nbigv)$,
we obtain
\begin{multline}
 \deg(\rho_{\ast}(\nbigp_{\ast}\nbigv))
=\deg(\rho_{\ast}(\nbigp_0\nbigv))
-\sum_{a\in D_1}
\sum_{-1<a\leq 0}
\sum_{j=0}^{m(P)-1}
\left(
\frac{a-j}{m(P)}
\right)
\dim\Gr^{\nbigp}_a(\nbigv_P)
 \\
=\deg(\nbigp_0\nbigv)
-\frac{\rank\nbigv}{2}\sum_{P\in X_1}(m(P)-1)
 -\sum_{P\in D_1}\sum_{-1<a\leq 0}
 \left(
 a-\frac{1}{2}(m(P)-1)
 \right)
 \dim\Gr^{\nbigp}_a(\nbigv_P)\\
 =\deg(\nbigp_{\ast}\nbigv)
 -\frac{\rank\nbigv}{2}\sum_{P\in X_1\setminus D_1}
 (m(P)-1).
\end{multline}
Thus, we are done.
\hfill\qed 

\begin{rem}
If there is no ramification point in $X_1\setminus D_1$,
we have $\deg(\nbigp_{\ast}\nbigv)=\deg(\rho_{\ast}\nbigp_{\ast}\nbigv)$.
We can also prove it as follows.
Let $h_0$  be a Hermitian metric of 
$\nbigv_{|X_1\setminus D_1}$
such that 
(i) $h_0$ is flat around any point of $D_1$,
(ii) $h_0$ is adapted to $\nbigp_{\ast}\nbigv$.
Let $R(h_0)$ be the curvature of the Chern connection of $h$.
Then, we have
\[
 \deg(\nbigp_{\ast}\nbigv)
=\frac{\sqrt{-1}}{2\pi}\int_{X_1\setminus D_1}\tr R(h_0).
\]
We have the induced metric $\rho_{\ast}(h_0)$
of $\rho_{\ast}(\nbigv)_{|X_2\setminus D_2}$.
It is flat around any point of $D_2$,
and it is adapted to $\rho_{\ast}(\nbigp_{\ast}\nbigv)$.
Hence, we have
\[
 \deg(\rho_{\ast}(\nbigp_{\ast}\nbigv))
=\frac{\sqrt{-1}}{2\pi}\int_{X_2\setminus D_2}\tr R(\rho_{\ast}h_0).
\]
Then, we obtain 
$\deg(\nbigp_{\ast}\nbigv)=\deg(\rho_{\ast}\nbigp_{\ast}\nbigv)$.
\hfill\qed
\end{rem}

\subsection{Dirichlet problem for wild decoupled harmonic bundles}

Let $Y$, $X$, $D$ and $(\nbigp_{\ast}\nbigv,\theta)$
be as in \S\ref{subsection;22.12.23.50}.
\begin{prop}
\label{prop;22.11.22.10}
Assume that $(V,\theta)$ is regular semisimple,
and that $\nbigp_{\ast}(\nbigv)$ is decomposable filtered extension.
Then, the harmonic metric $h$ in Theorem {\rm\ref{thm;22.12.20.4}}
is decoupled.
\end{prop}
\pf
It is enough to consider the case where
$\Sigma_{V,\theta}$ is connected.
Let $\Sigmatilde_{V,\theta}$ be the partial compactification
of $\Sigma_{V,\theta}$ as in \S\ref{subsection;22.12.23.40}.
Let $\Xtilde$ and $\Dtilde$ denote the inverse images of
$X$ and $D$ by the natural morphism
$\Sigmatilde_{V,\theta}\to Y$.
There exists the line bundle $L_V$ on $\Sigma_{V,\theta}$
corresponding to $(V,\theta)$.
Let $\nbigp_{\ast}\nbigl_V$ be the filtered line bundle
on $(\Sigmatilde_{V,\theta},\Dtilde)$
corresponding to $(\nbigp_{\ast}\nbigv,\theta)$.
There exists a Hermitian metric $h_0$ of $L_{V}$
such that
(i) $h_0$ is flat around any point of $\Dtilde$,
(ii) $h_0$ is adapted to $\nbigp_{\ast}\nbigl_V$,
(iii) $h_{0|\del \Xtilde}$ induces $h_{\del X}$.
Let $R(h_0)$ denote the curvature of the Chern connection of
$(L_V,h_0)$.
It vanishes around $\Dtilde$.
There exists a $\real$-valued $C^{\infty}$-function $\alpha$
on $\Xtilde$ such that
(i) $\delbar\del\alpha=R(h_0)_{|\Xtilde}$,
(ii) $\alpha_{|\del \Xtilde}=0$.
Then, $h_1=e^{-\alpha}h_0$ is a flat metric of
$L_{V|\Xtilde}$
adapted to $\nbigp_{\ast}\nbigl_V$
such that $h_{1|\del \Xtilde}=h_{0|\del \Xtilde}$.
Let $h_2$ be the decoupled harmonic metric of $(V,\theta)_{|X\setminus D}$
corresponding to $h_1$.
It is adapted to $\nbigp_{\ast}\nbigv$,
and it satisfies $h_{2|\del X}=h_{\del X}$.
By the uniqueness in Theorem \ref{thm;22.12.20.4},
we have $h=h_2$.
\hfill\qed

\section{Large-scale solutions with prescribed boundary value}

\subsection{Harmonic metrics of
regular semisimple Higgs bundles on a punctured disc}

\subsubsection{General case}

Let $U$ be a neighbourhood of $0$ in $\cnum$.
Let $U_0$ be a relatively compact open neighbourhood of $0$ in $U$
with smooth boundary $\del U_0$.
We set $U^{\ast}=U\setminus\{0\}$
and $U_0^{\ast}=U_0\setminus\{0\}$.

Let $(\nbigp_{\ast}\nbigv,\theta)$
be a good filtered Higgs bundle of rank $r$ on $(U,0)$
such that
$(V,\theta):=(\nbigv,\theta)_{|U^{\ast}}$
is regular semisimple.
Let $h_{\del U_0}$ be a Hermitian metric of
$V_{|\del U_0}$.
According to Theorem \ref{thm;22.12.20.4},
for any $t>0$,
there exists a unique harmonic metric $h_t$ of
$(V,t\theta)_{|U_0^{\ast}}$
such that
$h_{t|\del U_0}=h_{\del U_0}$
and that
$\nbigp^{h_t}_{\ast}(V)=\nbigp_{\ast}\nbigv$.
Note $\det(h_t)=\det(h_1)$ for any $t>0$.

\begin{prop}
\label{prop;22.11.7.12}
Let $t(i)$ be any sequence of positive numbers
such that $t(i)\to \infty$.
Then, there exists a subsequence $t'(j)$
such that the following holds.
\begin{itemize}
\item  $t'(j)\to \infty$.
\item  The sequence $h_{t'(j)}$ is convergent to
	a harmonic metric
       on any relatively compact open subset of $U_0^{\ast}$
       in the $C^{\infty}$-sense.
\end{itemize}
The limit $h_{\infty}$ is a decoupled harmonic metric of 
$(V,\theta)$
such that $\nbigp^{h_{\infty}}(V)=\nbigv$
and that $\det(h_{\infty})=\det(h_1)$.
\end{prop}
\pf
By taking the pull back via a ramified covering
map $\varphi_{\ell}$ as in \S\ref{subsection;22.12.22.10},
it is enough to consider the case where
there exist meromorphic functions
$\gamma(1),\ldots,\gamma(r)$ on $(U,0)$
and a decomposition
\[
(\nbigv,\theta)
=\bigoplus_{i=1}^r
(\nbigv_{i},\gamma(i)\,dz).
\]
Let $v_i$ be a frame of $\nbigv_i$ on $U$
such that $v_i$ is a section of $\nbigp_{<0}\nbigv$.
\begin{lem}
\label{lem;22.11.7.10}
There exists a constant $C>0$
such that
$h_{t}(v_i,v_i)\leq C$ for any $t>0$. 
\end{lem}
\pf
It is enough to consider the case where $\gamma(i)=0$.
We have $\theta(v_i)=0$.
Then, we have
$-\del_z\del_{\zbar}|v_i|_{h_t}^2\leq 0$
on $U_0^{\ast}$
(see a preliminary Weitzenb\"{o}ck formula
in \cite[Proof of Lemma 4.1]{s2}).
Because $v_i$ is a section of $\nbigp_{<0}\nbigv$,
$|v_i|_{h_t}^2$ is bounded for each $t$.
Hence, $|v_i|_{h_t}^2$ is subharmonic on $U_0$.
By the maximum principle,
we obtain
$|v_i|^2_{h_t}
\leq
\max_{\del U_0}|v_i|_{h_t}^2
=\max_{\del U_0}|v_i|_{h_1}^2$.
\hfill\qed

\vspace{.1in}

Let $\nbigv^{\lor}=\nhom_{\nbigo_U}(\nbigv,\nbigo_U(\ast 0))$
denote the dual of $\nbigv$.
We have the induced filtered bundle
$\nbigp_{\ast}(\nbigv^{\lor})$ on $\nbigv^{\lor}$.
We set
$(V^{\lor},\theta^{\lor})
=(\nbigv^{\lor},\theta^{\lor})_{|U^{\ast}}$.
The induced harmonic metric $h_t^{\lor}$
of $(V^{\lor},t\theta^{\lor})$
is adapted to $\nbigp_{\ast}(\nbigv^{\lor})$.

There exists the induced decomposition
$\nbigv^{\lor}
 =\bigoplus_{i=1}^r\nbigv_i^{\lor}$.
Let $v_i^{\lor}$ denote the section of $\nbigv^{\lor}_i$
such that $v_i^{\lor}(v_i)=1$.
There exists $m(i)\in\seisuu_{>0}$
such that
$z^{m(i)}v_i^{\lor}$ is a section of $\nbigp_{<0}(\nbigv^{\lor})$.
By Lemma \ref{lem;22.11.7.10},
we obtain the following lemma.

\begin{lem}
\label{lem;22.11.7.11}
There exists $C>0$ such that
$|z|^{2m(i)}h_t^{\lor}(v_i^{\lor},v_i^{\lor})\leq C$
for any $t>0$.
\hfill\qed
\end{lem}

Let $s_t$ be the automorphism of
$V_{|U_0^{\ast}}$
determined by $h_t=h_1\cdot s_t$.
Let $K$ be any relatively compact open subset of $U_0^{\ast}$.
By Lemma \ref{lem;22.11.7.10} and Lemma \ref{lem;22.11.7.11},
there exist $C_{K,1}>0$
such that the following holds for any $t>0$:
\begin{equation}
\label{eq;22.12.20.5}
 |s_t|_{h_1}+|s_t^{-1}|_{h_1}
 \leq
 C_{K,1}.
\end{equation}
By a variant of Simpson's main estimate
(see \cite[Proposition 2.3]{Decouple}),
there exist $t_{K,1},C_{K,2},C_{K,3}>0$
such that the following holds for any $t>t_{K,1}$
and for any local sections $u_{\beta(i)}$
and $u_{\beta(j)}$ of $V_{\beta(i)}$ and $V_{\beta(j)}$ on $K$
$(i\neq j)$:
\begin{equation}
 \label{eq;22.12.20.7}
  \bigl|h_t(u_{\beta(i)},u_{\beta(j)})\bigr|
  \leq
  C_{K,2}
  \exp(-C_{K,3}t)
  |u_{\beta(i)}|_{h_t}\cdot
  |u_{\beta(j)}|_{h_t}.
\end{equation}
There also exist $t_{K,2},C_{K,4},C_{K,5}>0$
such that the following holds on $K$ for any $t>t_{K,2}$
(see \cite[Theorem 2.9]{Decouple}):
\begin{equation}
\label{eq;22.12.20.6}
 \bigl|
 R(h_{t})
 \bigr|_{h_1}
 \leq
 C_{K,4}
 \exp\Bigl(
 -C_{K,5}t
 \Bigr).
\end{equation}
By (\ref{eq;22.12.20.5}) and (\ref{eq;22.12.20.6}),
it is standard to obtain the existence of
a convergent subsequence $h_{t'(j)}$.
By (\ref{eq;22.12.20.7}) and (\ref{eq;22.12.20.6}),
the limit is a decoupled harmonic metric.
By Lemma \ref{lem;22.11.7.10},
we obtain that $h_{\infty}(v_i,v_i)\leq C$.
Hence, $v_i$ are sections of
$\nbigp^{h_{\infty}}(V)$.
It implies that
$\nbigv\subset \nbigp^{h_{\infty}}(V)$.
Because both
$\nbigv$ and $\nbigp^{h_{\infty}}(V)$
are locally free $\nbigo_{U}(\ast 0)$-modules,
we obtain that 
$\nbigv=\nbigp^{h_{\infty}}(V)$.
\hfill\qed

\begin{prop}
\label{prop;22.12.20.11}
Let $h_{\infty}$ denote the limit
of a convergent subsequence
in Proposition {\rm\ref{prop;22.11.7.12}}.
Suppose the following condition.
\begin{itemize}
 \item For every $z_0\in \del U_0$,
       the eigen decomposition of $\theta$
       at $z_0$
       is orthogonal with respect to
       $h_{\del U_0}$.
\end{itemize}
Then, $h_{\infty|\del U_0}=h_{\del U_0}$.
\end{prop}
\pf
Let $U_1$ be a relatively compact open neighbourhood of $0$
in $U_0$ with smooth boundary $\del U_1$.
Because $h_{\infty}$ is a decoupled harmonic metric,
the following condition is satisfied.
\begin{itemize}
 \item For every $z_1\in \del U_1$,
       the eigen decomposition of $\theta$
       at $z_1$
       is orthogonal with respect to
       $h_{\infty}$.
\end{itemize}
We set $A=U_0\setminus \overline{U_1}$.
By Proposition \ref{prop;22.11.22.10},
there exists a decoupled harmonic metric
$h^{(1)}$
of
$(V,\delbar_V,\theta)_{|A}$
such that
$h^{(1)}_{|\del U_0}=h_{\del U_0}$
and
$h^{(1)}_{|\del U_1}=h_{\infty|\del U_1}$.
We note that
$h^{(1)}$ is a harmonic metric of
$(V,\delbar_V,t\theta)_{|A}$ for any $t>0$.
We also note that
$\det(h^{(1)})=\det(h_1)_{|A}$
because
$\det(h^{(1)})_{|\del A}=\det(h_1)_{|\del A}$.

Let $s_t$ be determined by
$h_t=h^{(1)}s_t$ on $A$.
We have
$-\del_z\del_{\zbar}\Tr(s_t)\leq 0$.
We have
$s_{t'(j)}\to \id$ on $\del U_1$
and $s_{t'(j)}=\id$ on $\del U_0$.
Hence, we obtain
$\bigl|
 \Tr(s_{t'(j)}-\id)
\bigr|\to 0$ as $t'(j)\to\infty$.
It implies the claim of the proposition.
\hfill\qed

\subsubsection{The irreducible case}

Suppose that the spectral curve is irreducible,
i.e., $\Sigma_{V,\theta}$ is connected.
We obtain the decomposable filtered bundle
$\nbigp^{\star}_{\ast}(\nbigv)$
determined by
$\det(\nbigp_{\ast}\nbigv)$
as in Proposition \ref{prop;22.12.20.1},
which is not necessarily equal to $\nbigp_{\ast}(\nbigv)$.

\begin{lem}
\label{lem;22.12.20.10}
Let $h_{\infty}$ denote the limit of a convergent subsequence
in Proposition {\rm\ref{prop;22.11.7.12}}.
Then, we have
$\nbigp^{h_{\infty}}_{\ast}(V)
=\nbigp^{\star}_{\ast}(\nbigv)$. 
\end{lem}
\pf
We have $\nbigp^{h_{\infty}}V=\nbigv$.
Because $h_{\infty}$ is a decoupled harmonic metric,
$\nbigp^{h_{\infty}}_{\ast}(\nbigv)$ is decomposable.
Because $\det(h_{\infty})=\det(h_1)$,
we obtain
$\det(\nbigp^{h_{\infty}}_{\ast}V)=\det(\nbigp_{\ast}\nbigv)$.
Then, the claim follows from the uniqueness of
$\nbigp^{\star}_{\ast}(\nbigv)$.
\hfill\qed

\vspace{.1in}

Let $h_0$ be any decoupled harmonic metric
of $(V,\delbar_V,\theta)$
such that $\nbigp^{h_0}(V)=\nbigv$
and that $\det(h_0)$ is adapted to
$\det(\nbigp_{\ast}\nbigv)$.
By the argument in the proof of Lemma \ref{lem;22.12.20.10},
we can prove 
$\nbigp^{h_0}_{\ast}(V)=\nbigp^{\star}_{\ast}(\nbigv)$.
Let $h_t$ $(t>0)$ be the harmonic metrics of
$(V,\delbar_V,t\theta)$
adapted to $\nbigp_{\ast}\nbigv$
such that $h_{t|\del U_0}=h_{0|\del U_0}$.

\begin{prop}
\label{prop;22.11.30.131}
The sequence $h_t$ is convergent to $h_0$ as $t\to\infty$
in the $C^{\infty}$-sense
on any relatively compact open subset of $U_0^{\ast}$.
\end{prop}
\pf
Let $t_i$ be any subsequence such that $t_i\to\infty$
and that $h_{t_i}$ is convergent.
Let $h_{\infty}$ denote the limit.
By Proposition \ref{prop;22.12.20.11},
we have
$h_{\infty|\del U_0}=h_{0|\del U_0}$.
We also have
$\nbigp_{\ast}^{h_{\infty}}(V)
=\nbigp^{\star}_{\ast}(\nbigv)
=\nbigp_{\ast}^{h_0}(V)$.
Hence, we obtain $h_{\infty}=h_0$.
It implies that
$h_t$ is convergent to $h_0$ as $t\to\infty$.
\hfill\qed

\subsubsection{Symmetric case}

We do not assume that the spectral curve is irreducible.
Instead, suppose that there exists
a perfect pairing $C$ of $(\nbigp_{\ast}\nbigv,\theta)$.
There uniquely exists a decoupled harmonic metric $h^C$
of $(V,\theta)$ which is compatible with $C$.
As in Lemma \ref{lem;22.12.20.30},
we have
$\nbigp^{h^C}_{\ast}(V)=
\nbigp^C_{\ast}\nbigv$.

Suppose that $h_{\del U_0}$
is compatible with $C_{|\del U_0}$.
Then,
$h_t$ $(t>0)$ are compatible with $C$
by Corollary \ref{cor;22.12.20.40}.
Let $s_t$ be determined by $h_t=h^Cs_t$.
We note that
$\det(h_t)=\det(h_1)=\det(h^C)$
by the compatibility with $C$.
The following proposition is a special case of
Corollary \ref{cor;22.12.20.31}.

\begin{prop}
\label{prop;22.11.30.202}
If $h_{\del U_0}$ is compatible with $C_{|\del U_0}$,
the sequence $h_t$ is convergent to $h^C$
in the $C^{\infty}$-sense
on any relatively compact subset $K$ of $U_0^{\ast}$.
Moreover,
there exists $t(K)>0$ such that the following holds
for any $\ell\geq 0$:
\begin{itemize}
 \item 
There exists $C(K,\ell)$ and $\epsilon(K,\ell)$
such that the norms of
$s_t-\id$ $(t\geq t(K))$
and their derivatives up to order $\ell$
are dominated by 
$C(K,\ell)\exp(-\epsilon(K,\ell)t)$.
\hfill\qed
\end{itemize}
\end{prop}

Let us also consider the case where 
$h_{\del U_0}$ is not necessarily compatible with $C_{|\del U_0}$,
but $\det(h_{\del U_0})$ is compatible with $\det(C)_{|\del U_0}$.
Because $\det(h_t)$ are compatible with $\det(C)$ on $U_0$,
we obtain $\det(h_t)=\det(h_1)=\det(h^C)$.
\begin{prop}
Let $h_{t(i)}$ be a convergent subsequence,
and $h_{\infty}$ denote the limit
as in Proposition {\rm\ref{prop;22.11.7.12}}.
Then, $\nbigp^{h_{\infty}}_{\ast}(V)=\nbigp^C_{\ast}(\nbigv)$. 
\end{prop}
\pf
Let $h'_t$ $(t>0)$ be harmonic metrics of $(V,t\theta)$
which are compatible with $C$,
such that $\det(h'_t)=\det(h_1)$.
We have already proved that
the sequence $h'_t$ is convergent to $h^C$.
We have $\det(h'_t)=\det(h_t)$.
Let $s_t$ be the automorphism determined by
$h_t=h_t's_t$.
Let $s_{\infty}$ be determined by
$h_{\infty}=h^Cs_{\infty}$.
The sequence $s_t$ is convergent to $s_{\infty}$.
Because $\det(s_t)=1$,
we have $\det(s_{\infty})=1$.
Because $\Tr(s_t)$ is subharmonic on $U_0$,
we obtain that
$\max_{U_0}\Tr(s_t)=
\max_{\del U_0}\Tr(s_t)
=\max_{\del U_0}\Tr(s_1)$.
We obtain that $\Tr(s_{\infty})$ is bounded.
Then, $s_{\infty}$ and $s_{\infty}^{-1}$
are bounded,
and we obtain
$\nbigp^{h_{\infty}}(V)=\nbigp^{C}_{\ast}(\nbigv)$.
\hfill\qed

\vspace{.1in}

Suppose that for every $z_0\in \del U_0$
the eigen decomposition of $\theta$
is orthogonal with respect to $h_{\del U_0}$.
There exists a decoupled harmonic metric
$\htilde$ of $(V,\theta)$
such that $\htilde_{|\del U_0}=h_{\del U_0}$
and $\nbigp^{\htilde}_{\ast}(V)=\nbigp^C_{\ast}(\nbigv)$.
\begin{cor}
The sequence $h_{t}$ is convergent to $\htilde$.
\hfill\qed
\end{cor}

\subsection{Local symmetrizability of Higgs bundles}

Let $U$ be a simply connected open subset in $\cnum$.
Let $D$ be a finite subset of $U$.
Let $(E,\delbar_E,\theta)$ be a Higgs bundle on $U$
such that $(V,\theta)=(E,\theta)_{|U\setminus D}$
is regular semisimple.
Let $\pi:\Sigma_{E,\theta}\to U$
denote the projection.
Let $\rho:\Sigmatilde_{E,\theta}\to\Sigma_{E,\theta}$
denote the normalization of
$\Sigma_{E,\theta}$.
We set $\Dtilde=(\pi\circ\rho)^{-1}(D)$.
We assume the following condition.
\begin{itemize}
 \item There exists a line bundle $L$
       on $\Sigmatilde_{E,\theta}$
       with an isomorphism
       $(\pi\circ\rho)_{\ast}L\simeq E$.
       Moreover,
       the Higgs field $\theta$ of $E$
       is induced by
       the $\nbigo_{T^{\ast}U}$-action
       on $\rho_{\ast}L$.
\end{itemize}

For any $P\in D$,
let $U_P$ be a simply connected neighbourhood of $P$ in $U$
such that $U_P\cap D=\{P\}$.
We set $U_P^{\ast}=U_P\setminus\{P\}$.
There exists the decomposition
\begin{equation}
\label{eq;23.1.28.2}
(V,\theta)_{|U_P^{\ast}}
=\bigoplus_{k\in S(P)}(V^{[k]}_{P},\theta^{[k]}_{P})
\end{equation}
such that
the spectral curves of $(V_P^{[k]},\theta_P^{[k]})$ are connected.
Because $E\simeq (\pi\circ\rho)_{\ast}L$,
(\ref{eq;23.1.28.2}) extends to the decomposition
\[
 (E,\theta)_{|U_P^{\ast}}
=\bigoplus_{i\in S(P)}(E^{[k]}_{P},\theta^{[k]}_{P}).
\]

Let $h$ be a decoupled harmonic metric 
of $(V,\theta)$.
The decomposition (\ref{eq;23.1.28.2}) is orthogonal
with respect to $h$.
Let $h_P^{[k]}$ denote the restriction of $h$
to $V^{[k]}_P$.
We consider the following condition.
\begin{condition}
\label{condition;23.1.28.32}
$\det(h^{[k]}_P)$ induces a flat metric of $\det(E^{[k]}_P)$,
and $\nbigp^hV=E(\ast D)$ holds.
\end{condition}

\begin{prop}
\label{prop;22.12.23.110}
Suppose that 
Condition {\rm\ref{condition;23.1.28.32}} is satisfied
at each $P\in D$.
Moreover, we assume that
each connected component of
$\Sigmatilde_{E,\theta}$ is simply connected.
Then, the following claims hold.
\begin{itemize}
 \item  There exists a non-degenerate symmetric pairing 
$C$ of $(E,\theta)$
such that $C_{|U\setminus D}$ is compatible with $h$.
 \item Let $C'$ be a non-degenerate symmetric pairing of
       $(V,\theta)$ which is compatible with $h$.
       Then, $C'$ induces a non-degenerate symmetric pairing
       of $E$.
 \end{itemize}
\end{prop}

\begin{rem}
If $\Sigma_{E,\theta}$ is
a simply connected complex submanifold
of $T^{\ast}U$,
we can apply Proposition {\rm\ref{prop;22.12.23.110}}
to $(E,\theta)$.

\hfill\qed
\end{rem}

\subsubsection{Special case}
\label{subsection;22.12.23.111}

Let us study the case that $D=\{0\}$,
and that $\Sigma_{V,\theta}$ is connected.
We set $\nbigv=E(\ast 0)$.
We use the notation in \S\ref{subsection;22.12.22.10}.
By choosing an $r$-th root of
$(\pi\circ\rho)^{\ast}(z)$ on $\Sigmatilde_{E,\theta}$,
we obtain
a holomorphic isomorphism
$\psi:\Sigmatilde_{E,\theta}\to U^{(r)}$
such that
$\varphi_r\circ\psi=\pi\circ\rho$.
There exists the decomposition (\ref{eq;22.11.7.1}) on $U^{(r)}$.
There exists the natural isomorphism
$\psi_{\ast}(L)(\ast 0)\simeq\nbigv_{\beta(1)}$.
Let $E_{\beta(1)}\subset\nbigv_{\beta(1)}$
denote the image of $L$.
We have $\varphi_{r\ast}(E_{\beta(1)})=E$.

Let
$C_{\beta(1)}:
\nbigv_{\beta(1)}\otimes \nbigv_{\beta(1)}
\lrarr\nbigo_{U^{(r)}}(\ast 0)$
be a non-degenerate symmetric pairing.
There exists the morphism
$\tr:\varphi_{r\ast}\nbigo_{U^{(r)}}(\ast 0)\to \nbigo_{U}(\ast 0)$
as in \S\ref{subsection;22.12.23.100}.
We obtain the induced symmetric pairing
$\Psi(C_{\beta(1)})=\tr\circ\varphi_{r\ast}(C_{\beta(1)})$
of $\nbigv=\varphi_{r\ast}(\nbigv_{\beta(1)})$.
There exists an integer $k$
such that
$C_{\beta(1)}(E_{\beta(1)}\otimes E_{\beta(1)})
=\nbigo_{U^{(r)}}(k\{0\})$.

\begin{lem}
\label{lem;22.12.24.1}
$\Psi(C_{\beta(1)})$
induces a symmetric pairing of $E$
if and only if $k\leq r-1$.
The induced pairing is non-degenerate if and only if
$k=r-1$.
\end{lem}
\pf
There exists a generator $v$ of $E_{\beta(1)}$
such that $C_{\beta(1)}(v,v)=\zeta^{-k}$.
The tuple $v,\zeta v,\ldots,\zeta^{r-1}v$
induces a frame of $E$.
Note that $\tr(\zeta^j)=0$ unless $j\in r\seisuu$.
It is easy to see that
$\tr\bigl(C_{\beta(1)}(\zeta^iv,\zeta^jv)\bigr)
=\tr(\zeta^{i+j-k})$ $(0\leq i,j\leq r-1)$
are holomorphic at $0$
if and only if $k\leq r-1$,
and that the induced pairing is non-degenerate at $0$
if and only if $k=r-1$.
\hfill\qed

\vspace{.1in}

Let $C_{0,\beta(1)}$ be a non-degenerate symmetric pairing of
$\nbigv_{\beta(1)}$
such that
$C_{0,\beta(1)}(E_{\beta(1)}\otimes E_{\beta(1)})
=\nbigo_{U^{(r)}}((r-1)\{0\})$.
We set
$C_{0}=\Psi(C_{0,\beta(1)})$
which is a non-degenerate symmetric pairing
of $(E,\theta)$.
Let $h_0$ be a decoupled harmonic metric
of $(V,\theta)$ compatible with $C_0$.
We note that $\det(h_0)$ is compatible with
$\det(C_0)$,
and hence it induces a Hermitian metric of $\det(E)$.

\vspace{.1in}

Let $h_1$ be any decoupled harmonic metric of $(V,\theta)$
such that 
$\nbigp^{h_1}(V)=\nbigv$
and that $\det(h_1)=\det(h_0)$.
By Corollary \ref{cor;22.12.25.10},
there exists a holomorphic function $\gamma_1$ on $U^{(r)}$
such that
(i)
$\varphi_r^{\ast}(h_1)_{|V_{\beta(1)}}
=\exp(2\Re\gamma_1)
\varphi_r^{\ast}(h_0)_{|V_{\beta(1)}}$,
(ii) $\sum_{\sigma\in\Gal(r)} \sigma^{\ast}\gamma_1=0$.
We set
\[
 C_{1,\beta(1)}=
\exp(2\gamma_1)C_{0,\beta(1)}.
\]
It is a non-degenerate symmetric pairing
of $\nbigv_{\beta(1)}$
satisfying
$C_{1,\beta(1)}(E_{\beta(1)}\otimes E_{\beta(1)})
=\nbigo_{U^{(r)}}((r-1)\{0\})$.
We obtain a non-degenerate symmetric pairing
$C_1=\Psi(C_{1,\beta(1)})$ of $(E,\theta)$
such that
$C_{1|U^{\ast}}$ is compatible with $h_1$.

Let $h$ be any decoupled harmonic metric of $(V,\theta)$
such that $\nbigp^h(V)=\nbigv$
and that $\det(h)$ induces a flat metric of $\det(E)$.
There exists a holomorphic function $\gamma_2$ on $U$
such that
$\det(h)=\exp(2\Re(\gamma_2))\det(h_1)$.
Then, $C=\exp(2\gamma_2)C_1$ is compatible with $h$,
and it induces a non-degenerate symmetric pairing of
$E$.

\begin{lem}
\label{lem;22.12.23.112}
Let $C'$ be a non-degenerate symmetric pairing
of $(V,\theta)$ compatible with $h$. 
Then, $C'$ induces a non-degenerate symmetric pairing of
$E$.
\end{lem}
\pf
There exist non-degenerate
symmetric pairings
$C_{\beta(1)}$ and $C'_{\beta(1)}$ of $\nbigv_{\beta(1)}$
such that
$\Psi(C_{\beta(1)})=C$
and 
$\Psi(C'_{\beta(1)})=C'$,
respectively.
Because both $C_{\beta(1)}$ and $C'_{\beta(1)}$
are compatible with
$\varphi_r^{\ast}(h)_{|V_{\beta(1)}}$,
there exists a constant $\alpha$ such that $|\alpha|=1$
such that
$C_{\beta(1)}'=\alpha C_{\beta(1)}$.
Hence,
$C'_{\beta(1)}(E_{\beta(1)}\otimes E_{\beta(1)})
=\nbigo_{U^{(r)}}((r-1)\{0\})$,
and hence $C'$ induces a non-degenerate symmetric pairing
of $E$.
\hfill\qed

\subsubsection{Proof of Proposition \ref{prop;22.12.23.110}}

It is enough to consider the case where
$\Sigma_{V,\theta}$ is connected,
which implies that $\Sigmatilde_{E,\theta}$ is connected.
Let $h_{L}$ denote the flat metric of
$L_{|\Sigma_{V,\theta}}$
corresponding to the decoupled harmonic metric $h$.
Let $P$ be any point of $D$.
By Proposition \ref{prop;22.12.20.2},
there exists a non-degenerate symmetric pairing of
$V_{|U_P^{\ast}}$
which is compatible with $h_{|U_P^{\ast}}$.
There exists a non-degenerate symmetric pairing of
$L$ on $(\pi\circ\rho)^{-1}(U_P^{\ast})$
which is compatible with $h_L$.
Hence, the monodromy of the Chern connection of $h_{L}$
around any point of $\Dtilde$ are $1$ or $-1$.
Because $\Sigmatilde_{E,\theta}$ is simply connected,
Lemma \ref{lem;22.12.23.20} implies that
there exists a non-degenerate symmetric pairing
$C_{L}$ of $L_{|\Sigma_{V,\theta}}$
compatible with $h_{L}$.
It induces a non-degenerate symmetric pairing $C$
of $(V,\theta)$ compatible with $h$.
By Lemma \ref{lem;22.12.23.112},
$C$ induces a non-degenerate symmetric pairing of $E$.
Thus, we obtain the first claim of Proposition \ref{prop;22.12.23.110}.
The second claim also follows from
Lemma \ref{lem;22.12.23.112}.
\hfill\qed

\subsection{A uniform estimate in the symmetric case}

\subsubsection{Setting}

For $R>0$,
we set $B(R)=\bigl\{z\in\cnum\,\big|\,|z|<R\bigr\}$.
Let $\nbigs\subset\cnum^n$ be a connected open subset
with a base point $x_0$.
Let $\nbigz_i$ $(i=1,2)$ be an open subset of
$\nbigs\times\cnum_{z_i}$.
For simplicity,
we assume that $\nbigz_2=\nbigs\times B(2)$.
Let $p_i:\nbigz_i\to \nbigs$ denote the projections.
We set $T^{\ast}(\nbigz_2/\nbigs)=\nbigs\times T^{\ast}B(2)$.
Let $\pi_2:T^{\ast}(\nbigz_2/\nbigs)\to\nbigz_2$
denote the projection.
Let $\Phi_0:\nbigz_1\to T^{\ast}(\nbigz_2/\nbigs)$
be a holomorphic map such that $p_1=p_2\circ \pi_2\circ\Phi_0$.
We set $\Phi_1:=\pi_2\circ\Phi_0:\nbigz_1\to \nbigz_2$.
We assume the following conditions.
\begin{itemize}
 \item $\Phi_1$ is proper and finite.
 \item There exists a complex analytic closed hypersurface
       $\nbigd\subset\nbigs\times B(R_1)\subset\nbigz_2$
       for some $0<R_1<1$
       such that
       (i) the induced map
       $\nbigz_1\setminus\Phi_1^{-1}(\nbigd)
       \lrarr
       \nbigz_2\setminus\nbigd$
       is a covering map,
       (ii) $\Phi_0$ induces
       an injection $\nbigz_1\setminus\Phi_1^{-1}(\nbigd)\to
       T^{\ast}(\nbigz_2\setminus \nbigd)$,
       (iii) $\nbigd\cap(\{x_0\}\times\cnum)=\{(x_0,0)\}$.
\end{itemize}
We set $r:=|\Phi_1^{-1}(P)|$ for any $P\in\nbigz_2\setminus\nbigd$.
We also set
$\nbigdtilde=\Phi_1^{-1}(\nbigd)$.

\begin{lem}
$\nbige=\Phi_{1\ast}(\nbigo_{\nbigz_1})$
is a locally free
$\nbigo_{\nbigz_2}$-module of rank $r$.
\end{lem}
\pf
By a change of local holomorphic coordinate system
on $\nbigz_1$,
it is enough to consider the case where
$\Phi_1^{\ast}(z_2)$ is expressed as a Weierstrass polynomial.
Then, it is reduced to
\cite[Chapter 2, \S4.2, Theorem]{Grauert-Remmert}.
\hfill\qed

\vspace{.1in}
Note that $\nbige=\pi_{2\ast}(\Phi_{0\ast}\nbigo_{\nbigz_1})$ is naturally
a $\pi_{2\ast}(\nbigo_{T^{\ast}(\nbigz_2/\nbigs)})$-module.
Hence, we obtain the relative Higgs field
$\theta:\nbige\to
\nbige\otimes \Omega^1_{\nbigz_2/\nbigs}$.
The following lemma is clear by the construction.
\begin{lem}
For any $P\in \nbigz_2\setminus\nbigd$,
there exist a neighbourhood $\nbigu$ of $P$
in $\nbigz_2\setminus\nbigd$
and a decomposition 
\begin{equation}
\label{eq;22.12.24.2}
 (\nbige,\theta)_{|\nbigu}
  =\bigoplus_{i=1}^{r}
  (\nbige_{P,i},\theta_{P,i}),
\end{equation}
where $\rank\nbige_{P,i}=1$,
and $\theta_{P,i}-\theta_{P,j}$ $(i\neq j)$
are nowhere vanishing.
\hfill\qed
\end{lem}
For any $x\in\nbigs$,
we set
$\nbigz_{i,x}=\nbigz_i\cap(\{x\}\times\cnum)$,
$\nbigdtilde_x=\nbigdtilde\cap(\{x\}\times\cnum)$
and $\nbigd_x=\nbigd\cap(\{x\}\times\cnum)$.
Note that $\nbigz_{2,x}=B(2)$ for any $x\in\nbigs$.
Let $\iota_x:\nbigz_{2,x}\to \nbigz_2$
denote the inclusion.
We obtain the Higgs bundles
$(\nbige_x,\theta_x):=\iota_x^{\ast}(\nbige,\theta)$
on $\nbigz_{2,x}$
which is regular semisimple outside $\nbigd_{x}$.

\subsubsection{A uniform estimate in the symmetric case}
\label{subsection;23.1.28.20}

Let $h^{\circ}_{x}$ $(x\in\nbigs)$
be decoupled harmonic metrics of
$(\nbige_{x},\theta_{x})
_{|B(1)\setminus\nbigd_{x}}$
such that they induce a $C^{\infty}$-metric of
$\nbige_{|\nbigz_2\setminus \nbigd}$.
Assume the following.
\begin{condition}
\label{condition;23.1.28.30}
For each $(x,P)\in\nbigd$,
Condition {\rm\ref{condition;23.1.28.32}}
is satisfied for $(\nbige_x,\theta_x,h^{\circ}_x)$ at $P$.
\hfill\qed 
\end{condition}

Let $h_{x,t}$ be harmonic metrics of
$(\nbige_x,t\theta_x)_{|B(1)}$
such that
$h_{x,t|\del B(1)}=h^{\circ}_{x|\del B(1)}$.
Let $s_{x,t}$ be the automorphism of $\nbige_{x|B(1)}$
determined by
$h_{x,t}=h^{\circ}_{x}\cdot s_{x,t}$.

\begin{prop}
\label{prop;23.1.28.40}
Let $R_1<R_2<1$.
Let $\nbigs'$ be a relatively compact open subset of $\nbigs$.
Then, there exists $t_0>0$ such that the following holds.
\begin{itemize}
 \item 
 For any $\ell\in\seisuu_{\geq 0}$,
there exist positive constants $C(\ell)$ and $\epsilon(\ell)$
such that
\[
 \bigl|
(s_{x,t}-\id)_{|B(R_2)\setminus B(R_1)}\bigr|_{L_{\ell}^2}
 \leq
 C(\ell)\exp(-\epsilon(\ell) t)
\]
for any $x\in\nbigs'$ and any $t\geq t_0$.
Here, we consider the $L_{\ell}^2$-norms
with respect to $h^{\circ}_x$ and
the standard Euclidean metric $dz_2\,d\zbar_2$.
\end{itemize}
\end{prop}
\pf
For $0<R\leq 2$,
we set
$\nbigz_{1,x}(R):=\Phi_1^{-1}(\{x\}\times B(R))
\subset \nbigz_{1,x}$.
\begin{lem} 
\label{lem;23.1.28.31}
If $R_1<R\leq 2$,
each connected component of $\nbigz_{1,x}(R)$
is diffeomorphic to a $2$-dimensional disc.
\end{lem}
\pf
Let us consider the case $R_1<R<2$.
We set
$\nbigz_{1}(R):=\Phi_1^{-1}(\nbigs\times B(R))
\subset \nbigz$.
It is a compact $C^{\infty}$-manifold
with smooth boundary.
The projection $\nbigz_{1}(R)\to \nbigs$
is submersive.
Each connected component of
$\nbigz_{1,x_0}(R)$ is diffeomorphic to a disc.
Because $\nbigs$ is connected,
we obtain that
each connected component of
$\nbigz_{1,x}(R)$ is diffeomorphic to a disc.
For $R_1<R<2$,
there exists a diffeomorphism
$\rho_R:B(R)\simeq B(2)$
whose restriction to $B(R_1)$ is the identity.
We can construct
a diffeomorphism
$\nbigz_{1,x}(R)\simeq\nbigz_{1,x}(2)$
by lifting $\rho_R$.
\hfill\qed

\begin{lem}
\label{lem;22.12.25.2}
There exist holomorphic non-degenerate symmetric pairings
$C_x$ $(x\in\nbigs)$ of $(\nbige_{x},\theta_{x})$
such that
the restrictions
$C_{x|B(1)\setminus\nbigd_{x}}$
are compatible with $h^{\circ}_x$
and continuous with respect to $x$.
\end{lem}
\pf
Let $h^{\circ}_{0,x}$ denote the flat metric of
$\nbigo_{\nbigz_{1,x}\setminus\nbigdtilde_{x}}$
corresponding to $h^{\circ}_x$,
which are continuous with respect to $x$.
Let $\nabla^{\circ}_{0,x}$ denote the Chern connection.
They are flat connections,
and continuous with respect to $x$.

By Proposition \ref{prop;22.12.23.110}
and Lemma \ref{lem;23.1.28.31},
for each $x\in\nbigs$,
there exists a holomorphic non-degenerate symmetric pairing
$C'_x$ of $(\nbige_{x},\theta_{x})$
such that
the restriction
$(C'_{x})_{|B(1)\setminus\nbigd_{x}}$
is compatible with $h^{\circ}_x$.
Let $C'_{0,x}$ denote the holomorphic non-degenerate
symmetric bilinear form 
of $\nbigo_{\nbigz_{1,x}\setminus\nbigdtilde_{x}}$
corresponding to $C'_x$,
which is compatible with $h^{\circ}_{0,x}$.

Let $z_1\in B(1)\setminus B(R_1)$.
There exists a continuous family of
non-degenerate symmetric pairings
$C^{\circ}_{0,(x,z_1)}$ of
the vector space
$\nbigo_{\nbigz_1|(x,z_1)}$
which are compatible with $(h^{\circ}_{0,x})_{|z_1}$.
We obtain $\alpha_x\in\cnum^{\ast}$
determined by
$C^{\circ}_{0,(x,z_1)}=\alpha_x (C'_{0,x})_{|z_1}$.
We set $C_{0,x}=\alpha_x C'_{0,x}$.
Because 
$C_{0,x}$ are $\nabla^{\circ}_{0,x}$-flat,
they are continuous with respect to $x$.
Let $C_{x}$ denote the non-degenerate
symmetric pairing of $(\nbige_x,\theta_x)$
corresponding to $C_{0,x}$.
(See Proposition \ref{prop;22.12.23.110}.)
Then, they satisfy the desired condition.
\hfill\qed

\vspace{.1in}
Because $h_{t,x|\del B(1)}=h^{\circ}_{x|\del B(1)}$
are compatible with $C_{x|\del B(1)}$,
we obtain that
$h_{t,x}$ are compatible with $C_x$.
Then, the claim of Proposition \ref{prop;23.1.28.40}
follows from Theorem \ref{thm;22.12.3.1}.
\hfill\qed

\vspace{.1in}
We also obtain the following proposition
from Theorem \ref{thm;22.12.3.1},
as in the proof of Proposition \ref{prop;23.1.28.40}.
\begin{prop}
\label{prop;23.1.28.10}
Let $R_1<R_2<2$.
Let $\nbigs'$ be a relatively compact open subset of $\nbigs$.
There exists $t_0>0$ such that the following holds.
\begin{itemize}
 \item
Let $h'_{x,t}$ be any harmonic metrics of $(\nbige_x,t\theta_x)$
$(x\in\nbigs')$
compatible with $C^{\circ}_x$.
Let $s'_{x,t}$ be determined by
$h'_{x,t}=h^{\circ}_x\cdot s'_{x,t}$.
Then, 
for any $\ell\in\seisuu_{\geq 0}$,
there exist positive constants $C(\ell)$ and $\epsilon(\ell)$
such that
\[
 \bigl|
(s'_{x,t}-\id)_{|B(R_2)\setminus B(R_1)}\bigr|_{L_{\ell}^2}
 \leq
 C(\ell)\exp(-\epsilon(\ell) t)
\]
for any $t\geq t_0$.
\hfill\qed
\end{itemize}
\end{prop}

\subsubsection{Examples of non-degenerate symmetric pairings
and decoupled harmonic metrics}
\label{subsection;22.12.24.120}

We obtain a holomorphic function
$G=\del_{z_1}(\Phi_1^{\ast}(z_2))$.
We have
$G^{-1}(0)\subset\nbigdtilde$.
We define the symmetric product
$C_0:\nbigo_{\nbigz_1}\otimes\nbigo_{\nbigz_1}
\lrarr G^{-1}\nbigo_{\nbigz_1}$
by
\[
 C_0(a\otimes b)=G^{-1}ab.
\]
We obtain the following lemma
by using Lemma \ref{lem;22.12.24.1}.
\begin{lem}
\label{lem;22.12.24.10}
$C_0$ induces a non-degenerate symmetric pairing $C_1$
of $\nbige$,
which induces a non-degenerate symmetric pairing of
$(\nbige_x,\theta_x)$ for any $x\in\nbigs$.
\hfill\qed 
\end{lem}

Let $h_0$ be the flat metric of $\nbigo_{\nbigz_1\setminus\nbigdtilde}$
defined as follows:
\[
 h_0(a,b)=|G|^{-1}a\overline{b}.
\]
\begin{lem}
\label{lem;22.12.24.11}
$h_0$ induces
a flat metric $h_1$ of $\nbige_{|\nbigz_2\setminus \nbigd}$.
For each $x\in\nbigs$,
the induced metric $h_{1,x}$
of $(\nbige_x,\theta_x)_{|\nbigz_{2,x}\setminus\nbigd_{x}}$
is a decoupled harmonic metric
such that $\det(h_{1,x})$ induces a flat metric of
$\det(\nbige_{x})$ for each $x\in\nbigs$.
\hfill\qed
\end{lem}
\begin{rem}
We shall use $h_0$ in {\rm\S\ref{subsection;22.12.25.100}}. 
\hfill\qed
\end{rem}

\section{Large-scale solutions on compact Riemann surfaces}

\subsection{Convergence in the locally irreducible case}

\subsubsection{Statement}

Let $X$ be a compact Riemann surface.
Let $\pi:T^{\ast}X\to X$ denote the projection.
For any $A\subset T^{\ast}X$,
the induced map $A\to X$ is also denoted by $\pi$.
Let $D\subset X$ be a finite subset.

Let $(\nbigp_{\ast}\nbigv,\theta)$
be a good filtered Higgs bundle of degree $0$
on $(X,D)$.
We obtain the Higgs bundle
$(V,\theta)=(\nbigv,\theta)_{|X\setminus D}$.
We assume the following.
\begin{condition}
$(V,\theta)$ is regular semisimple Higgs bundle on
$X\setminus D$.
\hfill\qed
\end{condition}

\begin{rem}
\label{rem;22.12.24.20}
If $(V,\theta)$ is generically regular semisimple,
there exists a finite subset $D'\subset X$
such that $(V',\theta')_{|X\setminus D'}$
is regular semisimple
and that $D\subset D'$.
We set $\nbigv'=\nbigv(\ast D')$.
For each $P\in D'\setminus D$,
we consider the filtered bundle
$\nbigp_{\ast}(\nbigv'_P)$ over $\nbigv'_P$
defined by
$\nbigp_a\nbigv'_P
=\nbigv_P([a]P)$,
where $[a]=\max\{n\in\seisuu\,|\,n\leq a\}$.
For harmonic metrics of $(V,t\theta)$
adapted to $\nbigp_{\ast}\nbigv$,
it is enough to study
harmonic metrics of $(V',t\theta')$
adapted to $\nbigp_{\ast}\nbigv'$.
\hfill\qed
\end{rem}

For any $P\in D$,
there exist a neighbourhood $X_P$ of $P$ in $X$
and a decomposition of the meromorphic Higgs bundle
\begin{equation}
\label{eq;23.1.28.11}
 (\nbigv,\theta)_{|X_P}
=\bigoplus_{i\in S(P)}
(\nbigv_{P,i},\theta_{P,i}),
\end{equation}
such that the spectral curves of
$(\nbigv_{P,i},\theta_{P,i})_{|X_P\setminus\{P\}}$
are connected.

\begin{condition}
We assume the following conditions.
\begin{itemize}
 \item $\Sigma_{V,\theta}$ is connected.
 \item For any $P\in D$,
       the decomposition {\rm(\ref{eq;23.1.28.11})} is compatible with
       the filtered bundle $\nbigp_{\ast}(\nbigv_P)$ over $\nbigv_P$,
       i.e.,
       $\nbigp_{\ast}(\nbigv_P)
       =\bigoplus_{i\in S(P)}\nbigp_{\ast}\bigl(
       (\nbigv_{P,i})_P\bigr)$.
\hfill\qed
\end{itemize}
\end{condition}

For each $P\in D$,
we obtain the filtered bundle
$\nbigp^{\star}_{\ast}(\nbigv_{P})
=\bigoplus_{i\in S(P)}
\nbigp^{\star}_{\ast}\bigl(
(\nbigv_{P,i})_P\bigr)$
over $\nbigv_{P}$
determined by
the filtered bundles $\det(\nbigp_{\ast}\nbigv_{P,i})$
as in Proposition \ref{prop;22.12.20.1}.
By patching
$\nbigp^{\star}_{\ast}(\nbigv_{P})$
$(P\in D)$ with $\nbigv$,
we obtain a decomposable filtered Higgs bundle
$(\nbigp^{\star}_{\ast}(\nbigv),\theta)$.

\begin{lem}
\label{lem;22.11.30.110}
$(\nbigp^{\star}_{\ast}(\nbigv),\theta)$
is stable of degree $0$.
As a result,
there exists a decoupled harmonic metric $h_{\infty}$
of $(V,\theta)$
adapted to $\nbigp^{\star}_{\ast}(\nbigv)$.
\hfill\qed 
\end{lem}
\pf
Because $\Sigma_{V,\theta}$ is connected,
there does not exist a non-trivial Higgs subbundle of
$(V,\theta)$.
Hence, $(\nbigp_{\ast}\nbigv,\theta)$ is stable.
Because
$\det(\nbigp^{\star}_{\ast}\nbigv)
=\det(\nbigp_{\ast}\nbigv)$,
we obtain
$\deg(\nbigp^{\star}_{\ast}\nbigv)=0$.
The second claim follows from
Proposition \ref{prop;22.12.21.20}.

\hfill\qed

\vspace{.1in}
Note that
$\det(h_{\infty})$ is a flat metric of $\det(V)$
adapted to
$\det(\nbigp_{\ast}\nbigv)=\det(\nbigp^{\star}_{\ast}\nbigv)$.
Because $\Sigma_{V,\theta}$  is connected,
$(\nbigp_{\ast}\nbigv,\theta)$ is stable of degree $0$
as in Lemma \ref{lem;22.11.30.110}.
Hence, for any $t>0$,
there exists a harmonic metric $h_t$
of $(V,t\theta)$ which is adapted to $\nbigp_{\ast}\nbigv$
such that $\det(h_t)=\det(h_{\infty})$.

\begin{thm}
\label{thm;22.11.30.130}
On any relatively compact open subset $K\subset X\setminus D$,
the sequence $h_{t}$ is convergent to $h_{\infty}$
in the $C^{\infty}$-sense. 
\end{thm}

\subsubsection{The case of locally and globally irreducible Higgs bundles}
\label{subsection;22.11.30.200}

We state Theorem \ref{thm;22.11.30.130}
in a special case for explanation
(see also Remark \ref{rem;22.12.24.20}).
Let $(E,\delbar_E,\theta)$ be
a generically regular semisimple Higgs bundle
of degree $0$ on $X$.
Let $\Sigma_{E,\theta}$ denote the spectral curve.
There exists the finite subset $D(E,\theta)\subset X$
such that the following holds.
\begin{itemize}
 \item $P\in D(E,\theta)$
       if and only if
       $|T_P^{\ast}X\cap\Sigma_{E,\theta}|<r$.
\end{itemize}

We impose the following condition.
\begin{condition}
\mbox{{}}\label{condition;22.11.30.40}
\begin{itemize}
 \item $\Sigma_{E,\theta}$ is irreducible,
       i.e.,
      $\Sigma_{E,\theta}\setminus \pi^{-1}(D(E,\theta))$ is connected.
 \item For any $P\in D(E,\theta)$,
       there exist a neighbourhood
       $X_P$ of $P$ in $X$
       and a decomposition
\begin{equation}
\label{eq;22.11.30.20}
       (E,\theta)_{|X_P}
       =\bigoplus_{i\in S(P)}
       (E_{P,i},\theta_{P,i})
\end{equation}
       such that the spectral curves
       $\Sigma_{E_{P,i},\theta_{P,i}}$
       are irreducible.
\hfill\qed
\end{itemize}
\end{condition}
We set
$D=D(E,\theta)$.
Let $\nbigp^{(0)}_{\ast}(E(\ast D)_P)$
be the filtered bundle over $E(\ast D)_P$
defined by
$\nbigp^{(0)}_a(E(\ast D)_P)=E_P([a]P)$,
where $[a]=\max\{n\in\seisuu\,|\,n\leq a\}$.
Because there exists the decomposition
\[
\nbigp^{(0)}_{\ast}(E(\ast D)_P)
=\bigoplus_{i\in S(P)}
\nbigp^{(0)}_{\ast}(E_{P,i}(\ast D)_P)
\]
induced by (\ref{eq;22.11.30.20}),
we obtain
the filtered bundle
$\nbigp^{\star}_{\ast}(E(\ast D)_P)$
determined by
$\det(\nbigp^{(0)}_{\ast}E_{P,i}(\ast D)_P)$
as in Proposition \ref{prop;22.12.20.1}.
By patching them with $(E(\ast D),\theta)$,
we obtain a filtered bundle
$\nbigp^{\star}_{\ast}\nbigv$ over
$\nbigv=E(\ast D)$.
The filtered Higgs bundle
$(\nbigp^{\star}_{\ast}(E(\ast D)),\theta)$
is decomposable.

As in Lemma \ref{lem;22.11.30.110},
there exists a decoupled harmonic metric
$h_{\infty}$ of $(E,\theta)_{|X\setminus D}$
such that
$h_{\infty}$ is adapted to $\nbigp^{\star}_{\ast}\nbigv$.
For any $t>0$,
there exists a unique harmonic metric $h_t$
of $(E,t\theta)$
such that $\det(h_t)=\det(h_{\infty})$.
As a special case of 
Theorem \ref{thm;22.11.30.130},
we obtain the following.
\begin{cor}
\label{cor;22.12.24.30}
On any relatively compact open subset $K\subset X\setminus D$,
the sequence $h_{t}$ is convergent to $h_{\infty}$
in the $C^{\infty}$-sense. 
\hfill\qed
\end{cor}

\begin{rem}
The second condition
in Condition {\rm\ref{condition;22.11.30.40}}
is satisfied
if $\Sigma_{E,\theta}$ is locally irreducible.
\hfill\qed 
\end{rem}

\subsubsection{Proof of Theorem \ref{thm;22.11.30.130}}
\label{subsection;22.11.30.201}

Let $P\in D$.
We set $X_P^{\ast}=X_P\setminus \{P\}$.
We set $V_{P,i}=\nbigv_{P,i|X_P^{\ast}}$,
and $r(P,i)=\rank V_{P,i}$.
Let $z_P$ be a holomorphic coordinate of $X_P$
by which $X_P\simeq \bigl\{z\in\cnum\,\big|\,|z|<2\bigr\}$.
We set $(h_{\infty})_{P,i}:=h_{\infty|V_{P,i}}$.
Let $h_{t,P,i}$ be a harmonic metric of
$(V_{P,i},t\theta_{P,i})$ such that
(i) the boundary value at $|z_P|=1$
is equal to that of $(h_{\infty})_{P,i}$,
(ii) $h_{t,P,i}$ is adapted to
$\nbigp_{\ast}\nbigv_{P,i}$.
We have $\det(h_{t,P,i})=\det((h_{\infty})_{P,i})$.
We obtain the following lemma by Proposition \ref{prop;22.11.30.131}.
\begin{lem}
The sequence $h_{t,P,i}$ is convergent to $(h_{\infty})_{P,i}$
as $t\to\infty$
in the $C^{\infty}$-sense
on any relatively compact open subset of $X_P^{\ast}$. 
\hfill\qed
\end{lem}

We regard $X_P$ as an open subset of $\cnum$
by $z_P$.
Let $\varphi_{P,r(P,i)}:\cnum\to\cnum$
be defined by
$\varphi_{P,r(P,i)}(\zeta_{P,i})=\zeta_{P,i}^{r(P,i)}$.
We set
$X_P^{(r(P,i))}=\varphi_{P,r(P,i)}^{-1}(X_P)$
and
$X_P^{(r(P,i))\ast}=\varphi_{P,r(P,i)}^{-1}(X_P^{\ast})$.
The induced maps
$X_P^{(r(P,i))}\to X_P$
and
$X_P^{(r(P,i))\ast}\to X_P^{\ast}$
are also denoted by $\varphi_{P,r(P,i)}$.

We define a Hermitian product $h_{t,P,i}^{(r(P,i))}$
of $\varphi_{P,r(P,i)}^{\ast}(V_{P,i})_{|X_P^{(r(P,i))\ast}}$ as follows.
We have the decomposition
\[
\varphi_{P,r(P,i)}^{\ast}
(V_{P,i},\theta_{P,i})_{|X_P^{(r(P,i))\ast}}
=\bigoplus_{p=1}^{r(P,i)}
 (V_{P,i,\beta(p)},\beta(p)\,d\zeta_{P,i}),
\]
where $\beta(p)$ are meromorphic functions on $X_P^{(r(P,i))}$.
Let $v_{\beta(1)}$ be a holomorphic frame of $V_{P,i,\beta(1)}$.
We obtain a frame
$v_{\sigma^{\ast}\beta(1)}=\sigma^{\ast}(v_{\beta(1)})$
of 
$V_{P,i,\sigma^{\ast}\beta(1)}$.
Let $\chi(\zeta_{P,i})$
be an $\real_{\geq 0}$-valued function
such that
(i) $\chi(\zeta_{P,i})$ depends only on $|\zeta_{P,i}|$,
(ii) $\chi(\zeta_{P,i})=1$ $(|\zeta_{P,i}|\leq 1/2)$,
$\chi(\zeta_{P,i})=0$ $(|\zeta_{P,i}|\geq 2/3)$.
For $p\neq q$,
we put
\[
 h_{t,P,i}^{(r(P,i))}(v_{\beta(p)},v_{\beta(q)})
 =\chi(\zeta_{P,i})
 \varphi_{P,r(P,i)}^{\ast}(h_{t,P,i})(v_{\beta(p)},v_{\beta(q)}).
\]
We define
$h_{t,P,i}^{(r(P,i))}(v_{\beta(p)},v_{\beta(p)})$ by
\begin{multline}
\log h_{t,P,i}^{(r(P,i))}(v_{\beta(p)},v_{\beta(p)})
 =\\
 \chi(\zeta_{P,i})
\log
\varphi_{P,r(P)}^{\ast}(h_{t,P,i})(v_{\beta(p)},v_{\beta(p)})
+(1-\chi(\zeta(P,i)))
\log
\varphi_{P,r(P)}^{\ast}\bigl(
(h_{\infty})_{P,i}^{(r(P,i))}
\bigr)(v_{\beta(p)},v_{\beta(p)}).
\end{multline}
Then, $h^{(r(P,i))}_{t,P,i}$ is $\Gal(r(P,i))$-invariant,
and we have
$h^{(r(P,i))}_{t,P,i}=\varphi_{P,r(P,i)}^{-1}(h_{t,P,i})$
on $\{0<|\zeta_{P,i}|<1/4\}$
and
$h^{(r(P,i))}_{t,P,i}
=\varphi_{P,r(P,i)}^{-1}\bigl((h_{\infty})_{P,i}\bigr)$
on $\{4/5<|\zeta_{P,i}|\}$.
There exists a Hermitian metric $\htilde_{t,P,i}$ of $V_{P,i}$
such that
$\varphi_{P,r(P,i)}^{-1}(\htilde_{t,P,i})
=h^{(r(P,i))}_{t,P,i}$
on $X_P^{(r(P,i))\ast}$.
We obtain a Hermitian metric
\[
\htilde_{t,P}=
\bigoplus_{i\in S(P)}
\htilde_{t,P,i}
\]
of $V_{|X_P^{\ast}}$.
By patching $\htilde_{t,P}$ and $h_{\infty}$,
we obtain Hermitian metrics $\htilde'_{t}$ of $V$.
We obtain the $C^{\infty}$-function $\alpha_t$
on $X\setminus D$
determined by
$\det(\htilde'_t)=e^{\alpha_t}\det(h_{\infty})$.
We set
$\htilde_t=e^{-\alpha_t/r}\htilde'_t$.
By the construction, the following lemma is clear.
\begin{lem}
There exists $t_0$ such that
$\htilde_{t}$ is positive definite
for any $t\geq t_0$.
Moreover, the following holds.
 \begin{itemize}
 \item The sequence $\htilde_t$ is convergent to $h_{\infty}$
      in the $C^{\infty}$-sense
      on any relatively compact open subset of $X\setminus D$.
      The support of
$R(\htilde_{t})+
[t\theta,(t\theta)^{\dagger}_{\htilde_{t}}]$
is contained in
$\{\left(\frac{1}{4}\right)^{\rank(E)}
\leq
|z_P|\leq
 \frac{4}{5}\}$ for $P\in D$.
In particular, 
\begin{equation}
\label{eq;22.11.30.30}
\int_X
 \Bigl|
 R(\htilde_{t})+[t\theta,(t\theta)^{\dagger}_{\htilde_{t}}]
 \Bigr|_{\htilde_{t},g_X}
\to 0 
\end{equation}
as $t\to\infty$.
\hfill\qed
 \end{itemize}
\end{lem}

Let $g_X$ be a K\"ahler metric of $X$.
Let $s_{t}$ denote the automorphism of $V$
determined by
$h_t=\htilde_{t}s_{t}$.
We have $\det(s_{t})=1$.
According to \cite[Lemma 3.1]{s1},
we obtain the following on $X\setminus D$:
\begin{equation}
\label{eq;22.11.30.31}
 \Delta_{g_X}\Tr(s_{t})
=
\Tr\Bigl(
\bigl(
R(\htilde_t)
+[t\theta,(t\theta)^{\dagger}_{\htilde_t}]
\bigr)
s_{t}\Bigr)
 -\bigl|
 \delbar(s_{t})s_{t}^{-1/2}
 \bigr|^2_{\htilde_{t},g_X}
 -\bigl|
 [t\theta,s_{t}]s^{-1/2}
 \bigr|^2_{\htilde_{t},g_X}.
\end{equation}

Note that
$\bigoplus_{i\in S(P)} h_{t,P,i}$
and $h_{t|X_P^{\ast}}$ are mutually bounded
for any $P\in D$.
Hence, $\Tr(s_t)$ is bounded.
We also note the following vanishing
(see Lemma \cite[Lemma 4.7]{Mochizuki-KH-Higgs}):
\begin{equation}
\label{eq;22.11.30.140}
 \int_X\Delta_X\Tr(s_t)\vol_{g_X}=0.
\end{equation}

We set $b_{t}=\sup_{X\setminus D}\Tr(s_{t})$.
Note that $b_{t}\geq \rank(E)$,
and $b_t=\rank(E)$ if and only if $s_t=\id_E$.
We set $u_{t}=b_{t}^{-1}\cdot s_{t}$.
There exists $C>0$,
which is independent of $t$
such that
$|u_{t}|_{\htilde_{t}}\leq C$.
By (\ref{eq;22.11.30.30}), (\ref{eq;22.11.30.31})
and (\ref{eq;22.11.30.140})
we obtain
\[
 \int_X\Bigl(
 |\delbar u_{t}|^2_{\htilde_{t}}
 +|[t\theta,u_{t}]|^2_{\htilde_{t}}
 \Bigr)
 \to 0
\]
as $t\to\infty$.

Let $t(i)>0$ be a sequence such that $t(i)\to\infty$ as $i\to\infty$.
By going to a subsequence,
$u_{t(i)}$ is weakly convergent in $L_1^2$
locally on $X\setminus D$.
In particular, it is convergent in $L^q$ for any $q\geq 1$
locally on $X\setminus D$.
Let $u_{\infty}$ denote the limit
which satisfies $\delbar u_{\infty}=[\theta,u_{\infty}]=0$.

\begin{lem}
$u_{\infty}\neq 0$.
\end{lem}
\pf
Note that
$\sup_X\Tr(u_{t(i)})=1$ for any $i$.
Let $0<\epsilon<1$.
Let $P_{i}\in X$ be points 
such that
$\Tr(u_{t})(P_{i})\geq \epsilon$.
By going to a subsequence,
we may assume that the sequence is convergent to
a point $P_{\infty}$.
Let us consider the case where
\[
P_{\infty}\not\in \bigcup_{P\in D}
\{|z_{P}|\leq 4/5\}=:W.
\]
Let $(X_{P_{\infty}},z)$ be a holomorphic coordinate neighbourhood
around $P_{\infty}$,
which does not intersect with $W$.
Because $F(\htilde_{t})=0$ on $X_{P_{\infty}}$,
we obtain
$\Delta_{g_X}\Tr(u_{t})\leq 0$.
By the mean value property
of the subharmonic functions,
there exists $C>0$ such that
\[
 C\epsilon
 \leq
 \int_{X_{P_{\infty}}}
 \Tr(u_{t(i)}).
\]
Because $u_{t(i)}$ is convergent to
$u_{\infty}$ in $L^p$ for any $p\geq 1$
on $X_{P_{\infty}}$,
we obtain that $u_{\infty}\neq 0$.

Let us consider the case
where $P_{\infty}\in \{|z_P|<4/5\}$ for some $P\in D$.
Let $(X_P,z_P)$ be a holomorphic coordinate neighbourhood
around $P$ as above.
By \cite[Lemma 3.1]{s1},
we have
\[
 \Delta_{g_X}
 \log\Tr(u_{t(i)})
 \leq
 \Bigl|R(\htilde_{t(i)})+[t\theta,(t\theta)^{\dagger}_{\htilde_{t(i)}}]
 \Bigr|_{\htilde_{t(i)},g_X}.
\]
There exist $C^{\infty}$-functions $\alpha_{i}$
on $X_{P}$
such that
(i)
$\Delta_{g_X} \alpha_{i}
 =\Bigl|R(\htilde_{t(i)})+[t\theta,(t\theta)^{\dagger}_{\htilde_{t(i)}}]
 \Bigr|_{\htilde_{t(i)},g_X}$,
(ii) $\alpha_{i|\del X_{P}}=0$,
(iii) there exists $C>0$
such that $|\alpha_{i}|\leq C$ for any $i$.
Because $\log\Tr(u_{t(i)})-\alpha_i$
is a subharmonic function on $X_P$,
the maximum principle allows us to obtain
\[
 \log \epsilon-C
 \leq
 \max_{P\in \del X_{P}}\bigl\{
 \log\Tr(u_{t(i)})
 -\alpha_{i}
 \bigr\}
= \max_{P\in \del X_{P}}\bigl\{
 \log\Tr(u_{t(i)})
 \bigr\}.
\]
Hence, there exists a sequence
$P'_i\in\del X_{P}$
such that
$\Tr(u_{t(i)})(P'_i)\geq \epsilon e^{-C}$.
By going to a subsequence,
we may assume that
the sequence $P'_i$ is convergent
to
$P'_{\infty}\in X\setminus W$.
Then, we can apply the result in the first part
of this proof.
\hfill\qed

\vspace{.1in}
Recall that $u_{\infty}\neq 0$ is an endomorphism of
$(V,\theta)$ such that
$\delbar u_{\infty}=[\theta,u_{\infty}]=0$.
At each point of $X\setminus D$,
an eigenspace of $\theta$ is
also an eigenspace of $u_{\infty}$.
Because each $u_{t(i)}$ is self-adjoint with respect to $\htilde_{t}$,
$u_{\infty}$ is self-adjoint with respect to $h_{\infty}$.
We obtain
$\del_{h_{\infty}}u_{\infty}=0$.
Hence, the eigenvalues of $u_{\infty}$ are constant.
Because
$\htilde_t(u_{t(i)}v,v)\geq 0$
for any local section $v$ of $V$,
we obtain
$h_{\infty}(u_{\infty}v,v)\geq 0$,
which implies that
the eigenvalues of $u_{\infty}$ are non-negative.
We also note that $\Sigma_{V,\theta}$ is connected.
Hence,
$u_{\infty}$ is a positive constant multiplication.
It implies that the sequence $b_{t}$ is bounded,
and that the subsequence $s_{t(i)}$ is convergent
to a positive constant multiplication.
Because $\det(s_{t})=1$,
the limit is the identity.
Because this is independent of
the choice of a subsequence,
we obtain the desired convergence.

\hfill\qed

\subsection{Order of convergence in a smooth case}

\subsubsection{Rough statement}

Let us study the order of the convergence
in the situation of \S\ref{subsection;22.11.30.200}
assuming the following stronger condition.
\begin{condition}
\label{condition;22.11.30.41}
Let $\rho:\Sigmatilde_{E,\theta}\to\Sigma_{E,\theta}$
be the normalization. 
There exists a line bundle $L$ on $\Sigmatilde_{E,\theta}$ 
with an isomorphism
$(\pi\circ\rho)_{\ast}L\simeq E$
such that  
$\theta$ is induced by
the $\nbigo_{T^{\ast}X}$-action on $\rho_{\ast}L$.
\end{condition}
Let $g(\Sigmatilde_{E,\theta})$ and $g(X)$
denote the genus of $\Sigmatilde_{E,\theta}$ and $X$,
respectively.
Then, we have
$\deg(L)=g(\Sigmatilde_{E,\theta})-rg(X)+r-1$.

\begin{rem}
If Condition {\rm\ref{condition;22.11.30.41}} is satisfied,
Condition {\rm\ref{condition;22.11.30.40}} is also satisfied.
Condition {\rm\ref{condition;22.11.30.41}} is satisfied,
if $\Sigma_{E,\theta}$ is smooth and connected.
\hfill\qed
\end{rem}

We set $(V,\theta)=(E,\theta)_{|X\setminus D}$.
Let $s(h_{\infty},h_t)$
be the automorphism of $V$
determined by
$h_t=h_{\infty}\cdot s(h_{\infty},h_t)$.
Let $g_X$ be a K\"ahler metric of $X$.
\begin{thm}
\label{thm;22.12.24.101}
For any relatively compact open subset $K$ of $X\setminus D$
and a non-negative integer $\ell$,
there exist positive constants
$C(K,\ell)$ and $\epsilon(k,\ell)$
such that the $L_{\ell}^2$-norm of
$s(h_{\infty},h_t)-\id_E$ on $K$
with respect to $h_{\infty}$ and $g_X$
are dominated by 
$C(K,\ell)e^{-\epsilon(k,\ell)t}$.
\end{thm}

\subsubsection{Refined statement}
We shall prove a refined statement.
For that purpose,
we refine the construction of $\htilde_t$
in the proof of Theorem \ref{thm;22.11.30.130}.
Let $P\in D$ and $i\in S(P)$.
\begin{lem}
\label{lem;23.1.28.50}
$\det((h_{\infty})_{P,i})$
induces a flat metric of
$\det(E_{P,i})$.
\end{lem}
\pf
It follows from the condition that
$\det((h_{\infty})_{P,i})$
is adapted to
$\det\nbigp^{\star}_{\ast}(E_{P,i}(\ast D)_P)
=\det\nbigp^{(0)}_{\ast}(E_{P,i}(\ast D)_P)$.
\hfill\qed

\vspace{.1in}

According to Proposition \ref{prop;22.12.23.110},
there exists a non-degenerate symmetric pairing $C_{P,i}$ of
$(E_{P,i},\theta_{P,i})$
such that $C_{P,i|X_P^{\ast}}$ is compatible with
$(h_{\infty})_{P,i}$. 
For $t>0$,
there exists a harmonic metric
$h_{t,P,i}$ of $(E_{P,i},\theta_{P,i})$
which is compatible with $C_{P,i}$
such that
its boundary value at $\del X_P$
is equal to that of $h_{\infty|E_{P,i}}$.
We construct the metric $\htilde_{t}$
by using $h_{t,P,i}$ as in the proof of Theorem \ref{thm;22.11.30.130}
(see \S\ref{subsection;22.11.30.201}).
By Proposition \ref{prop;22.11.30.202},
the following holds.

\begin{lem}
\label{lem;22.12.24.100}
Let $s(h_{\infty},\htilde_t)$ be
the automorphism of $E_{|X\setminus D}$
determined by
$\htilde_t=h_{\infty}\cdot s(h_{\infty},\htilde_t)$.
For any relatively compact open subset $K$ of
$X_P^{\ast}$ and for any $\ell\in\seisuu_{\geq 0}$,
there exist $C(K,\ell),\epsilon(K,\ell),t(K)>0$
such that
the $L_{\ell}^2$-norms of
$s(h_{\infty},\htilde_t)-\id$ on $K$
with respect to $h_{\infty}$ and $g_X$
are dominated by $C(K,\ell)\exp(-\epsilon(K,\ell)t)$
for any $t>t(K)$.
\hfill\qed 
\end{lem}

By Lemma \ref{lem;22.12.24.100},
we obtain
\begin{equation}
 \label{eq;22.12.24.31}
 \bigl|
 R(\htilde_{t})
+[t\theta,(t\theta)^{\dagger}_{\htilde_t}]
 \bigr|_{\htilde_t,g_X}\leq Ce^{-\epsilon t}.
\end{equation}
for some $\epsilon,C>0$.
Moreover, by the construction,
the support of (\ref{eq;22.12.24.31}) is contained in
$\bigcup_{P\in D}\{(1/4)^{\rank(E)}
\leq
|z_P|\leq
4/5\}$.

Let $s_t$ be the automorphism of $E$
determined by $h_t=\htilde_t\cdot s_t$.
We obtain Theorem \ref{thm;22.12.24.101}
from Lemma \ref{lem;22.12.24.100}
and the following theorem.

\begin{thm}
\label{thm;22.11.28.20}
For any $\ell>0$,
there exist $C(\ell),\epsilon(\ell)>0$
such that the $L_{\ell}^2$-norm of
$s_t-\id$ on $X$
with respect to $g_X$ and $\htilde_t$
are dominated by 
$C(\ell)\exp(-\epsilon(\ell) t)$.
\end{thm}
\pf
By \cite[Lemma 3.1]{s1} and (\ref{eq;22.12.24.31}),
there exist $C_1,\epsilon_1>0$ such that
\[
 \int_X
 \Bigl(
 \bigl|s_{t}^{-1/2}\del_{E,\htilde_{t}}(s_{t})
 \bigr|_{\htilde_{t}}^2
 +\bigl|
 [\theta,s_{t}]s_{t}^{-1/2}
 \bigr|_{\htilde_{t}}^2
 \Bigr)
\leq 
 C_1\exp(-\epsilon_1t).
\]
By Corollary \ref{cor;22.12.24.30},
$|s_{t}|_{\htilde_{t}}$ and
$|s_{t}^{-1}|_{\htilde_{t}}$ are uniformly bounded.
There exist $C_2,\epsilon_2>0$ such that
\begin{equation}
\label{eq;22.11.21.1}
 \int_X
 \Bigl(
 \bigl|\del_{E,\htilde_{t}}(s_{t})
 \bigr|_{\htilde_{t}}^2
 +\bigl|
 [\theta,s_{t}]
 \bigr|_{\htilde_{t}}^2
 \Bigr)
\leq 
 C_2\exp(-\epsilon_2t).
\end{equation}

Let $K$ be a relatively compact open subset
of $X\setminus D$.
By the variant of Simpson's main estimate
(\cite[Theorem 2.9]{Decouple})
and Lemma \ref{lem;22.12.24.100},
there exist $C_3(K),\epsilon_3(K)>0$
such that
the following holds on $K$:
\[
 \bigl|
 \delbar_E\bigl(
 s_{t}^{-1}
 \del_{E,\htilde_{t}}(s_{t})
 \bigr)
 \bigr|^2_{\htilde_{t}}
 \leq
 C_3(K)
 \exp(-\epsilon_3(K) t).
\]
Together with (\ref{eq;22.11.21.1}),
we obtain that
there exist $C_4(K),\epsilon_4(K)>0$
such that the following holds on $K$:
\begin{equation}
\label{eq;22.12.24.50}
 \bigl|
 \del_{E,\htilde_{t}}(s_{t})
 \bigr|_{\htilde_{t}}
\leq C_4(K) \exp(-\epsilon_4(K) t).
\end{equation}
Because $s_{t}$ is self-adjoint
with respect to $\htilde_{t}$,
we obtain the following on $K$:
\begin{equation}
\label{eq;22.12.24.51}
 \bigl|
 \delbar(s_{t})
 \bigr|_{\htilde_{t}}
\leq
C_4(K) \exp(-\epsilon_4(K) t).
\end{equation}

\begin{lem}
\label{lem;22.11.28.10}
There exist $C(K),\epsilon(K)>0$
such that 
the following holds on $K$:
\[
\bigl|
s_{t}-\id\bigr|_{\htilde_{t}}
\leq C(K)
 \exp(-\epsilon(K) t).
\]
\end{lem}
\pf Let $P$ be any point of $X\setminus D$.
Let $X_P$ be a simply connected neighbourhood of $P$
in $X\setminus D$.
There exists a decomposition into Higgs bundles of rank $1$:
\[
 (E,\theta)_{|X_P}
 =\bigoplus_{i=1}^{\rank(E)}
 (E_{P,i},\theta_{P,i}).
\]
We obtain the decomposition
$s_{t}=\sum (s_{t})_{j,i}$,
where
$(s_{t})_{j,i}:E_{P,i}\to E_{P,j}$.
By \cite[Proposition 2.3]{Decouple},
there exist $C_5(P),\epsilon_5(P)>0$
such that the following for $i\neq j$ on $X_P$:
\begin{equation}
\label{eq;22.11.30.50}
\bigl|(s_{t})_{j,i}\bigr|_{\htilde_{t}}
\leq
C_5(P)
 \exp(-\epsilon_5(P) t).
\end{equation}
By (\ref{eq;22.12.24.50}) and (\ref{eq;22.12.24.51}),
there exist $C_6(P),\epsilon_6(P)>0$
such that
\[
\bigl|
d(s_{t})_{i,i}
\bigr|
\leq
C_6(P)\exp(-\epsilon_6(P)t).
\]
Hence, 
there exist $C_7(P),\epsilon_7(P)>0$
such that
the following holds for any $P_1,P_2\in X_P$:
\[
\bigl|
 (s_{t})_{i,i}(P_1)
 -(s_{t})_{i,i}(P_2)
\bigr|
\leq C_7(P)
\exp(-\epsilon_7(P) t).
\]
Let $i\neq j$.
There exists a loop $\gamma$ in $X\setminus D$
such that
the monodromy of $\Sigma_{E,\theta}$
along $\gamma$ exchanges $E_i$ and $E_j$.
By taking a finite covering of $\gamma$
by relatively compact open subsets
and by applying the above consideration,
we obtain that
there exist $C_8(P),\epsilon_8(P)>0$
such that the following holds for any $P_1\in X_P$:
\begin{equation}
\label{eq;22.11.30.51}
\bigl|
 (s_{t})_{i,i}(P_1)
 -(s_{t})_{j,j}(P_1)
 \bigr|
\leq 
C_8(P)\exp(-\epsilon_8(P) t).
\end{equation}
By (\ref{eq;22.11.30.50}),
there exist $C_9(P),\epsilon_9(P)>0$
such that the following holds on $X_P$:
\begin{equation}
\label{eq;22.11.30.52}
\left|
 \prod_{i=1}^{\rank(E)} (s_{t})_{i,i}
 -1
\right|
\leq 
C_9(P)\exp(-\epsilon_9(P) t).
\end{equation}
By (\ref{eq;22.11.30.51}) and (\ref{eq;22.11.30.52}),
there exist $C_{10}(P),\epsilon_{10}(P)>0$ such that
\[
\bigl|
 (s_{t})_{i,i}-1
\bigr|
\leq
C_{10}(P)\exp(-\epsilon_{10}(P) t).
\]
Then, we obtain the claim of
Lemma \ref{lem;22.11.28.10}.
\hfill\qed

\vspace{.1in}
We obtain the estimate of
$|s_t-\id|_{\htilde_t}$ around $D$
by using Theorem \ref{thm;22.12.2.11}.
We can also obtain the estimate
for the higher derivatives
by using Theorem \ref{thm;22.12.2.11}.
\hfill\qed

\subsection{A family case}
\label{subsection;22.12.25.100}

\subsubsection{Setting}

Let $\nbigs$ be a connected complex manifold.
Let $\nbigy$ be a complex manifold
with a proper smooth morphism $p_1:\nbigy\to \nbigs$.
Let $p_2:\nbigs\times X\to \nbigs$
and $\pi_2:\nbigs\times T^{\ast}X\to\nbigs\times X$ denote the projections.
Let $\Phi_0:\nbigy\to \nbigs\times T^{\ast}X$
be a holomorphic map
such that $p_1=p_2\circ\pi_2\circ \Phi_0$.
We set
$\Phi_1=\pi_2\circ\Phi_0$.
We assume the following conditions.
 \begin{itemize}
 \item Each fiber of $p_1$ is connected and $1$-dimensional.
 \item $\Phi_1$ is proper and finite.
 \item There exists a closed complex analytic hypersurface
       $\nbigd\subset \nbigs\times X$
       such that
       (i) $\nbigd$ is finite over $\nbigs$,
       (ii) the induced map
       $\nbigy\setminus\Phi_1^{-1}(\nbigd)
       \to (\nbigs\times X)\setminus\nbigd$
       is a covering map,
       (iii) $\Phi_0$ induces an injection
       $\nbigy\setminus\Phi_1^{-1}(\nbigd)
       \lrarr \nbigs\times T^{\ast}X$.
 \end{itemize}
We set $r:=|\Phi_1^{-1}(P)|$
for any $P\in (\nbigs\times X)\setminus\nbigd$.
We set $\nbigdtilde:=\Phi_1^{-1}(\nbigd)$.
For any $x\in\nbigs$,
we set $\nbigy_x:=p_1^{-1}(x)$,
$\nbigdtilde_x:=\nbigy_x\cap\nbigdtilde$
and
$\nbigd_x:=p_2^{-1}(x)\cap\nbigd$.
Let $g(X)$ denote the genus of $X$.
Let $\gtilde$ denote the genus of $\nbigy_x$,
which is independent of $x\in\nbigs$.

Let $\nbigl$ be a line bundle on $\nbigy$
such that
\[
\deg(\nbigl_{|\nbigy_x})=\gtilde-rg(X)+r-1.
\]
We obtain the locally free $\nbigo_{\nbigs\times X}$-module
$\nbige=\Phi_{1\ast}\nbigl$.
It is equipped with the relative Higgs field
\[
 \theta:
 \nbige\to
 \nbige\otimes \Omega^1_{\nbigs\times X/\nbigs}
\]
induced by the $\nbigo_{\nbigs\times T^{\ast}X}$-action
on $\Phi_{0\ast}\nbigl$.
For any $x\in \nbigs$,
let $(\nbige_x,\theta_x)$ be the induced Higgs bundle
on $X\simeq\{x\}\times X$.
We obtain the following lemma by the construction.
\begin{lem}
Each $(\nbige_x,\theta_x)$
is stable of degree $0$.
\hfill\qed
\end{lem}

\subsubsection{Statement}

We obtain the holomorphic line bundle
$\det(\nbige)$ on $\nbigs\times X$.  
There exists a $C^{\infty}$-Hermitian metric
$h_{\det(\nbige)}$ of $\det(\nbige)$
such that
$h_{\det(\nbige),x}:=h_{\det(\nbige)|\{x\}\times X}$
are flat for any $x\in \nbigs$.

We have the decomposable filtered Higgs bundle
$(\nbigp^{\star}_{\ast}\nbige_x,\theta_x)$
on $(X,\nbigd_{x})$.
Let $h_{\infty,x}$ be the decoupled harmonic metric
of $(\nbige_x,\theta_x)_{|X\setminus\nbigd_{x}}$
such that $\det(h_{\infty,x})=h_{\det(\nbige),x}$.

\begin{lem}
$h_{\infty,x}$ $(x\in\nbigs)$
induce a $C^{\infty}$-metric of
$\nbige_{|(\nbigs\times X)\setminus\nbigd}$.
\end{lem}
\pf
It is enough to study locally around any point $x_0\in\nbigs$.
By using examples in \S\ref{subsection;22.12.24.120},
we can construct 
a $C^{\infty}$-Hermitian metric $h_0$ of
$\nbigl_{|\nbigy\setminus\nbigdtilde}$
such that
(i) $h_0$ is flat around $\nbigdtilde$,
(ii) $h_{0|\nbigy_x\setminus\nbigdtilde_{x}}$ is adapted to
$\nbigp^{\star}_{\ast}(\nbigl_{|\nbigy_x})$.
By using Lemma \ref{lem;22.12.25.1} below,
we can construct a $C^{\infty}$-function $f$ on $\nbigy$
such that
$h_{1,x}=e^{f}h_{0,x}$ $(x\in\nbigs)$
is a family of flat metrics 
$\nbigl_{|\nbigy_x\setminus\nbigdtilde_{x}}$.
It induces a family of decoupled harmonic metrics $h_{2,x}$ of
$(\nbige_x,\theta_x)_{|X\setminus\nbigd_{x}}$
such that
they give a $C^{\infty}$-Hermitian metric $h_2$ of
$\nbige_{|(\nbigs\times X)\setminus\nbigd}$.
Note that
$\det(h_{2,x})$ induces
a flat metric of $\det(\nbige_x,\theta_x)$.
For each $x\in \nbigs$,
because both
$\det(h_{2,x})$ and $h_{\det(\nbige),x}$
of $\det(\nbige_x)$,
we obtain $\alpha_x>0$
determined by
$\det(h_{2,x})=\alpha_xh_{\det(\nbige),x}$.
Because $\det(h_{2,x})$ $(x\in\nbigs)$
give a $C^{\infty}$-metric of
$\det(\nbige)_{|(\nbigs\times X)\setminus\nbigd}$,
we obtain that $\alpha_x$ $(x\in\nbigs)$
give a $C^{\infty}$-function on $\nbigs$.
Because $h_{\infty,x}=e^{-\alpha_x/r}h_{2,x}$,
we obtain
$h_{\infty,x}$
induces a $C^{\infty}$-metric of
$\nbige_{(\nbigs\times X)\setminus\nbigd}$.
\hfill\qed

\vspace{.1in}
Let $h_{t,x}$ be a harmonic metric of $(\nbige_x,t\theta_x)$
such that $\det(h_{t,x})=h_{\det(\nbige),x}$.
Let $(V_x,\theta_x):=(\nbige_x,\theta_x)_{|X\setminus\nbigd_{x}}$.
We obtain the automorphism $s(h_{\infty,x},h_{t,x})$ of $V_x$
determined by
$h_{t,x}=h_{\infty,x}\cdot s(h_{\infty,x},h_{t,x})$.

\begin{thm}
\label{thm;22.12.25.3}
Let $x_0\in \nbigs$.
Let $K$ be any relatively compact open subset 
in $X\setminus\nbigd_{x_0}$.
Let $\nbigs_0$ be a neighbourhood of $x_0$
such that $\nbigs_0\times K$ is relatively compact
in $(\nbigs\times X)\setminus\nbigd$.
For any $\ell\in\seisuu_{\geq 0}$,
there exist positive constants
$C(\ell,K)$ and $\epsilon(\ell,K)$
such that the $L_{\ell}^2$-norm of
$s(h_{\infty,x},h_{t,x})-\id$ $(x\in\nbigs_0,t\geq 1)$
on $K$ with respect to $h_{\infty,x}$ and $g_X$
are dominated by
$C(\ell,K)\exp(-\epsilon(\ell,K)t)$.
\end{thm}

\subsubsection{Refined statement}

Let $x_0\in\nbigs$.
For any $P\in\nbigd_{x_0}$,
let $(U_P,z_P)$ be a simply connected holomorphic coordinate
neighbourhood of $P$ in $X$
such that $U_P\cap \nbigd_{x_0}=\{P\}$
and that $z_P$ induces $(U_P,P)\simeq (B(2),0)$.
Moreover, we assume that $z_P$ induces
a holomorphic isomorphism between neighbourhoods
of the closures of $U_P$ and $B(2)$.
Let $\nbigs_{1,P}$ be a relatively compact open neighbourhood
of $x_0$ in $\nbigs$
such that
\[
\nbigd\cap (\nbigs_{1,P}\times U_P)
\subset \nbigs_{1,P}\times\{|z_P|\leq (1/4)^{\rank E}\}.
\]
Let $\nbigs_1$ be a connected open neighbourhood of
$x_0$ in $\bigcap_{P\in\nbigd_{x_0}}\nbigs_{1,P}$.

For $P\in\nbigd_{x_0}$ and $x\in\nbigs_1$,
let $h_{t,P,x}$ be the harmonic metric of
$(\nbige_x,\theta_x)_{|\{|z_P|<1\}}$
such that
$h_{t,P,x|\{|z_P|=1\}}
=h_{\infty,x|\{|z_P|=1\}}$.
We note that Condition \ref{condition;23.1.28.30}
is satisfied for $h_{\infty,x|U_P}$
by Lemma \ref{lem;23.1.28.50},
and we can apply 
Proposition \ref{prop;23.1.28.40} to $h_{t,P,x}$.
We construct 
Hermitian metrics $\htilde_{t,x}$ 
of $\nbige_{x}$ $(x\in\nbigs_1)$
from $h_{\infty,x}$ and $h_{t,P,x}$ $(P\in\nbigd_{x_0})$
as in \S\ref{subsection;22.11.30.201}.
Let $s(\htilde_{t,x},h_{t,x})$
be the automorphism of $\nbige_x$
determined by
$h_{t,x}=\htilde_{t,x}\cdot s(\htilde_{t,x},h_{t,x})$.
By using Proposition \ref{prop;23.1.28.40},
we obtain the following theorem
in the same way as Theorem \ref{thm;22.12.24.101},
which implies Theorem \ref{thm;22.12.25.3}.

\begin{thm}
\label{thm;22.12.25.4}
For any $\ell\in\seisuu_{\geq 0}$,
there exist positive constants $C(\ell)$ and $\epsilon(\ell)$
such that the $L^{2}_{\ell}$-norms of
\[
 s(\htilde_{t,x},h_{t,x})-\id\quad(x\in\nbigs_1,t\geq 1)
\]
with respect to $\htilde_{t,x}$ and $g_X$
are dominated by
$C(\ell)\exp(-\epsilon(\ell)t)$.
\hfill\qed
\end{thm}

\subsubsection{Appendix}

Let $M$ be a compact oriented $C^{\infty}$-manifold.
Let $S$ be a $C^{\infty}$-manifold.
Let $g_{S\times M}$ be a Riemannian metric of $S\times M$.
For each $x\in S$,
we set $M_x:=\{x\}\times M$.
Let $g_x$ and $\Delta_x$
denote the induced Riemannian metric
and the associated Laplacian of $M_x$.

\begin{lem}  
\label{lem;22.12.25.1}
Let $f_1$ be a $C^{\infty}$-function on $S\times M$
such that $\int_{M_x} f_{1}\vol_{g_x}=0$.
Let $f_2$ be a function on $S\times M$ determined by
the conditions
$\Delta_x(f_{2|M_x})=f_{1|M_x}$
and $\int_{M_x}f_{2|M_x}\vol_{g_x}=0$.
Then, $f_2$ is $C^{\infty}$.
 \end{lem}
\pf
We explain only a sketch of a proof.
For any $x\in S$,
let $f_{i,x}:=f_{i|M_x}$.
Let $S_0$ be a relatively compact open subset in $S$.
There exists a uniform lower bound of the first non-zero eigenvalue
of the operators $\Delta_{x}$ $(x\in S_0)$
(see \cite[Theorem 5.7]{Li-Peter-Book}).
There exists $C_0>0$ such that $\|f_{1,x}\|_{L^2}\leq C_0$ $(x\in S_0)$.
By $\Delta_x(f_{1,x})=f_{2,x}$,
for any $\ell\in\seisuu_{\geq 0}$
there exists $C_1(\ell)>0$ such that
$\|f_{1,x}\|_{L^2_{\ell}}\leq C_1(\ell)$
for any $x\in S_0$.
Let $x(i)\in S_0$ be a sequence convergent to $x(\infty)\in S_0$.
There exists a subsequence $x'(j)$ convergent to $x(\infty)$
such that the sequence $f_{1,x'(j)}$ is weakly convergent in $L_{\ell}^2$
for any $\ell\in\seisuu_{\geq 0}$.
The limit $f_{\infty}$ satisfies
$\Delta(f_{\infty})=f_{2,x(\infty)}$
and $\int_{M_{x(\infty)}}f_{\infty}\vol_{g_{x(\infty)}}=0$.
We obtain $f_{\infty}=f_{1,x(\infty)}$.
Hence, $f_{1,x}$ and their derivatives in the $M$-direction
are continuous with respect to $x\in S$.

Let $S_1$ be a relatively compact open subset of $S$
equipped with a coordinate system $(x_1,\ldots,x_n)$.
Let $[\del_{j},\Delta_x]$ be the differential operator on $S_1\times M$
defined by 
$[\del_{j},\Delta_x](f)=\del_{j}(\Delta_x(f))-\Delta_x(\del_{j}f)$.
It does not contain a derivative in the $S_1$-direction.
Note that $[\del_{j},\Delta_x](f_{1,x})$ and their derivative
in the $M$-direction are continuous with respect to $x\in S_1$.
Let $f_{1,x}^{(j)}$ be the solution of
the conditions
$\Delta_{x}(f^{(j)}_{1,x})
=\del_{j}f_{2,x}-[\del_{j},\Delta_x]f_{2,x}$
and
$\int_{M_{x}}f^{(j)}_{1,x}\vol_{g_{x}}=0$.
Choose $y=(y_1,\ldots,y_n)\in S_1$.
We define functions $F^{(j)}_x$ on $M_x$ by
$F^{(j)}_{x}=(x_j-y_j)^{-1}(f_{1,x}-f_{1,y})$
if $x_j\neq y_j$,
and
$F^{(j)}_x=f^{(j)}_{1,x}$
if $x_j=y_j$.
It satisfies
$\Delta_x(F_x^{(j)})
=(x_j-y_j)^{-1}(f_{2,x}-f_{2,y}-(\Delta_x-\Delta_y)f_{1,y})$
if $x_j\neq y_j$,
and
$\Delta_x(F^{(j)})_x=\del_jf_{2,x}-[\del_j,\Delta_x]f_{2,x}$
if $x_j=y_j$.
Then, by an argument in the previous paragraph,
we can prove that $F^{(j)}_x$ and their derivatives in the $M$-direction
are continuous with respect to $x$.
It implies that $f_{1,x}$ is $C^1$-with respect to $x$
and that $\del_jf_{1,x}=f^{(j)}_{1,x}$.
By a similar argument,
we can prove that $f_{1,x}$ and their derivatives in the $M$-direction
are $C^{\infty}$ with respect to $x$.
\hfill\qed

\end{document}